\documentclass[12pt]{amsart}
\usepackage{amscd,amsmath,amssymb,amsfonts}

\theoremstyle{plain}

\theoremstyle{definition}

\theoremstyle{remark}

\newcommand{\limto}{{\displaystyle\lim_{\longrightarrow}}}

\newcommand{\rightlim}{\mathop{\limto}}
\newcommand{\limfromn}{\leftlim\limits_{\raise3pt\hbox{$n$}}}

\newcommand{\leftlim}{\mathop{\displaystyle\lim_{\longleftarrow}}}
\newcommand{\rightlimit}[1]{\mathop{\lim\limits_{\longrightarrow}}\limits%
                    _{\raise3pt\hbox{$\scriptstyle #1$}}}
\newcommand{\leftlimit}[1]{\mathop{\lim\limits_{\longleftarrow}}\limits%
                    _{\raise3pt\hbox{$\scriptstyle #1$}}}
\newcommand{\iso}{\buildrel\sim\over\rightarrow}
\newcommand{\hra}{\hookrightarrow}

\newcommand{\Coker}{\text{Coker}}
\newcommand{\Ker}{\text{Ker}}

\newcommand{\fd}{{\frak d}}

\newcommand{\ff}{{\frak f}}

\newcommand{\fFr}{{\frak F\frak r}}

\newcommand{\fgl}{{\frak{gl}}}

\newcommand{\fr}{{\frak r}}
\newcommand{\fs}{{\frak s}}

\newcommand{\Ad}{{\text{Ad}}}

\newcommand{\Aut}{\text{Aut}}

\numberwithin{thm}{subsection}
\numberwithin{equation}{subsection}

\newcommand{\eq}[2]{\begin{equation}\label{#1}#2 \end{equation}}

\newcommand{\gr}{{\rm gr}}

\newcommand{\red}{{\rm red}}

\newcommand{\Hom}{{\rm Hom}}

\newcommand{\Spec}{{\rm Spec \,}}

\newcommand{\Tr}{{\rm Tr}}
\newcommand{\GL}{{\rm GL}}
\newcommand{\rO}{{\rm O}}

\newcommand{\Res}{{\rm Res }}

\newcommand{\sA}{{\mathcal A}}

\newcommand{\sC}{{\mathcal C}}
\newcommand{\sD}{{\mathcal D}}
\newcommand{\sE}{{\mathcal E}}

\newcommand{\sG}{{\mathcal G}}
\newcommand{\sH}{{\mathcal H}}

\newcommand{\sL}{{\mathcal L}}
\newcommand{\sM}{{\mathcal M}}

\newcommand{\sO}{{\mathcal O}}
\newcommand{\sP}{{\mathcal P}}
\newcommand{\sQ}{{\mathcal Q}}

\newcommand{\sS}{{\mathcal S}}

\newcommand{\sV}{{\mathcal V}}

\newcommand{\A}{{\mathbb A}}

\newcommand{\G}{{\mathbb G}}

\newcommand{\Q}{{\mathbb Q}}

\newcommand{\W}{{\mathbb W}}

\newcommand{\Z}{{\mathbb Z}}

\begin{document}

\title[$\sE$-factors for Gau\ss-Manin
Determinants]{$\sE$-factors
for Gau\ss-Manin Determinants } 
\author{Alexander Beilinson}
\address{Dept. of Mathematics,
University of Chicago,
Chicago, IL 60637,
USA}
\email{sasha@math.uchicago.edu}
\author{Spencer Bloch}
\address{Dept. of Mathematics,
University of Chicago,
Chicago, IL 60637,
USA}
\email{bloch@math.uchicago.edu}
\author{H\'el\`ene Esnault}
\address{
Universit\"at Essen, FB6, Mathematik, 45117 Essen, Germany}
\email{esnault@uni-essen.de}
\date{May 10, 2002}
\dedicatory{To Yuri Manin, with gratitude and
admiration.}
\begin{abstract}
We define $\varepsilon$-factors in the de Rham
setting and  calculate the determinant
of the Gau\ss-Manin connection for a family of
(affine) curves and a vector bundle equipped
with a flat connection.
\end{abstract}
\maketitle
\begin{quote}
``Ordentliche Leute pflegten ihren Schatten mit sich zu nehmen,
wenn sie in die Sonne gingen.'' \newline
A. v.
Chamisso, {\it Peter
Schlemihls wundersame Geschichte}
\vspace*{.5cm}
\noindent
\end{quote}

\medskip

\centerline{Table of contents}
\medskip

\hskip 2 cm 1. Introduction

\hskip 2 cm  2. Determinant lines of Fredholm
operators

\hskip 2 cm 3. The Heisenberg group and its
cousins

\hskip 2 cm 4. The $\varepsilon$-factors

\hskip 2 cm 5. Some formula (with apologies to
Charlotte)

\section{Introduction}

\subsection{}  Let us first recall briefly the
format of the classical theory of
$\varepsilon$-factors, which grew out of Tate's
thesis and the work of Dwork, Langlands,
Deligne, and Laumon.

We consider only the case of function
fields, so we fix a prime
$p$, a finite field $k$ of characteristic $p$,
and look at classical local
$k$-fields, i.e., topological $k$-fields
isomorphic to
$ k' ((t))$ where $k'$ is a finite extension
of $k$. For a local field $F$ we denote by
$\omega (F)$ the 1-dimensional $F$-vector
space of differentials; set
$\omega (F)^\times :=
\omega (F)\smallsetminus \{ 0\}$. Galois
modules and local systems have
coefficients in $\bar{\Q}_\ell$ for a
prime $\ell
\neq p$.

\medskip

A {\it  classical theory of
$\varepsilon$-factors} is a rule which
assigns to every pair $(F,V )$, $F$ is a
local $k$-field, $V$  a
Galois
module for $F$, a continuous function
$\varepsilon (F,V )=\varepsilon (V ):
\omega (F)^\times   \to \bar{\Q}_\ell^\times$,
$\nu
\mapsto
\varepsilon (V )_\nu$, which satisfies the
following properties:

(i) {\it Multiplicativity with respect to
$V$:} For a short exact sequence $0\to
V_1 \to V \to V_2 \to 0$ one has
\eq{1.1.1}{\varepsilon (F,V )_\nu
=\varepsilon (F,V_1 )_\nu \varepsilon
(F,V_2 )_\nu .}
Thus $\varepsilon (F,V)_\nu$ makes sense for
virtual Galois modules.

(ii) {\it Induction:} Let $F'/F$ be a finite
separable extension, $V'$ a
  Galois module for $F'$, $\nu
\in
\omega (F)^\times \subset \omega (F' )^\times$.
Then
\eq{1.1.2}{\varepsilon (F', V' )_\nu
=\varepsilon (F ,\text{Ind}V' )_\nu }
if $V'$ is a virtual Galois module
of rank 0.

(iii) {\it Product formula:} Let $X$ be a
  smooth projective curve over $k$,
$D\subset X$ a divisor,
$V$ a
local system  on
$U:=X\smallsetminus D$, $\nu$ an invertible
1-form on
$U$. For $x \in D$ let $F_{x }$
be the local field at $x$, and $V_x$,
$\nu_x$  the restrictions of
$V$,
$\nu$ to $F_x$.  Then\footnote{Here $Fr_k$ is the geometric
Frobenius for
$k$, $\bar{k}$ is an algebraic closure of
$k$.}
\eq{1.1.3}{
\mathop\prod\limits_{x \in D}\varepsilon
(F_{x} ,V_{x}
)_{\nu_x}    =  \prod_i \det (-Fr_k , H^{i}
(U_{\bar{k}}, V ))^{(-1)^{i+1}}.}

It follows from a deep theorem of Laumon
\cite{L} 3.2 that a theory of
$\varepsilon$-factors does exist. Namely,
for fixed non-trivial character $\psi : k\to
\bar{\Q}_\ell^\times$ and  $p^{1/2}
\in \bar{\Q}_\ell^\times$ the function
$\varepsilon (F,V)_\nu := q_F^{-rk V \cdot
v_F (\nu ) /2}\varepsilon_\psi (O_F ,Rj_* V,\nu
)$ is a theory of $\varepsilon$-factors in the
above sense. Here $q_F$ is the number of
elements of the residue field $k_F$, $v_F (\nu
)$ the valuation of $\nu$, $O_F$ the ring of
integers in $F$, and $\varepsilon_\psi$ is as
in \cite{L} (3.1.5.4).\footnote{In 
notation of \cite{D2} one has
$\varepsilon (F,V)_\nu = \epsilon (V,\psi_\nu
, \mu_\nu ,1) \det (-Fr_{k_F}, V^I (-1) )$
where
$\psi_\nu : F\to \bar{\Q}_\ell$ is an
additive character $f\mapsto \psi
\Tr\Res (f\nu )$, $\mu_\nu$  a Haar measure
on $F$ self-dual with respect to $\psi_\nu$,
and $I\subset Gal (\bar{F}/F)$ the inertia
subgroup.}

\medskip

{\it Remarks.} (i) Assume that our $(F,V,\nu )$
is such that $V$ is unramified and $v_F
(\nu )=0$.
Then\footnote{To see this compare (1.1.3) for
$U$ and
$U$ with one point removed.} (here  $(-1)$ is
the Tate twist)
\eq{1.1.4}{\varepsilon (F,V)_\nu =\det
(-Fr_{k_F}, V(-1)).}

(ii)  To determine
$\varepsilon $ one  invokes  the Gabber-Katz
theorem
\cite{GK} according to which every Galois
module for $k((t))$ can be extended to a
local system on $\G_m =\Spec k[t,t^{-1}]$
having tame ramification at
$t=\infty$. Therefore, by the product formula,
we know $\varepsilon (V )$ for arbitrary $V$
  if the global Frobenius determinants and
$\varepsilon (V )$ for tame $V$'s  happen
to be known.

(iii) The r.h.s.~in (1.1.3) is the (super)
trace of the Frobenius acting on the
determinant super line of the complex
$R\Gamma (U_{\bar{k}},V)[1]$.

(iv)
For a theory
$\varepsilon$ of $\varepsilon$-factors and
a fixed
$a\in k^\times$ the function $\nu \mapsto
\varepsilon (V )_{a\nu}$ is again a theory
of $\varepsilon$-factors. So the set $\mathbb E
(k)$ of theories of $\varepsilon$-factors
carries a
$k^\times$-action.

\subsection{} Unfortunately, the above story,
being global-to-local, does not tell much about
the nature of
$\varepsilon$-factors.  Therefore one is
tempted to look for a direct, purely local,
construction of
$\varepsilon$-factors. If available, such
construction would make it possible to use the
product formula  to compute  the global
Frobenius determinant.

A related dream is to have a {\it geometric}
theory of $\varepsilon$-factors. One considers
local fields over an arbitrary (not necessary
finite) base field $k$, and we look for a rule
$\sE$ that assigns to every
$(F,V )$ as above an invertible (super) local
system $\sE (V )$ on the
$k$-(ind-)scheme $\omega (F)^\times$. The
compatibilities 1.1(i)--(iii) now mean
canonical isomorphisms of invertible super
local systems. The classical
$\varepsilon$ is recovered as the ``trace of
Frobenius" function of $\sE$.\footnote{To be
more precise, $\sE (V)$ should rather be a
twisted local  system,
  an object of an appropriate
``$\varepsilon$-gerbe". A choice of
trivialization of the $\varepsilon$-gerbe
identifies then $\sE (V)$
  with a plain local system. One can
further speculate that in the gerbe setting
compatibility (1.1.2) holds for arbitrary $V'$,
and the reason behind the rank 0 condition
in (1.1.2) is the absence of a nice
trivialization of the $\varepsilon$-gerbe.
We are grateful to M.~Kapranov for this
suggestion.}

\subsection{} The aim of this
article is to present such geometric theory in
the de Rham setting. So our $k$ is now  a
field of characteristic 0,  and instead of
\'etale local systems we consider local
systems in the de Rham sense, that is vector
bundles equipped with flat connections. Recall
that there is a well-known, yet mysterious,
analogy between the phenomena of wild
ramification (\'etale setting) and irregularity
(de Rham setting). The role of local Galois
modules is played by de Rham local systems $V=
(V,\nabla )$ on a {\it formal} punctured
disc\footnote{Notice that the Stokes
structure, which is necessary to determine the
{\it analytic} local structure of a
connection, does not play any role here.}
$\Spec F$, where $F$ is a $k$-algebra
isomorphic to $k' ((t))$, $k'$ is a finite
extension of $k$.

So for every such
$V$  we define a super line bundle (its
degree is the irregularity of
$\nabla$) with flat connection
$\sE (V)$ on
$\omega (F)^\times$ together with data of
canonical isomorphisms parallel to
1.1 (i)--(iii). For
example,  the product formula  looks as
follows. Consider
$X, U,D, V$ as in 1.1(iv), so $V$ now is a
local system on $U$.
Let $\sE (F_D ,V)$ be the exterior
tensor product of  $\sE (F_x ,V_x
)$ for
$x\in D$. Then one has a
canonical isomorphism of  local systems on
$\omega (U)^\times$
\eq{1.3.1}{\sE (F_D
,V)|_{\omega
(U)^\times}\iso \mathop\otimes\limits_i
(\text{det}H^{i}_{dR}
(U,V))^{\otimes (-1)^{i+1}}_{\omega
(U)^\times}.} Here the l.h.s.~is the
restriction of $\sE (F_D ,V)$ to
$\omega (U)^\times \hra \prod \omega (F_x
)^\times$,  the r.h.s.~is a constant local
system.

{\it Remark.} Our $\varepsilon$-factors
differ from de Rham counterparts of
$\varepsilon$-factors from
\cite{L}. For example:

(i) For our de Rham
$\varepsilon$-factors compatibility (1.1.2)
holds for arbitrary $V'$ (regardless
of its rank).\footnote{ So in the de Rham
setting the hoped-for
$\varepsilon$-gerbe from ftn.~4 is
canonically trivialized.}

(ii)  If $V$ is a trivial local system of
rank 1 on $\Spec F$ then the local system $\sE
(V)$ has a non-trivial $\pm 1$ monodromy  on
the connected components of $\omega
(F)^\times$ which consist of $\nu$ with odd
order of zero.\footnote{In accordance with the
fermionic nature of $\sE$.}

\medskip

If our datum varies
nicely\footnote{Every variation is nice at
the generic point; for more details see  4.4.}
with respect to parameters
$S$ (this means, in particular, that
irregularity does not jump) then
$\sE (F_s ,V_s )$ form a super line bundle
$\sE (F/S ,V)$ on
$\omega (F/S)^\times$ equipped with an
$S$-relative connection. If the vector
bundle $V$ on the space of the family carries
an absolute flat connection  extending the
  $S$-relative one, then our
$S$-relative connection on $\sE (F/S,V)$
  extends canonically to an
absolute flat connection, i.e., $\sE (F/S)$
becomes a local system. The product formula
isomorphism (1.3.1) is compatible with the
absolute connections (the l.h.s.~carries the
Gau\ss-Manin one). So  (1.3.1) yields a
formula for the Gau\ss-Manin determinant, once
we are able to compute the
$\varepsilon$-connections.

\subsection{}
The $\varepsilon$-line $\sE (V)_\nu$ is defined
as follows.  Let
$\tau_\nu$ be the vector field
$\nu^{-1}$. Then $\nabla (\tau_\nu
)$ is a Fredholm operator on the
  infinite-dimensional
$k$-vector space $V$.\footnote{Its
index vanishes and the determinant line is
canonically trivialized, see e.g.~5.9(a)(iv).}
Let
$V=V^+
\oplus V^-$ be any decomposition of the $V$
into a sum of Fourier ``positive" and
``negative" parts.\footnote{Precisely, for
$F=k((t))$ this means that for certain
$F$-basis $\{ e_i \}$ of $V$ our $V^+$
contains $\Sigma k[[t]]e_i$ as a $k$-subspace
of finite codimension.} Consider the
component of $\nabla (\tau_\nu
)$ acting on
$V_-$. Our
$\sE (V)_\nu$ is the determinant line of this
Fredholm operator. The auxiliary choice of
$V^\pm$ is irrelevant.\footnote{The
determinant line of a Fredholm operator
$F:U\to U$ does not change if $U$
is modified by a finite-dimensional vector
space, or the operator is perturbed by an
operator with finite-dimensional image.}

To see the product formula isomorphism
notice that the composition $\Gamma (U,V)\to
\oplus V_x
\to \oplus V_x^-$ is Fredholm. Now (1.3.1)
is the corresponding identification of the
determinant lines for $\tau_\nu$  on
$\Gamma (U,V)$ and $\oplus V_x^-$.

The $\varepsilon$-connection on $\sE (V)$
is defined as follows. First, the
multiplicativity of determinant lines with
respect to product of operators shows that the
action of $F^\times$ on $\omega (F)^\times$
lifts to an action on
$\sE (V)$ of a Heisenberg group controlled by
$\det (V)$. The connection on $\det V$
identifies the corresponding Heisenberg Lie
algebra with the standard one (for trivialized
$V$). Finally, every $\nu \in \omega
(F)^\times$ defines on $F$ a non-degenerate
scalar product which yields a ``self-adjoint"
splitting of the standard Heisenberg Lie
algebra; this splitting specifies the
horizontal directions for the
$\varepsilon$-connection at $\nu$. In
fact,
for all directions but one (namely, along
the fibers of the map
$\nu
\mapsto
\Res
\,\nu$) the $\varepsilon$-connection comes from
the action of the group of infinitesimal
automorphisms of $F$ on our picture. This
provides the absolute $\varepsilon$-connection
refered to in 1.3.

The above argument for the product formula
is essentially Tate's
proof \cite{T2} of the residue formula.
Replacing the first
order differential operator $\nabla (\tau_\nu
)$ by a zero order one you get a proof
of Weil's reciprocity (cf.~\cite{ACK}); its
infinitesimal version is the residue formula.

\subsection{} The idea that
$\varepsilon$-factors have to do
with fermionic determinants is due
to Laumon:  the introduction to
\cite{L} opens with a discussion of
Witten's picture of classical Morse theory  as
a WKB approximation to  quantum mechanics of a
super particle. This passage was
considered by some as {\it
nezabudki}\footnote{Forget-me-nots from a
fable {\it ``Nezabudki i Zapyatki"} of  Koz'ma
Proutkoff$^{13}$ beloved by
I.~M.~Gelfand.}$^,$\footnote{A famous writer
and philosopher (1803--1863) whose project
``On introducing unanimity in Russia"
determined the  Russian Way for the XX--XXI
centuries.} for Laumon's construction itself
stems from the mere fact that the Fourier
transform of a compactly supported
distribution $f$ on a line is smooth, and its
asymptotics at
$\infty$ are controlled by the singularities of
$f$.

\subsection{}
The formula for the $\varepsilon$-connection
looks as follows. We assume that $S=\Spec K$,
$K$ is a field. The group of isomorphism
classes of line bundles with connection on
$S$ equals $\Omega^1_{K/k}/d
\text{log}(K^\times )$, so our formula
should determine an element in this group.

An absolute connection $(V,\nabla)$ on $K((t))$
is said to be {\it admissible} if it admits a
lattice $\sV_O = \oplus K[[t]] \subset V$ with
respect to which the connection has the form
(for some $m$)
$ \nabla = d + \sA(t) = d + g(t)dt/t^m
+\eta /t^{m-1}$,
$g
\in \GL_n (K[[t]])$, $\eta\in
\text{Mat}_n(\Omega^1_K)\otimes K[[t]]$.
When $(V,\nabla)$ is admissible  and
irregular,  we find\footnote{See 5.6 for a
more precise statement for general $\nu$.} that
the class of the connection on
$\sE_{t^{-m}dt}$ equals
\eq{1.6.1}{\text{Res}_t
\Tr(g^{-1}dg\wedge\eta/t^{m-1}) -
\frac{m}{2}d\log\det(g(0) ). }

The corresponding formula for the global
Gau\ss-Manin determinant was found (in case of
admissible singularities and a genus 0 curve)
in
\cite{BE3}; it was the starting point of the
present article.

The Levelt-Turritin theorem helps to reduce
(in non-effective way) computation of the
$\varepsilon$-connection for arbitrary
$\nabla$
to the  rank 1 situation (which is
always admissible).

\subsection{} A list of natural questions:

(a) How the above
picture is related to the rank 1 story of
\cite{BE2}?

(b) Can $\varepsilon$-factors be seen
microlocally? What is the relation between
the de Rham counterpart of Laumon's Fourier
approach \cite{L} and our picture?

(c) What could be a
higher dimensional generalization?

(d) Is there a de Rham version of Tate's
Thesis \cite{T1}? More generally, what about
the automorphic counterpart of the whole story?

(e) Can one treat families of
connections with jumping irregularity?

\subsection{} The article is organized as
follows.  Section 2 reviews the basic
formalism of determinant lines for families of
Fredholm operators in the algebraic setting.
The ideas here go back to \cite{KM} and
\cite{T2}. Much of the  subject is considered
from the analytic point of view in the book
\cite{PS} (or in the earlier articles
\cite{DKJM}, \cite{SW}), but we were not able
to find a convenient reference for the
algebraic situation.
Our exposition (which makes no pretense to
originality) is based on localization
statements 2.3(ii) and 2.12 borrowed from
\cite{Dr}.\footnote{A typical object we meet
is a projective $R((t))$-module of finite
rank, and 2.12 describes it as a
topological $R$-module  \'etale locally
on Spec$\, R$. } We tried to make signs take
care of themselves (cf.~\cite{ACK}) using the
language of ``super
extensions";\footnote{Borrowed from
\cite{BD}.} the appendix to sect.~2 (which is
a variation on theme of SGA 4 XVII 1.4) 
stores the required general nonsense (for a
thorough discussion of super subject we refer
the reader to \cite{BDM}). Section 3 begins
with a review of the Heisenberg (super)
extension of the group ind-scheme
$k((t))^\times$.  Its commutator pairing is
the parametric version  of the tame symbol from
\cite{C}. According to \cite{T2} or
\cite{ACK}, this format yields automatically
the Weil reciprocity. The Heisenberg group
approach to the 
Contou-Carr\`ere  symbol is a group-theoretic 
version of Tate's construction of residue
\cite{T2}; it was known to specialists for
quite a time (A.B.~learned it from P.~Deligne
about 10 years ago), but, apart from the
case of tame symbol proper, seems not to be
documented.  The Heisenberg action controls
the dependence of  our
$\varepsilon$-factors on $\nu$ and is
  responsible for the existence of the
$\varepsilon$-connection. The key fact here
(see 3.10) is that on the
$k((t))^\times$-torsor of invertible forms
$k((t))^\times dt$ every bundle equivariant
with respect to the action of the Heisenberg
group twisted by a local system on $\Spec
k((t))$  acquires automatically a plain flat
connection. In the fourth section we define
$\varepsilon$-factors and study their basic
properties. The final
section deals with explicit calculations of
$\varepsilon$-connections.  We explain  how
$\sE_\nu$ depends on $\nu$ in the general
situation. Then we treat the case of regular
singularities (cf.~\cite{BE1}), the one of
irregular admissible connections
(cf.~\cite{BE3}), and show that the
general computation can be reduced, in
principle, to the rank 1 case.

\medskip

{\it Acknowledgements:} It is a pleasure to
thank Takeshi Saito for discussions on the
analogy of the formula as written in
\cite{BE3} with the existing formulae for the
$\varepsilon$-factors of $L$-functions.
We also thank Gerd
Faltings for explaining to us his own version of a variation of
the global method of \cite{BE3}. We are most
greatful to Vladimir Drinfeld for
advice and help in understanding the
subject. We want to thank the referee for a
number of corrections and illuminating
remarks.  

Since this paper was written, we were informed
that a construction of $\varepsilon$-lines
and product formula isomorphism (1.3.1)
(but not the $\varepsilon$-connection) which
coincides essentially with 1.4, was given by
P.~Deligne in his unpublished seminar at IHES
in May-June 1984.

The authors were partially supported by,
respectively,  NSF grants DMS-0100108 and
DMS-0103765, and  DFG-Schwerpunkt ``komplexe
Mannigfaltigkeiten".

\bigskip

\section{Determinant lines of Fredholm
operators}

In this section we establish the basic structure of Fredholm determinants,
working insofar as possible in the category of functors on algebras. Of
particular importance is the material on Clifford algebras in
(2.14)--(2.18): it will be used in
section 3 to define a $\mu_2$-structure on a
certain determinant line over 
$\omega (F)^\times$ which, in turn, leads to a
connection on the epsilon line. 

  We write  ``$P\in\sC$" for
``$P$ is an object of a category $\sC$."

\subsection{Super lines and super extensions}
Let
$R$ be a commutative ring. We denote by
$\sM^s_R$ the tensor category of {\it super
$R$-modules.} So an object of $\sM^s_R$ is a
$\Z/2$-graded $R$-module $M=M^{\bar{0}}\oplus
M^{\bar{1}}$. The tensor product as well as 
the associativity
constraint, is
the usual tensor product of $\Z/2$-graded
modules,  and the commutativity constraint is
$a\otimes b = (-1)^{p(a)p(b)}b\otimes a$ where
$p(a)\in \Z/2$ is the degree (or parity) of
$a$, i.e., $p(a)=i$ for $a\in
M^i$. See e.g.~ch.1 of \cite{BDM} for
details.

  A {\it super
$R$-line} (or simply {\it super line}) is an
invertible object of $\sM^s_R$. Super lines
form a Picard
groupoid\footnote{See appendix to this
section, A1.}
$\sP ic_R^s$. Explicitly, a super line is a
pair
$(\sL ,p )$ where $\sL$ is an
invertible $R$-module and
$p$ (the parity of our super line)
is a locally constant function
$\Spec R \to \Z/2$.
  In notation we usually abbreviate $(\sL ,p)$
to  $\sL$ and write $p=p(\sL )$; we also
write
$\sL\cdot
\sL' :=\sL\otimes\sL'$,
$\sL /\sL' :=
\sL \otimes \sL^{' -1}$. The unit object
$(R,0)$ of $\sP ic^s_R$ is denoted by $1_R$.

  The group of isomorphism classes
of super lines
$\pi_0 (\sP ic_R^s )$ is $ Pic (R)\times\Z/2_R$
where $\Z/2_R$ is the group of locally
constant functions $p:\Spec R \to \Z/2$.  One has $\pi_1 (\sP ic_R^s ):=\Aut 1_R
=R^\times$.

The above picture is functorial with
respect to morphisms of rings: every $f:R\to
R'$ yields a base change morphism of tensor
categories $f^* :\sM_R^s \to \sM^s_{R'}$,
$M_R \mapsto M_{R'}:= M_R \mathop\otimes\limits_R
R'$, hence a morphism of Picard groupoids
$f^* :\sP ic^s_R \to \sP ic^s_{R'}$.
\medskip

{\it Variants:} (a) Replacing line bundles in
the above definition by $\mu_2$-torsors on
$\Spec R$ we get the Picard groupoid
$\mu_2$-{\it tors}$^s_R$  of {\it super
$\mu_2$-torsors}. So a super
$\mu_2$-torsor is the same as a super line
$\sL$ equipped with an isomorphism $\sigma :
\sL^{\otimes 2}\iso 1_R$ (we
identify $(\sL, \sigma )$ with the
$\mu_2$-torsor of sections
$\ell$ of $\sL$ that satisfy $\sigma
(\ell\otimes \ell )=1$, its parity is the
parity of
$\sL$).\footnote{The tensor product in this
format is
$(\sL ,\sigma )\otimes (\sL' ,\sigma'
)=(\sL \otimes \sL' ,\sigma \cdot \sigma'
)$ where $(\sigma \cdot \sigma' )((\ell_1
\otimes \ell'_1 )\otimes (\ell_2
\otimes\ell'_2 )):= \sigma (\ell_1
\otimes\ell_2 )\sigma' (\ell'_1 \otimes\ell'_2
)$. }

(b) Replacing $\Z/2$ in the above definition
by $\Z$ we get the Picard groupoid $\sP
ic^\Z_R$ of {\it
$\Z$-graded super lines}. So a $\Z$-graded
super line is the same as a pair $(\sL ,v )$
where $\sL$ is a super line and $v: \Spec R\to
\Z$ is a locally constant function such that
$p(\sL )=v\mod 2$. As above, we usually
abbreviate
$(\sL ,v)$ to
$\sL$; we call
$v=v(\sL )$ the {\it degree} or {\it index} of
the $\Z$-graded super line.

\medskip

  We will consider $\sP
ic^s_R$-extensions (see A2) of
various groupoids $\Gamma$ and refer to them
as {\it super $\sO^\times$-extensions} or
simply {\it super extensions.} Every super
extension
$\Gamma^\flat$ yields a
homomorphism
$\Gamma
\to \pi_0 (\sP ic_R^s )$ (see A2, Remark (ii)),
hence we have the parity homomorphism
$p :\Gamma \to\Z/2_R$.

Notice that $\sP ic^s_R$ coincides as a {\it
monoidal} category\footnote{i.e., we forget
about the commutativity constraint.} with the
product of $\sP ic_R$
and the discrete  groupoid $\Z/2_R$. By A4,
a $\sP ic_R$-extension of $\Gamma$ is the same as a
sheaf of central extensions of $\Gamma$
by $\sO^\times$ on $\Spec R$; we refer to
it as
a {\it plain
$\sO^\times$-extension}. Thus\footnote{Use A3
Remark (ii) and A2 Remark (ii).}  a super
$\sO^\times$-extension amounts to a plain
$\sO^\times$-extension of
$\Gamma$ together with a homomorphism $p
:\Gamma \to
\Z/2_R$.

If $\Gamma$ is a group and we have
its super extension $\Gamma^\flat$, then every
pair of commuting elements
$\gamma ,\gamma' \in \Gamma$ yields $\{
\gamma ,\gamma' \}^\flat
\in R^\times$ (see A5). Looking at the
corresponding plain extension we get $\{
\gamma ,\gamma' \}^{plain}\in R^\times$. One
has
\eq{2.1.1}{ \{\gamma,\gamma'\}^\flat =
(-1)^{p(\gamma)p(\gamma')}\{ \gamma ,\gamma'
\}^{plain}.}

  We leave it to the reader to
define the notion of  {\it super
$\sO^\times$-extension} of a group valued
functor on the category of commutative
$R$-algebras.

Below we refer to $\sP ic^s_R$-torsors (see
A6) as (neutral) {\it super
$\sO^\times$-gerbes} or simply {\it super
gerbes} on
$\Spec R$.\footnote{We skip the word
neutral since the subject of
  this article is
$R$-local, so we can assume that all gerbes
are neutral.}  As always, a {\it
trivialization} of a super gerbe is its
identification with
$\sP ic_R^s$. We also have the notion of {\it
super pre-gerbe} = pre $\sP ic^s_R$-torsor (see
Remark in A6).

Replacing in $\sP ic^s$ by $\sP ic^\Z$ or
$\mu_2${\it -tors}$^s$ we get the notions of
$\Z$-graded super $\sO^\times$-extension, resp.
super
$\mu_2$-extension. Same for $\Z$-graded super
$\sO^\times$-gerbes, super $\mu_2$-gerbes,
etc.

\subsection{Determinant lines } We denote by
$\sV_R$  the category of projective
$R$-modules; $\sV^\ff_R \subset
\sV_R$ is the subcategory of  modules of finite
rank.

For $M\in \sV^\ff_R$ let $\det M =\det_R
M\in\sP ic^\Z_R$ be the top exterior power of
$M$ placed in degree rk$M$. This
$\Z$-graded super line is functorial with
respect to isomorphisms of
$M$'s, so for
$f :N\iso M$ we have $\det f :\det N \iso \det
M$ or $\det f : 1_R \iso \det M /\det N$.

\medskip

The following standard compatibilities hold:

(i) For every finite
family
$\{ M_\alpha
\}$ there is a canonical
isomorphism
\eq{2.2.1}{\text{det}(\oplus M_\alpha
)=\otimes \text{det}M_\alpha .}

Recall that it is this
compatibility that forces us to consider $\det
M$ as a {\it super} line.

(ii) For $M$ equipped with a finite filtration
with  projective subquotients,
there is a canonical isomorphism
\eq{2.2.2}{\det M=\det (\text{gr} M)  .}

(iii)
There is a canonical isomorphism of super lines
\eq{2.2.3}{\det (M^* ) =(\det M )^{-1},} where
$M^* := \Hom (M,R)$ is the dual module, defined
by a pairing $e_M :\det M \cdot \det M^* \to R$
such that for a base\footnote{We consider our
picture locally on $\Spec R$.}
$m_1 ,.., m_n$ of
$M$ and the dual base $m_1^* ,.., m_n^*$ of
$M^*$ one has
$e ((m_1
\wedge ..\wedge m_n )\cdot (m_n^* \wedge
..\wedge m_1^* ))=1$.

Notice that (2.2.3) changes the $\Z$-grading
to the  opposite one, and one has
$e_{M^*}= e_M^m$ (see A1 for the notation).

\medskip

The following two useful facts about
projective modules of infinite rank are due,
respectively, to Kaplansky and
Drinfeld:

\subsection{Proposition}
(i)  If $R$ is a local ring then
every projective $R$-module is free.

(ii)  Let $M$ be a projective
$R$-module and
$K\subset M$ a finitely generated submodule.
Then Zariski locally on $\Spec R$ one can find
a finitely generated submodule $P\subset M$
which contains $K$  such that $M/P$ is
projective.

\begin{proof} (i) See \cite{K}.

(ii) See \cite{Dr} 4.2; we reproduce the
proof for completeness sake. Take any
$x\in\Spec R$; let $R_{(x)}$ be the
corresponding local ring. Then $M_{(x)}:=
R_{(x)}\otimes M$ is a free
$R_{(x)}$-module by (i).  So there is an
embedding $i_{(x)} : R^n_{(x)}
\hra M_{(x)}$ with projective cokernel whose
image contains  $K_{(x)}\subset M_{(x)}$.
Replacing
$\Spec R$ by an open affine neighbourhood of
$x$ we can assume that $i_{(x)}$ comes from
$i: R^n \to M$ whose image contains $K$, and
then that
$i$ is an embedding with a projective
cokernel (to see this represent $M$ as a
direct summand of a free module to reduce
the statement first to the case when $M$
is free, then to the case when $M$ is free of
finite rank).
\end{proof}

\medskip

{\it Remark.} If $R$ is a Noetherian
ring such that $\Spec R$ is connected
then, according to
\cite{Ba}, every projective $R$-module of
infinite rank is free. We will not use this
fact.

\subsection{Asymptotic morphisms and Fredholm
morphisms} For projective modules $M,N\in\sV_R$
let
$\text{Hom}^{\ff}_R
(N,M)\subset\text{Hom}_R (N,M)$ be the
$R$-submodule of morphisms  whose image lies in
a finitely generated $R$-submodule of
$M$. Set \eq{2.4.1}{\Hom^\infty_R (N,M):=\Hom_R
(N,M)/\Hom^\ff_R (N,M).} Elements of
$\Hom^\infty$ are  {\it asymptotic
morphisms}; for  $f:N\to M$ the corresponding
asymptotic morphism is denoted by $f^\infty$. The
composition of asymptotic morphisms is
well-defined.\footnote{The composition
of a morphism from
$\Hom^\ff$ with any morphism lies in
$\Hom^\ff$.} So we have an $R$-category
$\sV^\infty_R$ of projective $R$-modules
and
asymptotic morphisms together with the obvious
functor $\sV_R\to {\sV}^\infty_R$
which is the identity on objects. Notice that
$M\in\sV_R$ lies in
$\sV^{\ff}_R$ if and only if $id_M^\infty =0$.

A morphism $f$ is said to be {\it Fredholm} if
$f^\infty$ is invertible. Denote by
$\sV^{\fFr}_R$ the category of projective
modules and Fredholm morphisms. Let
$\sV^\times_R$, ${\sV}^{\infty\times}_R$ be
groupoids of invertible morphisms in $\sV_R$,
${\sV}^\infty_R$. Thus $\sV^\times_R \subset
\sV^{\fFr}_R$ and the morphism of groupoids
$\sV^\times_R \to
{\sV}^{\infty\times}_R$ extends to the
functor
$\sV^{\fFr}_R \to{\sV}^{\infty\times}_R$.

Any base change of a Fredholm morphism $f$ is
  Fredholm. So for every geometric point
$x$ of
$\Spec R$ we have a Fredholm operator $f_x :
N_x \to M_x$.\footnote{Here the
$k_x$-vector space $M_x  := k_x
\otimes M$ is the fiber of $M$ at $x$.} Denote
by
$i(f)_x =i(f_x )$ its index,
$i(f)_x := \dim \Coker f_x -\dim \Ker f_x$.
The index depends only on the corresponding
asymptotic morphism, so we write $i(f)_x
=i(f^\infty )_x$.

\subsection{Lemma} For a morphism $f: N\to M$
of projective modules the following conditions
are equivalent:

(i) $f$ is Fredholm,

(ii) $\Coker f$ is
finitely generated, for every geometric
point $x$ of $\Spec R$ the corresponding
operator $f_x : N_x \to M_x$ is Fredholm, and
the function $i(f): x\mapsto i(f_x )$ on $\Spec
R$ is locally constant.

(iii) $f$ can be written as composition
\eq{2.5.1}{N\buildrel{i}\over\to N\oplus Q
\buildrel{\tilde{f}}\over\to M\oplus P
\buildrel{p}\over\to M} where
  $P,Q\in\sV^\ff$, $i$ is
the embedding, $p$ the projection, and
$\tilde{f}$ is an isomorphism.

\begin{proof} It is clear that (iii) implies
(i) and (ii).

(i)$\Rightarrow$(ii): Since $f^\infty$ admits
a right inverse $\Coker f$ is finitely
generated. So there
exists $Q\in\sV^\ff$ and a morphism $\alpha
:Q\to M$ such that
$\pi :=(f,\alpha ):N\oplus Q \to M$ is
surjective. Then $P:=\Ker \pi \in\sV$, and (i)
assures that it has finite rank. Since $i(f)=
rk (Q)-rk (P)$, it is a locally constant
function.

(i)$\Rightarrow$(iii): Define
$\alpha$, $Q$, $P$ as above. A splitting of
$\pi$ yields a projector $\beta : N\oplus Q \to
P$. Then
$\tilde{f}:= (f,\alpha;\beta)  :N\oplus Q \iso
M\oplus P$ defines
  (2.5.1).

(ii)$\Rightarrow$(iii): We define $\tilde{f}$
as above. Our $P$ is projective,
so it remains to check that the conditions of
(ii) imply that it is finitely generated. We
know that every fiber $P_x$ is of finite
dimension, and this dimension is locally
constant with respect to $x$. Take any $x\in
\Spec R$; we want to find a Zariski
open $x\in U\subset \Spec R$ such that $P|_U$
is finitely generated. By 2.3(ii) we can
find $U$ such that $P|_U$ can be written as
$P'\oplus P''$ where $P'$ is finitely
generated and
$P'_x =P_x$, i.e., $P''_y =0$. Shrinking
$U$ if necessary to assure that $\dim P'_y$ is
constant on $U$, we see that
$P''_y =0$ for every $y\in U$. By 2.3(i) one
has
$P''=0$.
\end{proof}

{\it Remarks.}
(a) As follows from 2.5(ii), Fredholm
morphisms have local nature with respect to
the flat topology.

(b) If $I\subset R$ is a
nilpotent ideal then $f$ is Fredholm if and
only if $f_{R/I} :N/IN \to M/IM$ is Fredholm
(use 2.5(ii)).

(c) If $R$ is Noetherian
then the last condition in 2.5(ii)
is superfluous. Indeed, the proof of
(ii)$\Rightarrow$(iii) proceeds as above, but
we just notice that $P''_x =0$ implies, after
shrinking $U$ if necessary, that
$P''_y =0$ for every generic point of $U$.
So $P''
=0$ by 2.3(i).

\subsection{Relative determinant lines  } Of
course, there is no way to assign a determinant
line  to a projective module of infinite rank.
However if one has two such modules $M,N$ and
an asymptotic isomorphism $f^\infty$ between
them then one can use $f^\infty$  to
mutually cancel infinities in $\det M$, $\det
N$ so that the ratio $\det M /\det N$ (the
relative determinant line) is well-defined.
Let us explain this Dostoevskian
ansatz\footnote{``Let a reptile gobble up
another of its kind,"    a comment of a hero
of ``The brothers Karamazov" on  his
father and the elder
brother.}$^,$\footnote{One may prefer to
think of two sumo wrestlers, each of
infinite girth, locked in a Fredholm grip,
with the match determined by a mere finite
inbalance of forces.} in more details.

\medskip
Here is a list of data we look for:

(i) {\it A $\Z$-graded super
extension
${\sV}^{\infty\flat}_R$ of the groupoid
${\sV}^{\infty\times}_R$.} Thus for
every triple
$(M,N,f^\infty )$, where $M,N\in\sV_R$,
$f^\infty$ is an invertible asymptotic morphism
$N\to M$, we want to have a $\Z$-graded super
line
$\det (M,N,f^\infty )$  (the relative
determinant line), together with data of
composition isomorphisms
\eq{2.6.1}{c: \det
(L,M,g^\infty ) \cdot \det (M,N,f^\infty )\iso
\det (L,N, g^\infty f^\infty )}
satisfying the obvious associativity property.

(ii) {\it A splitting of the
pull-back  of
${\sV}^{\infty\flat}_R$ to $\sV^\times_R$.}
This means that for every isomorphism
$f: N\iso M$ in $\sV$ we have a canonical
trivialization
\eq{2.6.2}{\det f : 1_R \iso \det (M,N,
f^\infty )}  compatible with
composition of
$f$'s (via (2.6.1)).

(iii) {\it Compatibility with  sums.} For
every $(M_1 ,N_1 ,f^\infty_1 )$,
$(M_2 ,N_2 ,f^\infty_2 )$ we  have a canonical
isomorphism $\det (M_1 \oplus M_2 ,N_1 \oplus
N_2 ,f^\infty_1 \oplus f^\infty_2 )=\det (M_1
,N_1 ,f^\infty_1 )
\cdot\det (M_2 ,N_2 ,f^\infty_2 )$ compatible
with the associativity and commutativity
constraints. It is nicer to write it then as
  a canonical identification of
$\Z$-graded super lines \eq{2.6.3}{
\det (\oplus M_\alpha ,\oplus N_\alpha
,\oplus f^\infty_\alpha )=\otimes \det
(M_\alpha ,N_\alpha ,f_\alpha^\infty )}
defined for any
finite
family
$\{ (M_\alpha ,N_\alpha ,f^\infty_\alpha )\}$.
We want it to be compatible with the canonical
isomorphisms from (i), (ii).

(iv) {\it Compatibility with finite
determinants.} For every
$M,N\in
\sV^\ff$ there is a canonical
identification\footnote{We
omitted
$f^\infty$ in the l.h.s$.$ of (2.6.4) for
it bears no information: in the present
situation one has
$\Hom^\infty (N,M)=\text{Isom}^\infty
(N,M)=0$.}
\eq{2.6.4}{\det (M,N )=\det M /\det N} such
that (2.6.1)--(2.6.3) become the
usual compatibilities from 2.2.

  (v) {\it Base change.} Our structure should be
compatible with the base change in the obvious
manner.

Notice that (iv) and (ii) imply that the
degree of the super line $\det (M,N,f^\infty )$
is the index $i(f^\infty )$ (see 2.4).

\subsection{Proposition} Such a structure
exists and is unique (up to a unique
isomorphism).

\begin{proof} Here is an explicit
construction. We will not use the uniqueness
statement, so its proof is
left to the reader.\footnote{Hint: use 2.9.}

\medskip

  Let us construct the
determinant line $\det (M,N,f^\infty )$.
Assume for a moment that $M$ satisfies the
following assumption, which is fulfilled, e.g.,
if
$M$ is a free $R$-module:

$(*)$ For every finitely generated
$R$-submodule
$T\subset M$ there exists a finitely
generated $R$-submodule $P\subset M$ with
$M/P$ projective
  such that $T\subset
P$.

Let $f:N\to M$ be any
lifting of $f^\infty$. Since
$\Coker f$ is finitely generated (see
2.5.(ii)),  we can choose
$P
\subset M$ as in $(*)$ such that $P+
f(N)=M$. Set
$Q := f^{-1}(P )$. Then $P
,Q \in
\sV^\ff$.

\medskip

{\bf Lemma-definition.} The $\Z$-graded super
lines
$\det P /\det Q$ for all choices of
$f$,
$P$ are canonically identified, so they can
be considered as a single super line. This is
our
  $\det (M,N,f^\infty )$.

{\it Proof of Lemma.} We fix $f$ for a
moment. Let $P'$ be another submodule as
above, $Q':=f^{-1}(P')$. Let us
define the canonical identification
$\phi =\phi_{PP'}:\det P
/\det Q \iso
\det P'/\det Q'$.

If $P'\supset P$ then
$P'/P \in \sV^\ff$ and $f$ yields an
isomorphism $\bar{f} : Q'/Q\iso P'/P$. Our
$\phi$ is the image of $\det
\bar{f}$ by the canonical
identification $\det (P'/P)/\det
(Q'/Q)=(\det P' /\det Q' )/(\det P/\det Q)$
(see (2.2.2)).

If $P'$ is arbitrary then we can
choose   $P''$ as in $(*)$ such that $
P''\supset P+P'$, and set
$\phi_{PP'}:=
\phi^{-1}_{P'P''}\phi_{PP''}$. One
checks in a moment that this definition does
not depend on the choice of $P''$ and
$\phi$ so defined satisfy the
transitivity property. Thus all super
lines for fixed lifting $f$ of $f^\infty$ are
canonically identified.

Now let $f'$ be another lifting of $f^\infty$.
Choose $P$ as in $(*)$ so that
$P +f(N)=M$ and
$P \supset (f-f')(N)$. Then $P
+f'(N)=M$ and
$f^{-1}(P)=f^{'-1}(P)$. So the super
lines for
$f$ and $f'$ computed by means of such $P$
are simply equal. If we change $P$ the
identifications $\phi$  for $f$ and $f'$ are
also equal. These identifications of super
lines are transitive with respect to $f$. We
are done.
\hfill$\square$

\medskip

  In order to eliminate condition
$(*)$ notice that super lines are local
objects, so it suffices to define $\det$
locally on $\Spec R$. Now one can use
2.3(ii) instead of $(*)$.

{\it Remark.} One can also reduce the general
situation to the situation of $(*)$ directly
using the following trick: Choose $M'\in\sV$
such that
$M\oplus M'$ is free (and so satisfies $(*)$);
then
$\det (M,N,f^\infty )=\det (M\oplus
M',N\oplus N',f^\infty \oplus id_{M'}^\infty
)$.

\medskip

The canonical isomorphisms
(2.6.1)--(2.6.3) are defined as follows. We
use condition $(*)$ at will; in truth, one
should work locally on $\Spec R$ and use
2.3(ii).

To define the composition isomorphism $c$
of (2.6.1)  we choose some liftings $f$, $g$
and compute the determinant lines by means of
finitely generated $S
\subset L$ such that $L/S$ is projective
and $S + gf (N)=L$ and $P :=
g^{-1}(S )$, $Q := (gf)^{-1}(S
)=f^{-1}(P)$. Our $c$ is the standard
isomorphism $(\det S /\det P )\cdot (\det P
/\det Q )
\iso \det S /\det Q$. The independence of
choices is immediate.

To define $\det f$ of (2.6.2) one computes
$\det (M,N, f^\infty )$ using $f$ and $P
=0$.

To define (2.6.3) you compute the r.h.s.~using
some $f_\alpha$ and $P_\alpha$ and then compute
the l.h.s.~using
$f=\oplus f_\alpha$ and $\oplus
P_\alpha$; then use (2.2.1).

2.6(iv), 2.6(v) are inherent in the
construction. The mutual compatibilities of the
above canonical isomorphisms are obvious.
\end{proof}

\subsection{Remark}
  There is a more general
construction (which is not needed in this
article) that
assigns a determinant line to every pair
$(M^\cdot ,d^\infty )$ where
$M^\cdot$ is a finitely $\Z$-graded projective
module and $d^\infty$ is a differential of
degree 1 on
$M^\infty$ such that $(M^\cdot ,d^\infty )$
is a homotopically trivial complex. The
construction above corresponds to $M^\cdot$
supported in degrees $0,-1$.

\subsection{Compatibility with filtrations}
Let $M,N$ be projective modules equipped with
finite (increasing) filtrations $M_\cdot$,
$N_\cdot$ such that the associated graded
modules gr$M$, gr$N$ are projective (i.e., the
filtrations can be split). Let
$f^\infty :N\to M$ be an asymptotic morphism
compatible with filtrations such that
gr$f^\infty : \text{gr}N \to\text{gr}M$ is
invertible. Then $f^\infty$ is also
invertible.

{\bf Lemma.} There is a canonical
isomorphism
\eq{2.9.1}{\det (M,N,f^\infty )=\det
(\text{gr}M,\text{gr}N,\text{gr}f^\infty
)}  compatible with base change and such that
(2.9.1) is the identity map provided that our
filtered picture came from a graded one. Such
data of isomorphisms is unique.

{\it Proof.} Here is a construction of (2.9.1)
that makes use of general functoriality
from 2.6 only.

  Consider the Rees $R[t]$-module $M^\fr
:=\oplus M_i t^i$. This is a projective
$R[t]$-module equivariant with respect to
homotheties of $t$ whose fibers at $t=1$ and
$t=0$ are, respectively,
$M$ and gr$M$.
  Similarly, we have
$N^\fr$ and an invertible asymptotic
$R[t]$-morphism $f^{\infty\fr}: N^\fr \to
M^\fr$. The corresponding super $R[t]$-line
$D := \det
(M^\fr ,N^\fr ,f^{\infty\fr })$ is
equivariant with respect to homotheties of
$t$ (by base change), i.e., $D$ is a graded
$R[t]$-module. Its fibers at $t=1$ and $t=0$
are, respectively,
$\det (M,N,f^\infty )$ and $ \det
(\text{gr}M,\text{gr}N,\text{gr}f^\infty
)$. The fiber at $0$ has degree
$i=i(\text{gr}f^\infty )=i(f^\infty )$ with
respect to the action of homotheties. Let
$D^i$ be $i^{\rm th}$ component of $D$. This
is a super $R$-line, and $D=R[t]\otimes_R
D^i$. The restriction to $t=1,0$
identifies $D^i$ with $\det (M,N,f^\infty )$ and $ \det
(\text{gr}M,\text{gr}N,\text{gr}f^\infty
)$.
The composition of these identifications is
(2.9.1).

The properties mentioned in the statement of
Lemma, as well as the uniqueness, are
immediate.
\hfill$\square$

\medskip

{\it Remarks.} (i) Here is another, direct,
construction of (2.9.1). Let us construct
$\det (M,N,f^\infty ) $ as in the proof of 2.7
by means of some $f$ and
$P\subset M$. Choose
$f$ compatible with filtrations and $P$ such
that gr$P\subset$gr$M$  is such that
gr$M/$gr$P$ is projective
  and gr$P
+\text{gr}(f)(\text{gr} N)=\text{gr}
M$; then use (2.2.2) to define (2.9.1).
The independence of auxiliary choices is
immediate, as well as properties mentioned in
the statement of Lemma.

(ii) Isomorphisms (2.9.1) are compatible with
(2.6.1) and (2.6.2); if $M,N$ are of finite
rank then (2.6.4)  reduces them to (2.2.2).

\subsection{The Tate extension.} Below we will
consider set-valued functors on the
category of commutative $R$-algebras which
commute with  finite inverse limits. The
``category" of such functors (called
{\it ``$R$-spaces"}) is denoted  by
$\sS_{R}$; it is closed under finite inverse
limits. For an $R$-algebra $B$ the
corresponding representable functor is denoted
by $\Spec B\in\sS_R$. We have a ring object
$\sO =\Spec R[t] \in\sS_R$, $\sO (R')=R'$, so
one can consider
$\sO$-modules, $\sO$-algebras in 
$\sS_R$, etc. For $X\in\sS_{R}$ and an 
$R$-flat commutative $R$-algebra $B$,
$\Spec B=:Y$, we have an $R$-space
$\sH om (Y ,X)$, $\sH om
(Y ,X)(R'):= X(R'\mathop\otimes\limits_R
B)$. For $B=R[h]/h^2$ we get the {\it tangent
bundle} $\Theta_X $ of $X$. Since
$R[h]/h^2$ is naturally an $R$-module in the
category of commutative $R$-algebras $A$
equipped with a morphism $A\to R$, the tangent
bundle is an $\sO$-module over $X$, so for any
$x\in X (R)$ the tangent space
$\Theta_{X,x}$ (:= the fiber of $\Theta_X$ at
$x$) is an $\sO$-module in
$\sS_R$. For a group $R$-space $G$ its
tangent space at $1\in G(R)$ is denoted by
Lie$G$; this is a Lie $\sO$-algebra in
$\sS_R$.\footnote{See SGA3 t.1 exp.II,
pp.26-27. Here is a sketch of an argument: (a)
Check that the addition law on Lie$G$ comes
from the product on $G$. (b) To define the Lie
bracket of $\alpha ,\beta\in
$ (Lie$G) (R)$ consider them as elements of,
respectively, $G(R[s]/s^2 ),G(R[t]/t^2 )
\subset G(R[s,t]/(s^2 ,t^2 )$; then the
commutator
$(\alpha,
\beta ):=\alpha\beta\alpha^{-1}\beta^{-1}$
belongs to $G (R[h]/h^2 )\subset
G(R[s,t]/(s^2,t^2 ))$ where $h:=st$; this is
$[\alpha ,\beta ]$. (c) The Lie bracket is
obviously skew-symmetric. To see  Jacobi
 consider $\alpha ,\beta, \gamma\in
$ (Lie$G) (R)$ as elements of, respectively, 
$G(R[s]/s^2 ), G(R[t]/t^2 ), G(R[u]/u^2 )
\subset G(R[s,t,u]/(s^2 ,t^2, u^2 ))$. Then
$[\alpha, [\beta,\gamma ]]$ corresponds to
$(\alpha ,(\beta ,\gamma ))\in G(R[v]/v^2)$
where $v:= stu$, and Jacobi follows since
$(\alpha ,\beta )$ commutes with $(\alpha
,\gamma )$, $\alpha$,  $\beta$. }

Every projective $R$-module $M$ defines an
$\sO$-module in $\sS_{R}$ which we denote
also by $M$, $M(R') :=
M_{R'}:=M\mathop\otimes\limits_R R'$. For
$M,N\in\sV_R$ we have  $\sO$-modules $\Hom
(N,M),\Hom^\infty
(N,M)$ in $\sS_R$, $\Hom
(N,M)(R'):=\Hom_{R'} (N_{R'},M_{R'})$,
  $\Hom (N,M)^\infty
(R'):=\Hom^\infty_{R'} (N_{R'},M_{R'})$. So
for $M\in\sV_R$ we have the corresponding
$\sO$-algebras of endomorphisms and the groups
of invertible endomorphisms $\GL (M),
\GL^\infty (M)=\GL (M^\infty )$. There is
  a canonical homomorphism $\GL (M)\to \GL
(M^\infty )$; denote by
$\GL^\ff (M)$  its kernel.

By 2.6(v) the $\Z$-graded super
extensions
$\Aut^{\infty\flat}_{R'}(M_{R'})$ form the {\it
Tate} super
$\sO^\times$-extension $\GL^{\flat} (M^\infty
)$ of
$\GL (M^\infty )$.  For $g^\infty \in \GL
(M^\infty )(R')$ we denote the corresponding
  $\Z$-graded super line in $\GL^\flat (M^\infty
)(R')$
  by
$\det (M^\infty , g^\infty )$ or simply
$\lambda_{g^\infty}$.

Here
is a list of basic properties of the Tate super
extension:

\medskip

(i) According to Remark (iii) in A3, the Tate
extension depends functorially on
$M$ {\it considered  as an object of}
${\sV}^{\infty\times}_R$, i.e., every
asymptotic isomorphism $f^\infty :N\to M$
yields a canonical identification of $\Z$-graded super
extensions
\eq{2.10.1}{\text{Ad}_{f^\infty}^\flat :
\GL^{\flat} (N^\infty )\iso
\GL^{\flat} (M^\infty ).} To write
it down explicitly, consider
$\sV_R^{\infty\flat}$ as a plain
$\sO^\times$-extension of
$\sV_R^{\infty\times}$ (see 2.1) and choose
(locally on $\Spec R$) a generator
$\tilde{f}$ of the line
$ \det (M,N,f^\infty )$. Now for
$g^\infty \in \GL (N^\infty )$ the map
$\Ad_{f^\infty}^\flat : \lambda_{g^\infty}\to
\lambda_{\Ad_{f^\infty}(g^\infty )}$ is the
adjoint action of $\tilde{f}$ (in our
plain $\sO^\times$-extension) multiplied by
$(-1)^{p(f^\infty )p(g^\infty )}$ where $p$ is
parity (i.e., the index mod 2).

In particular, the adjoint action of
$\GL (M^\infty )$ lifts canonically
to an action of $\GL (M^\infty )$ on
the Tate extension $\GL^{\flat} (M^\infty )$.

{\it Remark.} The groupoid
${\sV}^{\infty\times}$ (as opposed to
$\sV^\times$) does not satisfy the flat
descent property (see \cite{Dr} for a
discussion of this subject). One  can sheafify
it formally, and the above property assures
that automorphism groups of objects of this
stack have a canonical super
$\sO^\times$-extension.

\medskip
(ii) The map
$\GL (M)\to \GL (M^\infty )$, $g\mapsto
g^\infty$, lifts canonically to a homomorphism
  \eq{2.10.2}{\GL (M)\to
\GL^{\flat} (M^\infty ),
\quad g\mapsto \det (M, g)\in \det
(M^\infty ,g^\infty ) .} Namely, for $g\in
\GL (M)$ one has $\det (M^\infty , g^\infty
)=\det (M,M,g^\infty )$, and $\det (M,g)$ is
  $\det g$ from (2.6.2).

{\it Exercises.} (a) Show that the restriction
of (2.10.2) to the subgroup $\GL^\ff (M)$
is the usual determinant
$\det : \GL^\ff (M)\to\sO^\times$.

(b) Assume that $f^\infty
,g^\infty \in \Aut^\infty (M)$ commute and
$f^\infty$ can be lifted to
$f\in \Aut (M)$. Then $f$ acts by
functoriality on the determinant line $\det
(M,M,g^\infty )$. Show that this action is
multiplication by $\{ f^\infty ,g^\infty
\}^\flat $.\footnote{See A5 for notation.}

\medskip

(iii) For a finite collection $\{
M_\alpha \}$ the restriction of
$\GL^{\flat} (\oplus M_\alpha^\infty
)$ to the subgroup $\Pi
\GL (M_\alpha^\infty )\subset \GL
(\oplus M_\alpha^\infty )$ identifies
canonically with the Baer product of the
Tate super extensions $\GL^{\flat}
(M_\alpha^\infty )$ (see 2.6(iii)). In other
words, for every
$\alpha$ the embedding $\GL
(M_\alpha^\infty )\hra
\GL (\oplus M_\alpha^\infty )$ lifts
canonically to an embedding of super
extensions
$\GL^{\flat} (M_\alpha^\infty ) \hra
\GL^{\flat} (\oplus M_\alpha^\infty )$, and for
different $\alpha$'s the images of these
embeddings commute.\footnote{i.e.,
for $\gamma \in \GL^{\flat}
(M_\alpha^\infty )$, $\gamma' \in \GL^{\flat}
(M_{\alpha'}^\infty )$, $\alpha\neq\alpha'$,
one has $\{ \gamma ,\gamma' \}^\flat =1$, see
A5.}

{\it Remark.} This property is no longer
true if we consider our super extensions as
plain $\sO^\times$-extensions (see 2.1).

\medskip

(iv) Assume that an object $M^\infty \in
{\sV}^\infty_R$ has  filtration which admits a
splitting; let
$B\subset \GL (M^\infty )$ be the subgroup
preserving this filtration. Then the
restriction $B^\flat$ of $\GL^\flat (M^\infty
)$ to $B$ is canonically isomorphic to the
pull-back of $\GL^\flat (\text{gr}M^\infty )$
by the map $B \to \GL (\text{gr}M^\infty )$
by 2.9.

\medskip

Consider the Lie algebras $\fgl (M),
\fgl^\ff (M),\fgl (M^\infty ),\fgl^\flat
(M^\infty )$ of our groups
$\GL (M)$, etc. The first three of them are
just the corresponding algebras of
endomorphisms considered as Lie algebras, so
$\fgl (M^\infty )=\fgl (M)/\fgl^\ff (M)$, and
$\fgl^{\flat}(M^\infty)$ is a central
extension of $\fgl (M^{\infty})$ by $\sO$. By
(ii) the projection $\fgl (M)\to
\fgl (M^\infty )$ lifts canonically to a
morphism of Lie algebras $\fgl (M)\to
\fgl^{\flat}(M^\infty )$; on $\fgl^\ff (M)$
this morphism is the trace map $tr :\fgl^\ff
(M)\to \sO$.\footnote{Indeed, $tr$ is tangent
  to $\det$ from Exercise (a) above.} Therefore
one has:

(v) $\fgl^{\flat}(M^\infty )$ is the
push-forward of the extension $0\to \fgl^\ff
(M)\to\fgl (M)\to \fgl (M^\infty )\to 0$ by
the ad-invariant morphism $tr$.

\medskip

(vi) Since $\fgl^{\flat}(M^\infty)$ is a central
extension of $\fgl (M^{\infty})$ we can
rewrite its bracket as a pairing $[\, ,\,
]^\flat :\fgl (M^\infty)\times
\fgl (M^\infty)\to\fgl^\flat (M^{\infty})$.
It satisfies the {\it  cyclic  identity}: for
every
$a,b,c\in\fgl (M^\infty )$ one has
\eq{2.10.3}{[ab,c]^\flat +[ca,b]^\flat
+[bc,a]^\flat =0.} Here $ab$ is the
product of $a$, $b$ as endomorphisms of
$M^\infty$.\footnote{Lift $a, b,c$ to
$\tilde{a},\tilde{b},\tilde{c}\in\fgl (M)$.
By (v), it suffices to show that
$[\tilde{a}\tilde{b},\tilde{c}]
+[\tilde{c}\tilde{a},\tilde{b}]
+[\tilde{b}\tilde{c},\tilde{a}]=0$. This is
immediate.}

\subsection{The topological setting: Tate
$R$-modules}  An unpleasant feature of the
above formalism is the absence of duality: the
dual of a projective module of an infinite
rank is not there. A way to recover the
duality is to consider the setting of Tate
modules.

If $R$ is
a commutative ring then a {\it topological
$R$-module} is an
$R$-module  equipped with a topology\footnote{The topologies we consider
are always compatible with the additive
structure; ``base" means ``base
of neighbourhoods of 0."}  which
has a base formed by
$R$-submodules; we always assume that the
topology is complete and separated.
Topological $R$-modules form an additive
$R$-category. A morphism
$R\to R'$ defines an obvious functor from the
category of topological
$R'$-modules to that of $R$-modules; its
left adjoint is the base change functor
$F_R \mapsto
F_{R'}:= F_R \mathop{\hat{\otimes}}\limits_R
R'$.

We can consider every plain $R$-module $N_R$ as
a topological module with discrete
topology. If $M_R$ is a {\it projective}
$R$-module then its dual $M_R^* =\Hom_R
(M_R ,R)$ equipped with the weak topology
(formed by annihilators of finite subsets in
$M_R$) is a topological
$R$-module. The canonical map $M_R \to (M_R^*
)^* :=
\Hom^{cont}(M_R^* ,R)$ is an isomorphism, and
the functor $M_R \mapsto M^*_R$ is fully
faithful. For every
$R\to R'$ one has
$M_{R'}^* =(M_R^*)_{R'}$.\footnote{To see this
represent $M_R$ as a direct summand of a free
$R$-module.}

{\it Remark.} For any $N_R$, $M_R$ as above
the image of every continuous morphism $M^*_R
\to N_R$ is a finitely
generated $R$-module.\footnote{To see this
notice that we can assume that $M_R$ is a free
$R$-module.}

\medskip
We call a projective $R$-module considered as
  discrete topological
$R$-module a {\it discrete Tate}
$R$-module. The corresponding duals are called
  {\it compact Tate} $R$-modules.
A topological $R$-module $F_R$ is
  {\it special Tate} $R$-module if it can be
represented as the direct sum of a discrete and
a compact Tate $R$-modules. Following
Drinfeld \cite{Dr}, we define a {\it Tate}
$R$-module as a topological
$R$-module which is a direct summand of a
special Tate $R$-module. Notice that a
direct summand of a compact Tate module is a
compact Tate module; same is true for discrete
Tate modules.  The base change sends
Tate modules to  Tate modules. All base changes
of a given Tate
$R$-module  considered simultaneously (see
2.10) form a {\it Tate
$\sO$-module} in $\sS_R$.

\medskip

{\it Examples.} (i) $R((t))$ equipped with
topology with base $t^n R[[t]]$, $n\ge 0$, is a
special Tate $R$-module.

(ii) Let $F$ be a finitely generated
projective $R((t))$-module. Then $F$ carries a
canonical topology  whose base is formed
by $R[[t]]$-submodules of $F$ which generate
$F$ as an $R((t))$-module.
$F$ is a Tate
$R$-module.

(iii) Let $k$
be a commutative ring and  $R\subset k[x]$  the
subalgebra of polynomials $f=f(x)$ such that
$f(1)=f(0)$. Let $F_R \subset k[x]((t))$ be
the $R$-submodule of elements $g=g(x,t)$ such
that $g(1,t)=tg(0,t)\in k((t))$. It carries a
topology with base  $F \cap t^n k[x]((t))$,
$n\ge 0$. Then $F$ is a
Tate $R$-module which is not special.
Moreover, it is not special  Zariski
locally on $\Spec R$.

  \medskip

The following result is due to Drinfeld
  \cite{Dr} 8.1, 4.1:

\subsection{Proposition} (i) Let $M$ be a
projective $R$-module. Then every asymptotic
projector
$\pi^\infty
\in
\text{End}_R^\infty (M)$, $(\pi^\infty )^2
=\pi^\infty$, can be lifted
to a true projector
$\pi \in \text{End}_R M$
  \'etale locally on
$\Spec R$.

(ii) Every Tate
$R$-module is special  \'etale
locally on
$\Spec R$.

{\it Remark.} The proof actually shows that the
statements hold Nisnevich locally.

\begin{proof} (i) Choose any lifting of
$\pi^\infty$ to $\text{End}_R M$. Let us
consider $M$ as an
$R[t]$-module where $t$ acts on $M$ by this
lifting. We will find (after an appropriate
localization of $R$) a polynomial $p=p(t)\in
R[t]$ such that (a) $p(p-1)$ kills
$M$, and (b) $p(0)=0$, $p(1)=1$. Let $\pi$ be
  the action of $p$ on $M$. This is a
projector by (a) which lifts
$\pi^\infty$ because of (b).\footnote{Indeed,
(b) means that $p(t)-t$ is
divisible by $t(t-1)$, so the corresponding
asymptotic endomorphism of
$M$ vanishes (since the one for $t(t-1)$
does).}

Choose a finitely generated
$R$-submodule
$L\subset M$ which contains $t(t-1
)(M)$. There is a monic $f(t)\in
R[t]$ which kills $L$, so
$t(t-1) f(t)$ kills $M$. Localizing $\Spec R$
with respect to \'etale\footnote{In fact,
Nisnevich} topology we can assume that
$f(t)=g(t)h(t)$ where $g(1)$ and $h(0)$ are
invertible, and
$g(t), h(t)$ generate the unit ideal. Then $M$
is supported on the union of
non-intersecting subschemes
$tg(t)=0$ and $(t-1)h(t)=0$. Our $p(t)$ is any
polynomial which vanishes on the first
subscheme and equals to 1 on the second one.

\medskip

(ii) Let $F$ be a Tate $R$-module. Choose a
special Tate $R$-module $G$ such that $F$ is
a direct summand of $G$. We will find, after
an appropriate localization of $R$, an open
submodule $L\subset G$ such that (a) $L$ is a
compact Tate module, and (b)  $P:=
F/F\cap L$ is a projective $R$-module.
Then $F$ is a special Tate module. Indeed,  a
section $\gamma : P\to F$ yields $F\iso
P\oplus F\cap L$. Now $F\cap L$ is a direct
summand of $L$,\footnote{A projector $\Pi
:G\twoheadrightarrow F$ yields a projector
$L \twoheadrightarrow F\cap L$, $\ell
\mapsto \Pi (\ell )-\gamma (\Pi (\ell
)\text{mod}L)$.}  hence it is a compact Tate
module.

Let
  $\Pi \in
\text{End}G$ be a projector such that
$\Pi (G)=F$, and $G= M\oplus N^*$ be any
decomposition of $G$ into a sum of a discrete
and compact Tate modules. Let
$\phi$ be the composition $M\hra
G\buildrel{\Pi}\over\to G\twoheadrightarrow
M$. Then $\phi^\infty$ is an asymptotic
projector (see Remark in 2.11). By (i), after
an appropriate
  localization of $\Spec R$, one can
find a projector $\pi \in \text{End}_R M$ such
that $\phi^\infty =\pi^\infty$. The  images of the composition $N^* \hra G
\buildrel{\Pi}\over\to G\twoheadrightarrow M$
and of $\pi -\phi$
are contained in a finitly generated submodule
$K\subset M$ (see Remark in 2.11). By 2.3(ii)
applied to the images of $K$ in $\Ker \pi$ and
in Im$\pi$ one can find finitely generated
$R$-submodules $L' \subset \Ker \pi$ and $L''
\subset \text{Im}\pi$ such that $\Ker \pi /L'$
and Im$\pi /L''$ are projective $R$-modules.
Set $L:= N^* \oplus L' \oplus L'' \subset G$.
Since $L' ,L''$ are projective $L$ satisfies
(a) above; since
$F/F\cap L =
\text{Im}\pi /L''$ it satisfies (b).
  We are done.
\end{proof}

\medskip

All constructions of this article
are local with respect to the flat topology of
$\Spec R$, so we will often tacitly
assume that the Tate modules we deal with are
  special.

{\it Remark.} According to \cite{Dr}8.3,
Tate $R$-modules, as opposed to special Tate
$R$-modules, satisfy the flat descent property.

\subsection{Duality, lattices, and the Tate
extension}
Let $F$ be a Tate
$R$-module. We say that an $R$-submodule
$L\subset F$ is {\it bounded} if for every
discrete $R$-module $P$ and a continuous
morphism $F\to P$ the image of $L$ in
$N$ is contained in a finitely generated
$R$-submodule of $P$. If $F$ is realized as a
direct summand in a special Tate module
$G=M\oplus N^*$, $M,N\in\sV_R$, then $L$ is
bounded if and only if the image of $L$ in
$M=G/N^*$ is contained in a finitely generated
$R$-submodule.\footnote{This follows from
Remark in 2.11.} Notice that bounded open
submodules of $F$ form a base of the topology
of $F$.

For a Tate $R$-module $F$ its {\it
dual}
$F^*$ is
$\Hom^{cont}_R (F,R)$ equipped with a topology
whose base is formed by orthogonal complements
to bounded $R$-submodules of $F$. Then $F^*$
is again Tate module, and the canonical
morphism $F\to (F^* )^*$ is an isomorphism.
To see this notice that duality commutes with
(finite) direct sums, and for a special Tate
module $F=M\oplus N^*$, $M,N\in \sV_R$, one
has $F^* = N\oplus M^*$.
So the duality
functor is an anti-equivalence of the
category of Tate $R$-modules. It
commutes with the base change.

\medskip

For a  Tate $R$-module $F$ a {\it c-lattice} in
$F$ is an open $R$-submodule $L\subset
F$ which satisfies either of the following
conditions (to check  their equivalence
is an exercise for the reader):

(i) $L$ is a compact Tate module and the
projection $F\twoheadrightarrow F/L$ admits an
$R$-linear section, 

(ii) $L$ is a bounded submodule of $F$ such
that $F/L$ is a projective $R$-module.

Of course, a c-lattice exists if and only if
$F$ is a special Tate module.

  A {\it
d-lattice} in $F$ is an $R$-submodule
$M\subset F$ complementary to some c-lattice.

{\it Remarks.} (i) c-lattices need not form a
base of the topology of $F$. However, by
2.12(ii) and 2.3(ii), they form a base
Nisnevich locally.

(ii) If $L\supset L'$ are c-lattices then
$L/L' \in \sV^\ff_R$.

(iii) Every d-lattice is a
projective $R$-module. If $M\supset M'$ are
d-lattices then $M/M' \in \sV^\ff_R$.

\medskip

  Let $L,L'$ be c-lattices in a special Tate
module $F$. Every splitting $F /L \hra F$
yields a morphism
$F /L\to F /L'$. The corresponding
asymptotic morphism (see 2.4) is an
isomorphism that does not depend on the choice
of splitting. Therefore the objects $(F
/L)^\infty \in\sV_R^\infty$, $L$ is any
c-lattice in $F$, are canonically identified.
They form a {\it canonically defined object} of
the category
${\sV}^\infty_R$ which we denote by
$F^\infty$.  Our $F^\infty$ depends on
$F$ in a functorial way. This construction is
compatible with the base change.

  For  $L, L'
$ as above set \eq{2.13.1}{\det (L:L'):=\det
(F/L' ,F/L , id_{F^\infty} ).} One
has canonical composition isomorphisms of
$\Z$-graded super lines $\det
(L:L')\cdot
\det (L' :L'') \iso \det (L:L'' )$, so we have
defined a super pre-gerbe structure on the set
of c-lattices. The corresponding
$\Z$-graded super gerbe is called  the {\it
gerbe of c-lattices.}

For every $L\supset L'$ there is a
canonical identification $\iota :\det (L:L'
)\iso \det (L/L')$ such that for $L\supset L'
\supset L''$ the composition becomes the
standard isomorphism
$\det (L/L' )\cdot \det (L'
/L'')\iso \det (L/L'' )$. 

{\it Remark.}  A super pre-gerbe
structure on the set of c-lattices together
with identifications
$\iota$ and compaitibility with base change is
unique up to a unique isomorphism (use Remark
(i) above).

\medskip

Consider the semigroup of endomorphisms $g$ of
$F$ such that $g^\infty$ is invertible.  
Pulling back the Tate super extension by the
 homomorphism $g\to g^\infty$ to
$\Aut (F^\infty )$ one gets a super
extension of our semigroup. In
particular, we have  a super extension
$\Aut^\flat (F )$ of $\Aut (F )$. 

Notice that if $L$ is a c-lattice
such that
$L+g(F)=F$ then \eq{2.13.2}{\lambda_g
:= \lambda_{g^\infty}=
\det (L:g^{-1}(L)). } If $g\in\Aut (F)$ then
any $L$ will do, and, acting by $g$ on the
r.h.s., we get a canonical isomorphism
$\lambda_g =\det (g(L):L)$. 

{\it Remark.} The group $\Aut (F)$ acts on the
set of c-lattices and on the datum of super
lines $\det (L:L')$, hence on the gerbe of
c-lattices. The above formula means that
$\Aut^\flat (F)$ is the super
extension defined by this action, see A6.

\medskip

If $F$ is an arbitrary Tate $R$-module then
the object $F^\infty$ is well-defined only
after certain \'etale (or Nisnevich)
localization of
$\Spec R$ (see 2.12(ii)). Since super lines
have
\'etale local nature this suffices to define
the Tate super extension of the above semigroup
hence of $\Aut (F)$.

This construction is compatible with base
change, so we  have a group object
$\GL (F)$ of $\sS_R$,
$\GL (F)(R'):= \Aut (F_{R'})$, and its Tate
super extension $\GL^\flat (F)$. Let
$\fgl (F)$ be the Lie algebra of $\GL (F)$, so
$\fgl (F)(R')$ is the Lie algebra of
endomorphisms of the Tate $R'$-module
$F_{R'}$; the Lie algebra $\fgl^\flat (F)$ of
$\GL^\flat (F)$ is a central extension of $\fgl
(F)$ by
$\sO$.

\medskip

  Let $L,M\subset F$ be,
respectively, a c- and a d-lattice in a
(special) Tate module $F$. Denote by
$\GL (F,L),$ $\GL (F,M)\subset \GL (F)$ the
parabolic subgroups of transformations
preserving $L,M$. One has the {\it  standard
splittings}
\eq{2.13.3}{s_L^c : \GL (F,L)\to \GL^\flat (F),
\quad s_M^d : \GL (F,M)\to \GL^\flat (F),}
defined as $s_L^c (g):= \det (F/L , g
)$, $s_M^d (g') := \det (M, g')$, see
(2.10.2). If $L\oplus M \iso F$ then $s^c_L (g)
=s^d_M (g)$ for $g\in \GL (F,L)\cap \GL (F,M)$.

\medskip

The central extension
$\fgl^\flat (F)$ can be described as follows.
Set $\fgl_c (F):= \Ker (\fgl (F)\to\fgl
(F^\infty ))$ and let $\fgl_d (F) \subset \fgl
(F)$ be the submodule of endomorphisms with
open kernel. Both $\fgl_c
(F),\fgl_d (F)$ are ideals in $\fgl
(F)$,\footnote{ Moreover, they are stable with
respect to the adjoint action of $\GL (F)$.}
and their sum equals
$\fgl (F)$. Set $\fgl^\ff (F):=\fgl_c
(F)\cap\fgl_d (F)$. There is a canonical
Ad-invariant trace functional
$tr:
\fgl^\ff (F)\to \sO$. Namely, to
compute
$tr (a)$ for
$a \in \fgl^\ff (F)$ you find
(after possible Zariski localization of $R$)
c-lattices $L\subset L'$ such that $a(L)=0$
and $a(F)\subset L'$; then $tr (a)$ is the
trace of the induced endomorphism of $L'/L$.

One has $\fgl (F^\infty )=\fgl
(F)/\fgl_c (F)=\fgl_d
(F)/\fgl^{\ff}(F)$. As follows from 2.10(v),
the central extension $\fgl^\flat (F^\infty )$
is the push-out of $0 \to \fgl^\ff (F) \to
\fgl_d (F)\to\fgl (F^\infty )\to 0$ by  $tr$.
Our $\fgl^\flat (F)$ is the pull-back of
$\fgl^\flat (F^\infty )$ by the projection
$\fgl (F) \to\fgl (F^\infty )$.

Here is a
c-d-symmetric description. As follows from
the above construction, $\fgl^\flat (F)$
admits canonical ad-invariant splittings
$s_c , s_d$ over the ideals $ \fgl_c (F),
\fgl_d (F)\subset\fgl (F)$ such that on
$\fgl^\ff (F)$ one has $s_d - s_c =tr$. This
structure determines $\fgl^\flat (F)$ uniquely.

\medskip

One has a canonical bijection $L\mapsto
L^\bot$ between the sets of c-lattices in $F$
and its dual $F^*$. There is a canonical
isomorphism of super lines  
\eq{2.13.4}{\det
(L:L')\iso\det (L^\bot :{L'}^\bot )}
compatible with composition,
i.e., the pre-gerbes of c-lattices in $F$ and
$F^*$ are canonically identified. To define
(2.13.4) it suffices, by Remark after
(2.13.1), to consider the case
$L\supset L'$. Here l.h.s.~is $\det
(L/L')$, r.h.s.~is inverse to $\det
((L/L')^*)$, and (2.13.4) is (2.2.3) for
$M=(L/L')^*$. The compatibility with
composition is immediate.

The isomorphism $\GL (F)\iso \GL (F^* )$,
 $g\mapsto {}^t g^{-1}$, is compatible with
the actions on the pre-gerbes of c-lattices.
So, by Remark after (2.13.2), it lifts
canonically to super
extensions
\eq{2.13.5}{\GL^{\flat } (F)\iso \GL^{\flat
}(F^* ).}

Notice that (2.13.5) changes the sign of the
$\Z$-grading to the opposite.
It interchanges
the ideals $\fgl_c ,\fgl_d \subset \fgl$ and
the  splittings $s_c$, $s_d$.
 
\medskip

{\it Remark.} The above compatibility with
duality can be extended to a larger semi-group
of  Fredholm
endomorphisms as follows. We say that a
morphism $g: N\to M$ of Tate $R$-modules is
{\it Fredholm} if both
$g^\infty : N^\infty \to M^\infty$, ${}^t
g^\infty : {M}^{*\infty}\to N^{*\infty}$ are
invertible. Equivalently, this means that (locally on
$\Spec R$) one can find c-lattices
$L\subset M, P\subset N$ such that $g(N)+L=M$,
$P\cap
\Ker g =0$ and $g^{-1}(L)\subset N$,
$g(P)\subset M$ are again c-lattices. Notice
that any pair of isomorphisms $N^\infty \iso
M^\infty$, $M^{*\infty} \iso N^{*\infty}$ can
be lifted to a Fredholm morphism $N\to M$.

A Fredholm $g$ defines a
$\Z$-graded super line $\det (M,N, g)$.
To define it choose $L$, $P$ as above such
that $L\supset g(P)$; then $\det (M,N,g):=
\det (L/g(P))/\det (g^{-1}(L)/P)$. The
arguments of Lemma-Definition in 2.7 show
that our super line does not depend on the
auxiliary choice of $L$, $P$ and actually
depends only on $(g^\infty , g^{*\infty})$. It
is compatible with composition, and there is a
canonical isomorphism $\det (M,N,g)=\det (N^*
,M^* ,{}^t g)^{-1}$ coming from (2.2.3). If
$g$ is an isomorphism then our super line is
canonically trivialized by $\det g\in
\det (M,N,g)$.

Let $g$ be a Fredholm endomorphism of $F$;
  set $\det (F,g):=\det
(F,F,g)$. There is a canonical isomorphism 
\eq{2.13.6}{\lambda_g
\cdot\lambda_{{}^t g} = \det (F,g).} Namely,
for c-lattices $L$, $P$ as above one 
has $\lambda_g =
\det (L: g^{-1}(L))$, $\lambda_{{}^t g}= \det
({P}^\bot : ({}^t g)^{-1} ({P}^\bot ))=\det
(P :g (P))$ (see (2.13.2) and (2.13.4)), and
(2.13.6) is the product of these
isomorphisms. If $g$ is invertible then the
r.h.s.~of (2.13.6) is trivialized by $\det g$,
and (2.13.6) amounts  to (2.13.5)

\subsection{ Clifford modules} Here is another
description of the Tate extension which is
inherently self-dual.

Fix a Tate module $F$.  The
Tate module
$F\oplus F^*$ carries a standard hyperbolic
symmetric bilinear form. Denote by $W=W_F$ the
same Tate module considered as a super object
placed in {\it odd} degree. The above form
considered as a bilinear form on $W$ is {\it
skew}-symmetric; we denote it by $<\, , \, >$.
So
$\sC := W\oplus R$ is a Lie super $R$-algebra
with respect to a bracket whose only non-zero
component is $<\, , \, >$.

A {\it Clifford module} for $W$ is a
$\sC$-module\footnote{in the tensor
category $\sM^s_R$ of super $R$-modules, see
2.1.} $\sQ$ such that the
$\sC$-action is continuous with respect to the
discrete topology of $\sQ$ and $1\in R = \sC^0$
acts as
$id_Q$.  The category of Clifford modules is
denoted  by
$\sM^s_\sC$. For morphisms of $R$'s one has the
obvious base change functors. It is clear that
Clifford modules are local objects with
respect to the flat topology of $\Spec R$.

For a super
$R$-module $M$ and a Clifford module $\sQ$ the
tensor product
$M\otimes
\sQ$ is a Clifford module in the obvious
way.  A
Clifford module
$\sP$ is {\it invertible} if the functor
$\sM_R^s\to\sM_\sC^s$, $M\mapsto M\otimes \sP$,
is an equivalence of categories.  Denote by
$\sP ic_\sC^s$ the groupoid of invertible
$\sC$-modules. For $\sL\in\sP
ic^s_R$, $\sP\in\sP ic_\sC^s$ one has
$\sL\otimes\sP \in \sP ic_\sC^s$, so $\sP
ic_\sC^s$ carries a canonical action of the
Picard groupoid $\sP ic_R^s$.

\medskip

Assume that $F$ is a special Tate
$R$-module. Let
$L_W
\subset W$ be a c-lattice,
$L_W^\bot
\subset W$ its $<\, , \, >$-orthogonal
complement; this is again a c-lattice. If $L_W
\subset L^\bot_W$ then $<\, , \, >$ yields a
non-degenerate form on $\bar{W}:= L_W^\bot /L_W
\in
\sV^\ff$. Let $\bar{\sC} =\bar{W} \oplus
R$ be the corresponding Clifford Lie super
algebra. In other words, $L_W^\bot \oplus R$
is the normalizer of $L_W \subset \sC$, and
$\bar{\sC}$ is the  subquotient algebra of
$\sC$. Let $\sM_{\bar{\sC}}^s$ be the category
of the corresponding Clifford
modules.\footnote{i.e., an object of
$\sM_{\bar{\sC}}^s$ is a $\bar{\sC}$-module
such that
$1\in R\subset \bar{\sC}$ acts as identity.}
Now the functor
\eq{2.14.1}{  \sM_{\sC}^s \to
\sM_{\bar{\sC}}^s ,\quad \sQ \mapsto
\bar{\sQ}=\sQ^{L_W},} is an equivalence of
categories\footnote{If
$F$ is a finitely generated projective
$R$-module then (2.14.1) is the usual
Morita equivalence. The general case reduces
to this one by 2.3(ii) (see Remark (i) in
2.13) since Clifford modules are Zariski local
objects.} (its inverse is the induction
functor). In particular, if
$L_W^\bot =L_W$ then $\bar{\sC} =R$, hence
$\sM_{\sC}^s
\iso
\sM_{R}^s$. This happens when $L_W = L
\oplus L^\bot$ where $L \subset F$ is
any c-lattice and $L^\bot := (F/L)^* \subset
F^*$ is its orthogonal complement.

Equivalences (2.14.1) commute
with functors $M\otimes \cdot$.
The groupoid  $\sP
ic_\sC^s$ is a
$\sP ic_R^s$-torsor, i.e., a super
$\sO^\times$-gerbe. Indeed, for every $L_W$ as
above such that $L_W^\bot =L_W$  (2.14.1)
identifies  $\sP ic_\sC^s$ with $\sP ic^s_R$.

\medskip

Let $\rO (W)\subset \GL (W)$ be the subgroup of
automorphisms preserving $<\, ,\, >$. It acts
on $\sC$ in the obvious way, hence $\rO (W)$
acts on the category of Clifford modules.
Explicitly, for
$g\in \rO (W)$ and $\sP \in \sM_\sC^s$ the
Clifford module $g\sP$ equals $\sP$ as a super
$R$-module, and the $\sC$-action on it is
$w,p \mapsto g^{-1}(w)p$. This
action obviously commutes with functors
$M\otimes \cdot$. Therefore $\rO (W)$ acts on
$\sP ic^s_\sC$ as on a super
$\sO^\times$-gerbe. Let
$\rO^\flat (W)$ be the  super
$\sO^\times$-extension {\it opposite} to
the one defined by this action (see A6), or,
equivalently, the super extension defined by
the action on the  opposite super
gerbe. So for
$g\in \rO (W)$ its super line in $\rO^\flat
(W)$ is
$\sP /g\sP $, $\sP \in \sP ic^s_\sC$.
Equivalently,
$\rO^\flat (W)$ consists of pairs
$(g,g^\flat )$ where $g\in \rO (W)$ and
$g^\flat$ is an automorphism of $\sP$,
considered as a plain $R$-module, such that
for $w\in W , p\in
\sP$ one has $g^\flat (wp)=g(w)g^\flat (p)$.

If $F$ is an arbitrary (not necessary
special) Tate module then we define 
$\rO^\flat (W)$ using 2.12(ii).\footnote{ As
we did to define
$\GL^\flat (F)$.}

\medskip

{\it Remarks.} (i) For $L_W$ as above set $\rO
(W,L_W ):= \rO (W)\cap \GL (W,L_W)$; there is
an obvious projection $\rO (W,L_W) \to
\rO (\bar{W})$. Then (2.14.1) identifies the
restriction of
$\rO^\flat (W)$ to $\rO (W,L_W )$ with the
pull-back of the extension $\rO^\flat
(\bar{W})$. In particular, if $L_W =L_W^\bot$
then there is a canonical splitting $s_{L_W}^c
: \rO (W,L_W )\to \rO^\flat (W)$.

(ii) The above discussion for c-lattices has an
immediate d-version. Namely, let
$M_W
\subset W$ be a d-lattice (see 2.13) such that
$<\, ,\, >$ vanishes on
$M_W$. Then $M_W^\bot \subset W$ is a
d-lattice and
$\bar{W}:= M_W^\bot /M_W \in\sV^\ff$ carries
the induced non-degenerate form; let
$\bar{\sC}$ be the corresponding Clifford Lie
super algebra. Thus $M_W^\bot \oplus R$ is the
normalizer of $M_W \subset \sC$, and $\bar{\sC}
$ is the subquotient Lie algebra. We have an
equivalence of categories
\eq{2.14.2}{\sM^s_\sC \iso \sM^s_{\bar{\sC}},
\quad \sQ \mapsto \sQ_{M_W}.} It identifies
the restriction of $\rO^\flat (W)$ to a
parabolic subgroup $\rO (W,M_W )$ $:= \rO
(W)\cap
\GL (W,M_W)$ with the pull-back of the
extension $\rO^\flat (\bar{W})$ by the
projection $\rO (W,M_W )\to \rO (\bar{W})$. In
particular, if $M_W =M_W^\bot$ then there is a
canonical splitting $s_{M_W}^d : \rO (W,M_W
)\to
\rO^\flat (W)$.

\subsection{Proposition} The embedding $\GL (F) \hra
\rO (W)$ which assigns to $g\in \GL (F)$ the
diagonal matrix $g^o$ with components $g, {}^t
g^{-1}$
lifts canonically to a morphism of  super
extensions
\eq{2.15.1}{\GL^\flat (F) \hra \rO^\flat (W).}

\begin{proof} By Remark after
(2.13.2), $\GL^\flat F)$ is the super extension
defined by the action of $\GL (F)$ on the
gerbe of c-lattices. The super extension
induced from $\rO^\flat (W)$ is defined by the
action of $\GL (F)$ on the gerbe opposite
to $\sP ic^s_\sC$. It
remains to identify these two gerbes in a way
compatible with the $\GL (F)$-actions. 

By Remark in A6 we have to assign to every
$c$-lattice $L$ an invertible Clifford
module $\sQ_L$ and define identifications
$\det (L:L')=\sQ_{L'} /\sQ_{L}$ which satisfy
the transitivity property. By Remark  after 
(2.13.1) it suffices to establish the latter
identifications for $L\supset L'$.

We define $\sQ_L$ as the Clifford module such
that $(\sQ_L )^{L_W}=R$ where $L_W := L\oplus
L^\bot$ (see (2.14.1)). 
Then $\sQ_{L'}/\sQ_{L}=(\sQ_{L'} )^{L_W}$, so
the promised identification can be rewritten
as a canonical isomorphism
\eq{2.15.2}{\det (L:L' )= \sQ^{L_W}/\sQ^{L'_W}
} valid for any $\sQ\in\sP
ic^s_\sC$. If $L\supset L'$ then the super
lines $(\sQ_{L'} )^{L_W},
\sQ_{L'}^{L'_W}$ lie in $\sQ_{L'}^P$, where
$P:= L' \oplus L^\bot \subset W$, which is a
module for the Clifford algebra $\bar{\sC}$ for
$\bar{W}= P^\bot /P = L/L'
\oplus {L'}^\bot /L^\bot$. Now $\det
(L:L')=\det (L/L')$ lies in $\bar{\sC}$, and
its action transforms
$\sQ^{L'_W}$ to $\sQ^{L_W}$. This is 
(2.15.2). The transitivity property  and
compatibility with the $\GL (F)$-action are
obvious.
\end{proof}

\medskip

{\it Remarks.} (i) We see that $\GL^\flat (F)$
acts canonically on every Clifford module
$\sP$. The action of $\fgl^\flat (F)$ can be
described explicitly as follows. For
$a^\flat \in \fgl^\flat (F)$ over $a
\in \fgl (F)$ its action on $\sP$ is compatible
with the action of $W\subset \sC$: one has
$a^\flat (wp)- w(a^\flat p) = a(w) p$ for $w\in
W$, $p\in\sP$. This condition alone does not
determine the action of $a^\flat$; an extra
normalization is needed. Since $\fgl
(F)=\fgl_c (F)+\fgl_d (F)$\footnote{See
  the discussion at the end of
2.13.} it suffices to consider the cases
$a^\flat =s_c (a)$, $a\in \fgl_c (F)$, and
$a^\flat =s_d (a)$, $a\in \fgl_d (F)$. The
action of $s_c (a)$ is uniquely determined by
property that it kills $\sP^L$ for a
sufficiently large c-lattice $L\subset F$; the
action of $s_d (a)$ is
uniquely determined by property that it kills
$\sP^{L^\bot}$ for a sufficiently small
c-lattice $L\subset F$.

(ii) Let $L,M \subset F$ be a c- and
d-lattice. Set $L_W := L\oplus L^\bot ,M_W :=
M\oplus M^\bot \subset W$; these are  c- and
  d-lattices in $W$ such that $L_W^\bot =L_W$,
$M_W^\bot =M_W$, and (2.15.1) identifies
sections $s_L^c$, $s_M^d$ of (2.13.3) with
restrictions of sections $s_{L_W}^c$,
$s_{M_W}^d$ from Remarks (i), (ii) of 2.14.

(iii) Recall that $\GL^\flat (F)$ is a
{\it $\Z$-graded} super extension. This
$\Z$-grading can be described in Clifford
terms as follows. Consider a
$\Z$-grading on $W$ with $F$ in degree 1 and
$F^*$ in degree $-1$; then $\GL (F)\subset \rO
(W)$ is the group of orthogonal automorphisms
preserving this
$\Z$-grading. One defines $\Z$-graded Clifford
modules in the obvious way. The  picture of
2.14 remains valid in the $\Z$-graded setting,
so we have a
$\Z$-graded super $\sO^\times$-gerbe $\sP
ic_\sC^\Z$ equipped with an action of $\GL (F)$
and the corresponding $\Z$-graded super
extension  of $\GL (F)$. This defines a
canonical $\Z$-grading on the pull-back of
$\rO^\flat (W)$ to $\GL (F)\subset \rO (W)$.
It is easy to see that (2.15.1) is compatible
with the
$\Z$-gradings.

(iv) The Clifford algebras for $F$ and $F^*$
coincide, hence the corresponding groups $\rO
(W)$ and their super extensions $\rO^\flat
(W)$ are equal. On the subgroups $\GL (F)$,
$\GL (F^* )$ this identification is the
standard isomorphism $g\mapsto{}^t g^{-1}$;
 its lifting to the Tate super
extensions via (2.15.1) is (2.13.5).

\subsection{Scalar products and
super $\mu_2$-torsors.} We suggest
  reading this subsection simultaneously with
3.6.

Let
$\Phi =\Phi (F)\in \sS_R
$ be the object of all symmetric
non-degenerate bilinear forms on
$F$, or, equivalently, that of symmetric
isomorphisms
$\phi :F\iso F^*$ (the form corresponding to
$\phi$ is
$(a ,b )_\phi := (\phi a, b)$).

For $\phi \in\Phi$  let $\phi^o \in \rO (W)$ be
the anti-diagonal matrix with components $\phi
,\phi^{-1}$. Let $\lambda_{\phi}$
be the super line of $\phi^o$ in $\rO^\flat
(W)$. One has
$\phi^{o2}=1$, so the corresponding
identification
$\lambda_{\phi}^{\otimes 2}\iso
\sO_\phi$ yields a super $\mu_2$-torsor
$\mu_\phi$ (see 2.1). We have defined a
canonical super $\mu_2$-torsor $\mu_\Phi$ on
$\Phi$.

For  $\phi ,\phi' \in \Phi (F)$
set
$g_{\phi' ,\phi}:=
{\phi'}^{-1}\phi
\in \GL (F)$. One has
  ${\phi'}^{o}=g_{\phi' ,\phi}^o \phi^o
\in \rO (W)$. By
  (2.15.1) this provides a canonical
isomorphism of super lines
\eq{2.16.1}{\lambda_{g_{\phi'
,\phi}}\cdot\lambda_{\phi}\iso \lambda_{\phi'}.
}
  Notice that $g_{\phi''
,\phi}=g_{\phi'' ,\phi'}g_{\phi' ,\phi}$ and
identifications (2.16.1) are
transitive.

One can rephrase this as follows.
Let $\sG_\Phi$
be the simply transitive groupoid on $\Phi$,
so for every $\phi ,\phi' \in \Phi$ there
is a single arrow $\phi \to \phi'$ in
$\sG_\Phi$. We have a homomorphism of groupoids
$g: \sG_\Phi
\to \GL (F)$ which sends the above arrow to
$g_{\phi' ,\phi}$. Let
$\sG_\Phi^{\flat}$ be the $g$-pull-back
of the Tate super extension. Now (2.16.1) is a
$\sG_\Phi^{\flat}$-action on
$\lambda_\Phi$, i.e., a
$\lambda_\Phi$-splitting of
$\sG_\Phi^{\flat}$ in the terminology of
A2.

Therefore (see A2)
$\sG_\Phi^{\flat}$ comes from a
$\mu_\Phi$-split super $\mu_2$-extension
$\sG_\Phi^{\mu}$ of
$\sG_\Phi$. If
$1/2
\in R$ then any $\mu_2$-extension is
\'etale, hence it is canonically
trivialized over the formal completion
$\hat{\sG}_\Phi$ of
$\sG_\Phi$.\footnote{For
$\phi ,\phi' \in \Phi (R')$ the arrow
$\phi \to \phi'$ lies in
$\hat{\sG}_\Phi \subset
\sG_\Phi$ if and only $\phi$ equals $\phi'$
modulo  some nilpotent ideal
$I\subset R'$.} The corresponding
trivialization of
$\sG_\Phi^{\flat}$  over
$\hat{\sG}_\Phi$ is a {\it  canonical
  formal rigidification} of
$\sG_\Phi^{\flat}$.

\medskip
{\it Remarks.} (i) The above constructions are
compatible with direct sums: if $F=\oplus
F_\alpha$, $\phi =\oplus \phi_\alpha$, then
$\mu_\phi =\otimes \mu_{\phi_\alpha }$, and
duality:
  $\mu_{\phi^{-1}}=\mu_\phi $.

(ii)  Fix any $\phi \in \Phi (F)$. Then the
morphism
$g_\phi : \Phi (F) \to \GL (F)$, $\phi' \mapsto
g_{\phi' ,\phi}$, is injective, and its image
consists of all
$g\in \GL (F)$ which are self-adjoint with
respect to $(\, ,\,)_\phi$. The involution
$\fs_\phi :=$Ad$_{\phi^o}$ preserves
$\GL (F),\Phi \subset \rO (W)$, so it induces
involutions of $\GL (F)$ and $\Phi$ which we
denote again by $\fs_\phi$. For
$g\in \GL (F)$, $\phi' \in \Phi$ one has
$\fs_\phi (g )= \phi^{-1}({}^t
g)^{-1}\phi $ (which is the adjunction with
respect to $(\, ,\, )$) and $\fs_\phi
(\phi' )=\phi\phi^{'-1}\phi $. Thus
$g_\phi (\fs_\phi (\phi' ) )=(g_\phi (\phi'
))^{-1}$. The action of $\fs_\phi$ on
$\lambda_\Phi$ is compatible with the
trivialization $\lambda_\Phi^{\otimes 2}\iso
\sO_\phi$, so $\fs_\phi$ acts on $\mu_\Phi$.

(iii) The action of $\fs_\phi$ on $\fgl (F)$
interchanges the ideals $\fgl_c (F),\fgl_d
(F)$ $\subset \fgl (F)$ and the sections $s_c
,s_d$ (see (2.13.5) and Remark (iv) in 2.15).
Our
$g_\phi$ identifies the tangent space to $\phi
\in\Phi$ with the submodule
$\fgl (F)^\phi \subset \fgl (F)$ of
$(\, ,\,)_\phi$-self-adjoint operators. If
$1/2\in R$ then the canonical formal
rigidification of $\sG_\phi^\flat$
yields a splitting $\nabla_\phi^\flat :\fgl
(F)^\phi
\to
\fgl^\flat (F)$ which is
$\fs_\phi$-invariant; it is
uniquely defined by this property.

(iv) Let $L,L_W ,M, M_W$ be as in Remark (ii)
of 2.15. Set
\eq{2.16.2}{\Phi_L :=
\{\phi\in
\Phi :\phi (L)=L^\bot \},\,\,  \Phi_M := \{\phi
\in\Phi :\phi (M)=M^\bot \}.} Since $\Phi_L^o
\subset\rO (F,L_W)$, $\Phi_M^o \subset \rO
(W,M_W )$ the sections $s_{L_W}^c$,
$s_{M_W}^d$ from Remarks (i), (ii) of 2.14
trivialize the restrictions of our
$\mu_2$-torsor $\mu_\Phi$ to
$\Phi_{L_W}$ and $\Phi_{M_W}$. We denote these
trivializations also by $s_L^c$,
$s_{M}^d$.  Notice that for $\phi ,\phi'
\in\Phi_{L}$ the action (2.16.1) identifies
$s_L^c (g_{\phi',\phi})s_L^c (\phi  )$ with
$s_L^c (\phi' )$ (see  Remark (ii) of 2.15).
Same for
$L$ replaced by $M$.

\subsection{} Let us describe $\mu_\phi$
explicitly assuming that
$R$ is a field.

(i) Case $p(\mu_\phi ) =0$.\footnote{Here
$p$ is the parity of our super torsor, see
2.1.} Then one can find a c-lattice
$L$ such that the
$\phi$-orthogonal complement $L^\bot_\phi$
equals
$L$, i.e., $\phi \in \Phi_L$. As above, it
  yields the
trivialization $s_L^c (\phi )$ of $\mu_\phi$.
If
$F\neq 0$ and we live in characteristic
$\neq 2$  then the Grassmannian of such
$L$'s has 2 components. We leave it to the
reader to check that our trivialization
changes sign as we switch the component. Thus
$\mu_\phi
$ is the set of components considered as a
$\mu_2$-torsor.

(ii) Case $p(\mu_\phi )=1 $. Then one
can find a c-lattice $L$ such that $L^\bot_\phi
\supset L$ and $\bar{L} :=L^\bot_\phi /L$ has
dimension 1. We have the induced $\bar{\phi}
:\bar{L} \iso \bar{L}^*$.  Let us define
a canonical identification
\eq{2.17.1}{\mu_\phi  =\{
\bar{l}\in\bar{L}:(\bar{l},\bar{l})_{\bar{\phi}}
=1\}.} Indeed,
$\phi^o$ preserves $L_W := L\oplus \phi (L)$
and
$\bar{W}:= L_W^\bot /L_W = \bar{L}\oplus
\bar{L}^*$, so, by Remark (i) in 2.14, we have
$\mu_\phi =\mu_{\bar{\phi}}$. Take any
invertible Clifford module
$\bar{\sP}$ for $\bar{W}$. An
element
$\sigma
\in\mu_{\bar{\phi}}$ is an automorphism
  of $\bar{\sP}$ as a
$k$-vector space such that
$\sigma^2 =1$ and
$\sigma (wp) =\bar{\phi}^o
(w)\sigma (p)$ for $w\in
\bar{W}$, $p\in\sP$. The restriction of such
$\sigma$ to $
\bar{\sP}^{\bar{L}^*}\subset \bar{\sP}$ is
multiplication by
$\bar{l}\in \bar{L}$ such that
$(\bar{l},\bar{l})_{\bar{\phi}}
=1$.  This provides (2.17.1).

\bigskip

\centerline{\bf Appendix }

\medskip

  A1.\footnote{See SGA
4 XVIII 1.4.} A {\it Picard groupoid}
$\sP$ is a symmetric monoidal
category\footnote{We always assume that $\sP$
is equivalent to a small category.} such that
every object in
$\sP$ is invertible, as well as every
morphism. We denote the operation in $\sP$ by
$\cdot$  and the unit object by
$1_\sP$. The commutative group of 
isomorphism classes of objects in $\sP$ is
denoted by $\pi_0 (\sP )$. Set $\pi_1 (\sP ):= 
\Aut (1_\sP )$; this is a commutative group,
and for every $P\in\sP$ the identification
$1\cdot P =P$ yields a canonical isomorphism
$\pi_1 (\sP )\iso \Aut P$, $a\mapsto a\cdot
id_P$.

For $P\in\sP$ the corresponding element of
$\pi_0 (\sP )$ is denoted by $|P|$. For a finite family $\{ F_\alpha \}$ of
objects of $\sP$ we usually denote their
product by $\otimes F_\alpha$.

{\it Remark.} For any commutative group $A$
the category $A$-$tors$ of $A$-torsors is a
Picard groupoid in the obvious way. For every
Picard groupoid $\sP$ its Picard
subgroupoid $\sP^0 $ of objects isomorphic to
$1_\sP$ is canonically equivalent to
$\sA_\sP$-$tors$ via $\sP^0 \iso
\pi_1 (\sP )$-$tors$, $P\mapsto \Hom (1_\sP ,P
)$.

For two Picard groupoids $\sP$, $\sP'$ a
{\it morphism} $\phi :\sP \to\sP'$ is the same
as a symmetric monoidal functor. All morphisms
form a Picard groupoid $\Hom (\sP ,\sP' )$:
namely, the product of two morphisms
$\phi,\psi$ is
$(\phi\cdot\psi )(P)=\phi (P)\cdot \psi (P)$.

\medskip

 For $P\in\sP$ its {\it
inverse} $P^{-1}$ is an object of $\sP$
together with an identification $e: P\cdot
P^{-1} \iso 1_{\sP}$; the inverse is
determined uniquely up to a unique isomorphism.
There are {\it two} natural identifications
$P\iso (P^{-1})^{-1}$ defined by pairings
$e^c , e^m : P^{-1}\cdot P \iso 1_{\sP}$ where 
$e^c$ is defined using commutativity, and $e^m$
is determined by the condition that  $e\cdot
id_P = id_P \cdot e^m : P\cdot P^{-1} \cdot P
\iso P$ (or, equivalently, that $e^m \cdot
id_{P^{-1}}=id_{P^{-1}} \cdot e : P^{-1}\cdot
P\cdot P^{-1}
\iso P^{-1}$). Notice that the definition
of $e^m$ does {\it not} use commutativity of
$\cdot$, and 
$e^m /e^c =\alpha
$ where $\alpha =\alpha (P)\in \pi_1 (\sP )$
is the action of the commutativity symmetry on
$P\cdot P$. {\it Unless stated
explicitly otherwise, we use identification
$e^c$.}

Using commutativity, one extends naturally the
map $P\mapsto P^{-1}$ to an
auto-equivalence of the Picard groupoid
$\sP$.

For
$P,Q\in\sP$ we write
$P/Q := P\cdot Q^{-1}$.

\medskip

  A2. For a category $\Gamma$, a {\it
$\sP$-extension} $\Gamma^\flat$ is a rule that
assigns to every arrow $\gamma$ in $ \Gamma$ an
object
$P^\flat_\gamma
\in\sP$ and to every pair of composable
arrows $\gamma ,\gamma'$ a composition
isomorphism
$c_{\gamma , \gamma'}: P^\flat_\gamma \cdot
P^\flat_{\gamma'}\to  P^\flat_{\gamma\gamma'}$.
We demand that
$c$ be associative in the obvious sense.

Assume we have a $\sP$-bundle $\sL$ on the
set of objects of $\Gamma$ which is
a rule which assigns to every object $x$ in
$\Gamma$ an object $\sL_x $ of $\sP$. An
{\it action} of $\Gamma^\flat$ on $\sL$ is
a rule that assigns to every arrow $\gamma :
x\to x'$ an isomorphism $P^\flat_\gamma
\cdot \sL_x \iso \sL_{x'}$ which satisfies an
obvious transitivity property.

Notice that for a given $\sL$
there is a unique (up to a unique isomorphism)
$\sP$-extension $\Gamma^\flat_\sL$ of
$\Gamma$ acting on $\sL$. Indeed, the action
amounts to a datum of
isomorphisms
$P_{\sL\gamma}^\flat :=
\xi_{x'}/\xi_x$ which identifies the
composition in $\Gamma^\flat$ with the obvious
product. We call
$\Gamma^\flat_\sL$ the {\it 
$\sL$-split}
$\sP$-extension. So for an
arbitrary $\sP$-extension $\Gamma^\flat$ of
$\Gamma$ we also refer to its action on $\sL$
as an {\it $\sL$-splitting} of $\Gamma^\flat$.

\medskip

{\it Remarks.} (i) For every object $x$ in
$\Gamma$ there is a canonical identification
$e_x : P^\flat_{id_x} \iso 1_\sP$ defined by
$c_{id_x ,id_x}: P^\flat_{id_x} \cdot
P^\flat_{id_x}
\iso P^\flat_{id_x}$. This isomorphism
identifies every
$c_{id_x ,\gamma'}$ and
$c_{\gamma,id_x}$ with the 
canonical isomorphisms
$1_\sP \cdot P^\flat_{\gamma'} \iso
P^\flat_{\gamma'}$,  $P^\flat_\gamma \cdot
1_\sP
\iso P^\flat_\gamma$. For every invertible
arrow
$\gamma :x \to x'$ in $\Gamma$ there is a
canonical identification
$P^\flat_{\gamma^{-1}} \iso (P^\flat_\gamma
)^{-1}$ defined by the isomorphism $e_\gamma
:=  e_{x'} c_{\gamma ,\gamma^{-1}}  
:P^\flat_\gamma
\cdot P^\flat_{\gamma^{-1}} \iso 1_\sP$.
Notice that $e_{\gamma^{-1}} = e_\gamma^m$.

(ii) A $\sP$-extension $\Gamma^\flat$
yields a homomorphism\footnote{For a
group $A$ a homomorphism
$a:\Gamma \to A$ assigns to every arrow
$\gamma$ an element $a(\gamma )\in A$ so that
$a(\gamma
\gamma')=a(\gamma )a(\gamma' )$ for every
composable $\gamma$, $\gamma'$.} $\Gamma \to
\pi_0 (\sP )$,
$\gamma \mapsto |P^\flat_\gamma |$. If
$\sP$ is discrete, i.e.,
$\sA_\sP =1$, then $\sP$-extensions are the
same as homomorphisms $\Gamma \to \pi_0 (\sP
)$.

\medskip

A3. Let $\Gamma^\flat$, $\Gamma^{\flat'}$ be
two
$\sP$-extensions of $\Gamma$. A {\it morphism}
$\theta :\Gamma^\flat \to\Gamma^{\flat'}$ is a
system of morphisms $\theta_\gamma :
P_\gamma^\flat \to P^{\flat '}_\gamma$
compatible with composition; it is clear how
to compose morphisms of $\sP$-extensions. The
{\it Baer product} of
$\Gamma^\flat$,
$\Gamma^{\flat'}$ is a $\sP$-extension $\gamma
\mapsto P_\gamma^\flat \cdot P^{\flat
'}_\gamma$ with the obvious composition rule;
the Baer product is associative and
commutative in the obvious way. Therefore for
small $\Gamma$ its $\sP$-extensions form a
Picard groupoid which we denote by $\sE xt
(\Gamma ,\sP )$. Its unit object $1_\sE$ is
 the {\it  trivial split extension} $\gamma
\mapsto 1_\sP$. One has $\pi_1 (\sE xt
(\Gamma ,\sP ))=\Hom (\Gamma
,\pi_1 (\sP ))$.

For a functor $\phi : \Gamma'
\to \Gamma$ and a $\sP$-extension
$\Gamma^\flat$ it is clear what is the
pull-back of $\Gamma^\flat$ by $\phi$. Thus
for small $\Gamma$, $\Gamma'$ we have a
morphism of Picard groupoids $\phi^* : \sE xt
(\Gamma ,\sP )\to\sE xt
(\Gamma' ,\sP )$.

{\it Remarks.} (i) The notion of
$\sP$-extension of
$\Gamma$ depends on the {\it isomorphism}
class of
$\Gamma$, and {\it not} on its equivalence
class: if
$\phi :\Gamma \to\Gamma'$  is an
equivalence  then $\phi^*$ need
not be an equivalence.\footnote{Consider the
case when $\Gamma$ is a trivial groupoid and
$\Gamma'$ is an equivalent groupoid
with 2 objects. }

(ii) The notion of $\sP$-extension depends
only on {\it monoidal} structure of $\sP$: the
commutativity constraint is irrelevant.
However, the commutativity constraint is used
in the definition of the Baer product.

(iii) Let $\Gamma^\flat$ be a $\sP$-extension
of a {\it groupoid} $\Gamma$. Then for  $x\in
\Gamma$ we have a $\sP$-extension
$\Aut^\flat (x)$ of the group $\Aut (x)$. It
depends {\it functorially} on $x$: for every
$\gamma : x\to x'$ the identification
$Ad_\gamma :\Aut (x) \iso\Aut (x')$
lifts  canonically to an isomorphism of
$\sP$-extensions
$\Ad_\gamma^\flat :\Aut^\flat (x)
\iso\Aut^\flat (x')$. Namely, for $g\in \Aut
(x)$ the corresponding identification
 $\Ad_\gamma^\flat :P^\flat_g \iso
P^\flat_{\Ad_\gamma (g)}$ is the composition
$$P^\flat_g \iso P^\flat_g \cdot
(P^\flat_\gamma
\cdot P_{\gamma^{-1}}^{\flat})\iso
P^\flat_\gamma
\cdot P^\flat_g
\cdot P_{\gamma^{-1}}^{\flat} \iso
P^\flat_{Ad_\gamma (g)}.$$ Here the first
arrow is $id_{P^\flat_g}\cdot
(e_\gamma )^{-1}$, the second one is the
commutativity constraint, the third  is the
composition map.\footnote{The compatibility of
$\Ad_\gamma^\flat$ is with product
(composition) follows from the fact that the
composition map $P_{\gamma^{-1}}^\flat \cdot
P_\gamma^\flat \to 1_\sP$ is equal to
$e_\gamma^m$, see Remark (i) in A2. } For
composable
$\gamma,
\gamma'$ one has $\Ad_{\gamma \gamma'}^\flat
=\Ad_\gamma^\flat \Ad_{\gamma'}^\flat$.

\medskip

A4.  Let $A$ be a
commutative group. A {\it central
$A$-extension} of
$\Gamma$ is a category $\tilde{\Gamma}$
together with a functor $\pi
:\tilde{\Gamma}\to\Gamma$ and a homomorphism
from $A$ to the automorphisms of the identity
functor of
$\tilde{\Gamma}$. We demand that $\pi$ is
bijective on objects, for every $x,x'$ the
projection
$\pi: \Hom_{\tilde{\Gamma}}(x,x')\to\Hom_\Gamma
(x,x')$ is surjective, and the action of $A$
$a(\tilde{\gamma})=a_{x'}\tilde{\gamma}=\tilde{\gamma}a_x$
is simply transitive along the fibers of this
projection. If $\Gamma$ is a group (i.e., a
groupoid with single object), then its central
$A$-extensions  coincide with central
extensions of $\Gamma$ by $A$.

Any $\tilde{\Gamma}$ as above yields an
$A$-$tors$-extension $\Gamma^\flat$ of
$\Gamma$: namely, $P^\flat_\gamma$ is the
$A$-torsor $\pi^{-1}(\gamma )$, the
composition $c$ is the composition of arrows
in $\tilde{\Gamma}$.

This way we see that $A$-$tors$-extensions are
the same as central
$A$-extensions.

\medskip

  A5. Suppose that $\Gamma$ is a group; set
$A:= \pi_1 (\sP )$. As we have seen in Remark
(iii) of A3, the adjoint action of
$\Gamma$ on itself lifts canonically to a
$\Gamma$-action $\Ad^\flat$ on any
$\sP$-extension $\Gamma^\flat$. If
$\gamma ,\gamma' \in\Gamma$ commute then
$\Ad_\gamma^\flat : P^\flat_{\gamma'} \iso
P^\flat_{\gamma'}$ is multiplication by an
element of $\sA$ which we denote by $\{
\gamma ,\gamma' \}^\flat$. 
Equivalently,
$\{ \gamma ,\gamma '\}^\flat \in A$ is the
composition 
$$P^\flat_\gamma \cdot
P^\flat_{\gamma'}
\buildrel{c_{\gamma,\gamma'}}\over\rightarrow
P^\flat_{\gamma\gamma'} =P^\flat_{\gamma'
\gamma}
\buildrel{c_{\gamma',\gamma}^{-1}}\over\rightarrow
P^\flat_{\gamma'}
\cdot P^\flat_{\gamma}\iso P^\flat_\gamma \cdot
P^\flat_{\gamma'}$$ where the last arrow is the
commutativity constraint. If $\Gamma$ is
commutative then $\{ \, ,\,
\}^\flat :\Gamma\times \Gamma \to \sA$ is a
bimultiplicative skew-symmetric pairing. If
$\sP =A$-$tors$ then $\{ \, ,\, \}^\flat$ is
the usual commutator pairing for the
corresponding central
$A$-extension.

\medskip

A6. A {\it $\sP$-action} on a category $\sC$
is a functor $\cdot :\sP \times \sC \to \sC$
equipped with an associativity constraint $P
\cdot (P' \cdot C )\iso (P\cdot P')\cdot C$
and $1_\sP \cdot C \iso C$ which satisfies the
obvious compatibilities. We say that $\sC$ is
a {\it $\sP$-torsor} if $\sC$ is non-empty and
for every $C\in \sC$ the functor $\sP \to
\sC$, $P\mapsto P\cdot C$, is an equivalence
of categories. For a $\sP$-torsor $\sC$ one has
a canonical functor $\sC \times \sC \to \sP$,
$C,C' \mapsto C/C'$, together with natural
identifications
$(C/C')\cdot C' = C$; such functor is unique.

{\it Remark.} Let $\sC$ be a $\sP$-torsor. We
have then a $\sP$-extension of the simply
transitive groupoid on its set of objects: for 
 $C', C \in\sC$ the object $P_{C,C'}^\flat
\in\sP$ corresponding to the (unique) arrow
$C' \to C$ is $ C/C'$, the composition is
obvious. Conversely, suppose we have a set
$K$ together with a $\sP$-extension $P^\flat_K$
of the simply transitive groupoid on
$K$.\footnote{Which is a rule that assigns to 
$k, k' \in K$ an object $ P^\flat_{k,k'}\in
\sP$ and composition maps $P^\flat_{k,k'}\cdot
P^\flat_{k' ,k''} \to P^\flat_{k,k''}$ which
satisfy the associativity property.} Then
there is a $\sP$-torsor $\sC$ together with a
map from $K$ to the objects of $\sC$ and its
lifting to a morphism of $\sP$-extensions.
Such datum is unique in the obvious sense. We
call $\sC$ the $\sP$-torsor generated by {\it
pre
$\sP$-torsor} $(K,P^\flat_K )$.

\medskip

Let $\Gamma$ be a group. Assume that it acts
on a category $\sC$, i.e., we have a rule that
assigns to every $\gamma \in \Gamma$ a
functor $\gamma :\sC \to \sC$, $C\mapsto \gamma
C$, together with natural isomorphisms
$\gamma (\gamma' C)=(\gamma \gamma')C$,
$1C=C$; this datum should satisfy the usual
compatibilities. If $\sC$ is also equipped
with a $\sP$-action then we say that the
$\Gamma$ and $\sP$-actions {\it commute} if we
are given a system of natural isomorphisms
$\gamma (P\cdot C)=P\cdot (\gamma C)$ that
  satisfy the obvious
compatibilities.

\medskip

For example, assume that $\sP$ acts on $\sC$
and we have a $\sP$-extension $\Gamma^\flat$ of
$\Gamma$. Then $\gamma C:= P_\gamma^\flat
\cdot C$ is a
$\Gamma$-action on $\sC$ which commutes with
the $\sP$-action in the obvious way.

If $\sC$ is a $\sP$-torsor then every
$\Gamma$-action on $\sC$ commuting with the
$\sP$-action (we say then that $\Gamma$ acts
on $\sC$ {\it as on a $\sP$-torsor}) arises as
above. Precisely, there is a unique, up to a
unique isomorphism,
$\Gamma^\flat$ together with an identification
of the corresponding action on $\sC$ with our
action which is compatible with the
constraints. Namely, one has $P_\gamma^\flat =
(\gamma C )/C$ (for various $C\in\sC$ these
objects are canonically identified in the
obvious way) and the composition law is
$P_\gamma^\flat \cdot P_{\gamma'}^\flat
=(\gamma' (\gamma C)/\gamma C)\cdot (\gamma C
/C)=\gamma'\gamma C/C
=P^\flat_{\gamma\gamma'}$. We call
$\Gamma^\flat$ the $\sP$-extension
defined by the action of $\Gamma$ on $\sC$.

For example, suppose we have a pre
$\sP$-torsor $(K, P_K^\flat )$ as in Remark
above equipped with a $\Gamma$-action (so
$\Gamma$ acts on the
simple transitive groupoid
on $K$,\footnote{In all
examples we consider $\Gamma$ actually acts on
$K$ as on a mere set.} and this action is
lifted to
$P_K^\flat$). Then
$\Gamma$ acts on the
$\sP$-torsor $\sC$ generated by $(K,P^\flat_K
)$, which yields a
$\sP$-extension $\Gamma^\flat$.

\section{The Heisenberg group and its
cousins}

The purpose of this section is to study 
 various Heisenberg central
extensions which arise by pullback from the Tate central extension in
2.10. In 3.5, the superline 
$\lambda_{\omega (F)^\times}$ on
$\omega (F)^\times$ is constructed and the
structure group is reduced to $\mu_2$. In
particular, this line has a canonical
connection. For $V$ a projective
$F$-module of rank $n$, a connection on $\det_FV$ leads in 3.8 to an
identification of formal groupschemes $\widehat{F}^{\times n\flat}
\cong\widehat{F}^{\times \flat}_{(V)}$. For $\sE$ any superline bundle on
$\omega (F)^\times$ with an action of
$\widehat{F}^{\times \flat}_{(V)}$, this
isomorphism enables us to transfer the
connection from
$\lambda_{\omega (F)^\times}$ to $\sE$. 

\subsection{The setting} Let $R$
be a commutative ring and $F_R$ a
commutative topological
$R$-algebra which is a Tate $R$-module (see
2.11). Let $F_R^\times$ be the multiplicative
group of $F_R$.  By base
change we have a group valued functor
$F^\times \in \sS_R$, $F^\times (R'):=
F_{R'}^\times$. The Lie algebra of $F^\times$
is  $F$: its
$R'$-points is the additive group of
$F_{R'}$.

Similarly, let $\Aut (F_R )$ be the group of
continuous $R$-automorphisms of the algebra
$F$; we have a group functor $\Aut
(F)\in\sS_R$,
$\Aut (F)(R'):=\Aut (F_{R'})$, which acts on
$F^\times$. The
$R$-points of the Lie algebra of $\Aut (F)$ is
the Lie algebra $\Theta (F_R )$ of continuous
$R$-derivations of $F_R$.

Let $F^\tau_R$ be our $F_R$ considered as
a Tate $R$-module, so we have a group functor
$\GL (F^\tau )$, $R'\mapsto \GL (F^\tau_{R'})$
(see 2.13).  One has the obvious embeddings
\eq{3.1.1}{F^\times ,\,\Aut (F)
\hra \GL (F^\tau ).} They
identify the action of
$\Aut (F)$ on $F^\times$ with
the adjoint action of $\GL (F^\tau )$.

The pull-back of the Tate super extension
of $\GL (F^\tau )$ from 2.13 to $F^\times$ is
the {\it Heisenberg super extension}
$F^{\times\flat}$ of $F^\times$. Its Lie
algebra is a central extension
$F^\flat$ of $F$ by $\sO$. One  deduces an
explicit description of $F^\flat$ from
2.10(iv).
Notice that the action of $\Aut (F)$ on
$F^\times$ lifts canonically to the (adjoint)
action on
$F^{\times\flat}$.

The Heisenberg super extension defines the
commutator pairing\footnote{See 2.1, A5.}
\eq{3.1.2}{ \{\, ,\,\}^\flat :
F^{\times}\times F^{\times} \to
\sO^\times }  and its Lie algebra version
\eq{3.1.3}{ [\, ,\, ]^\flat : F\times F \to
\sO .} So the adjoint action of $f\in
F^\times$ on $F^{\times\flat}$ and $F^\flat$
is\footnote{The second formula comes from the
description 2.10(v) of $F^\flat$: if
$\tilde{a}$ is an operator with open kernel
which represents $a^\flat$ then $\Ad_f
(a^\flat )$ is represented by $f
\tilde{a}f^{-1}$, so $\Ad_f
(a^\flat )-a^\flat =tr ( f
\tilde{a}f^{-1} -\tilde{a})=tr([f,
\tilde{a}f^{-1}])=[f, af^{-1}]^\flat=[f, f^{-1}a]^\flat$.}
\eq{3.1.4}{\Ad_f (g^\flat )=\{f ,g\}^\flat
g^\flat , \quad \Ad_f (a^\flat ) =
[f,f^{-1}a]^\flat + a^\flat .}

Consider the Tate dual $F^*_R :=
\Hom_R^{cont}(F_R^\tau ,R)$ as an
$F_R$-module. As follows from (2.10.3) the
$R$-linear morphism \eq{3.1.5}{ d :F \to F^*
,\quad da (b):= [a,b]^\flat ,} is a
continuous derivation. It yields a canonical
$F_R$-linear morphism
$\omega_{F_R } \to F^*_R$ where
$\omega_{F_R}=\omega_{F_R /R}$ is the module of
continuous differentials relative to $R$. Let
$\Res^\flat  :
\omega_{F_R } \to R$ be its composition with
evaluation at $1\in F_R$, i.e., $
\Res^\flat (bda):=[a,b]^\flat$. Thus the
adjoint action of $F^\times$ on $F^\flat$ is
given by a cocycle \eq{3.1.6}{F^\times \times
F\to\sO, \quad f,a \mapsto \Res^\flat (ad\log
f).}

{\it Exercises.} (i) Let
$A\to R$ be a morphism of commutative algebras
such that $R$ is finite and flat over $A$.
Then $F_R$ is also a Tate $A$-module, so we
have the two pairings $ \{\, ,\,\}^\flat_R  :
F^{\times}_R \times F^{\times}_R \to R^\times$
and $\{\, ,\,\}^\flat_A :
F^{\times}_R \times F^{\times}_R \to
A^\times$. Show that $\{\, ,\,\}^\flat_A
=Nm_{R/A}\{\, ,\,\}^\flat_R$.\footnote{Here
$Nm_{R/A}:A^\times \to R^\times$ is the norm
map.}

(ii) Assume that $F$ is a product of
finitely many algebras $ F_\alpha$ as above.
Show that for $f=(f_\alpha ),$ $
g=(g_\alpha )\in F^\times =\Pi
F^\times_\alpha$ one has $\{ f,g\}^\flat_F
=\Pi \{ f_\alpha ,g_\alpha
\}^\flat_{F_\alpha}$.

\subsection{The case of $F_R =R((t))$} An
important example of the above $F$ is
the algebra
$F_R =R((t))$ of Laurent power
series equipped with the usual
topology.\footnote{A base of neighbourhood of 0
is $t^n R[[t]]$,
$n\ge 0$.} Let us discuss it in more
detail.

The language of ind-schemes is convenient
here. For us, an {\it ind-scheme} is a
functor $X$ on the category of
commutative rings which can be
presented as a filtering inductive
limit of functors $R\mapsto X_\alpha (R)$ where
$X_\alpha$ are schemes and
transition maps $X_\alpha \to
X_\beta$ are closed embeddings.  We say that $X$ is {\it
ind-affine} if all $X_\alpha$ are affine
schemes, $X$ is {\it formally smooth} if  the
usual Grothendieck property is satisfied;
$X$ is {\it reduced} if one can choose
$X_\alpha$ to be reduced schemes. To every $X$
there corresponds the reduced ind-scheme
$X_\red :=\rightlim (X_{\alpha})_{\red}$.

  For a commutative ring $R$  set
  $F_R =R((t))$. It is well-known
that the multiplicative group $F_R^\times$
is a direct product of the following four
subgroups:
(i) the subgroup $\Z_R$ of elements $t^d$;
(ii) $R^\times$;
(iii) the subgroup of elements
$1+ a_1 t +a_2 t^2 +..$, $a_i \in R$; and
(iv) the subgroup of elements $1+ b_1 t^{-1}
+b_2 t^{-2} +..$ where $b_i$ are {\it
nilpotent} elements of $R$ and almost all $b_i
=0$.

The functor $F^\times : R\mapsto F^\times (R):=
F_R^\times$ is a commutative group
ind-scheme which is ind-affine and
formally smooth. The subgroup (i) is
represented by the discrete group $\Z$, (ii) by
$\G_m$,   (iii) by the group scheme $\W$ of
(big) Witt vectors, and (iv) by the formal
completion $\hat{\W}$ of
$\W$ (via the change of variable $t\mapsto
t^{-1}$). Therefore \eq{3.2.1}{F^\times =
\Z\times\G_m
\times \W\times
\hat{\W}.}

One has $F^\times_\red = \Z\times\G_m \times
  \W$. The projection $ F^\times \to
\Z$ has connected kernel; for $d\in\Z$  let
$F^{\times d}\subset F^\times$ be the preimage
of
$d$.

  The functor $\Aut
(F) : R\mapsto \Aut (F)(R):=\Aut (F_R )$ is a
group ind-scheme.  The map $\Aut (F)\to
F^{\times 1} , \phi \mapsto \phi (t)$, is an
isomorphism of ind-schemes.

The next proposition shows that (3.1.2)
coincides in the present situation  with the
inverse to a ``parametric" version  of the
tame symbol introduced by Contou-Carr\`ere
\cite{C} (see also \cite{D3} 2.9); the usual
tame symbol is its restriction to
$F^\times_{\red}$.

\subsection{Proposition} (i) $\{\,
,\,\}^\flat$ is $\Aut
(F)$-invariant.

(ii)
  With respect to decomposition (3.2.1)  the
non-zero components

of $\{\,
,\,\}^\flat$ are:

(a) the  $\Z \times \Z$-component: one has
$\{ t^m ,t^n \}^\flat =(-1)^{mn}$.

(b) the $\G_m \times\Z$-component (and its
transpose):  $\{ r, t^m \}^\flat = r^{-m}$.

(c) the $\W\times\hat{\W}$-component (and its
transpose): for $f(t)\in \W$ and

  $g(t)\in
\hat{\W}$ one has $\{ f(t),g(t^{-1})\}^\flat =
r(f*g)$. Here $*$ is the product

  of
Witt vectors\footnote{Recall (see
e.g.~\cite{Mu} Lecture 26) that
$\W$ is naturally a commutative ring scheme;
the ring structure is uniquely characterized by
property that every map $w_i : \W (R) \to R$, 
$\Sigma w_i t^{i} := -t\partial_t \log f$,
is a ring homomorphism. } which sends
$\W\times\hat{\W}$ to
$\hat{\W}$, and $r: \hat{\W}\to\G_m$ is 

the morphism $h(t)\mapsto
h(1)^{-1}$.\footnote{Here $\hat{\W} (R)$ is
identified with units in $1+tR[t]$, i.e., with
polynomials of the form $1+a_1 t +...$ where
$a_i$ are nilpotent.}

(iii) The morphism $\Res^\flat$  is
the usual residue $\Res$ at $t=0$, so the
canonical morphism $\omega_F \to F^*$ is an
isomorphism.

\begin{proof} (i) is clear since  $\Aut (F)$
acts on $F^{\times\flat}$.

We will prove (ii), (iii) simultaneously.

(a) follows from (2.1.1) since $\{ f,f
\}^{plain}=1$ for any $f$, and
(b) follows from Exercise (b) in 2.10(ii).

For the remainder of the proof, we
assume we are working in the
category of $\Q$-algebras. We can do
this since it suffices to check our statements
for the universal $R$; it has no torsion, so
we  loose no information tensoring  by
$\Q$.

Let us check that $\Z\times\G_m$ is orthogonal
to
$\W\times\hat{\W}$. Since we are in char 0 and
$\W\times\hat{\W}$ is connected we can
replace this group by its Lie algebra  which
is the subspace
$F^0 \subset F$ of formal power series $\Sigma
a_i t^i$ with
$a_0 =0$. The action of
homotheties $\subset \Aut (F)$ via
$t\mapsto at$ preserves
$F^0$ and has  trivial coinvariants there.
Since homotheties fix $\G_m \subset \Z\times
\G_m$ we see, by (i), that $\G_m$  is
orthogonal to
$F^0$. The action of homotheties on the
quotient $\Z =\Z\times\G_m /\G_m$ is trivial,
so the same argument shows that $\Z$ is also
orthogonal to
$F^0$, q.e.d.

It remains to compute the restriction of $\{\,
,\,\}^\flat$ to $\W\times\hat{\W}$. Again, it
suffices to check the formulas on the level of
Lie algebras.

The group pairing $f(t),
g(t^{-1})\mapsto r(f*g)$ from (c) looks as
follows. For $f(t)=\exp (
\alpha_1 t +\alpha_2 t^2 +..)$,
$g(t)=\exp(\beta_1 t +\beta_2 t^2 +..)$ one
has
$(f*g)(t)=\exp (-\sum i\alpha_i
\beta_i t^i )$, so $r(f*g)=\exp (\sum
i\alpha_i \beta_i )$. It is clear that the
corresponding infinitesimal pairing is
$a,b
\mapsto
\Res (bda)$.

  So it remains to
check (iii), i.e., to show that $[ t^m ,t^n
]^\flat =m\delta_{m,-n}$, which is an
easy  calculation (use presentation
2.10(iv) for
$F^\flat$).
\end{proof}

\medskip

{\it Remarks.} (i) Here is another
convenient formula.\footnote{We are grateful
to referee for pointing it out to us.} Let
$f\in R[[t]]$ be a formal power series with
invertible constant term, and $g\in
R[t,t^{-1}]$ a Laurent polynomial with one
coefficient invertible and the other ones
nilpotent. Then $\text{div} g$ is a Cartier
divisor supported at $t=0$ and, as
follows from Exercise (b) in 2.10(ii), one has
\eq{3.3.1}{\{ f ,g\}^\flat =
\text{N}_{\text{div} g}(f)^{-1}.}  

(ii) According to \cite{C}
the pairing $\{ \, ,\, \}^\flat$ is
non-degenerate, i.e., it identifies
the group ind-scheme $F^\times$ with its
Cartier dual. For an analytic version of this
statement see
\cite{D3}, sect.~4.

\subsection{The standard splittings and the
Weil reciprocity} Let
$R$,
$F_R$ be as in 3.1. Suppose we have
$R$-subalgebras
$O_R ,A_R
\subset F_R$ which are, respectively, c- and
  d-lattices in $F$
(see 2.13). By base change we can consider it
as an embedding $O, A\hra F$ of $\sO$-algebras
in $\sS_R$ (see 2.10). Then
$O^\times , A^\times
\subset F^\times$ and (2.13.3) gives
the standard splittings
\eq{3.4.1}{ s_O^c : O^\times \to
F^{\times\flat},\quad s_A^d :A^\times
\to F^{\times
\flat}.} In particular, the restriction of
$\{\, ,\, \}^\flat$ to $O^\times$ and
$A^\times$ vanishes.

\medskip

This situation occurs in geometry as follows.
   Let $X$ be a
smooth projective family of curves over $\Spec
R$, $D\subset X$ a  relative divisor such that
$U=X\smallsetminus D$ is affine. Let
$O_D$ be the algebra of functions on the formal
completion of $X$ at
$D$ and
$F_D$ its localization with respect to an
equation of $D$.\footnote{So $O_D :=
\leftlim \Gamma (X,\sO/\sO (-nD))$,
$F_D :=\rightlim\leftlim
  \Gamma (X, \sO (mD)/\sO
(-nD))$.} Then
$F_D$ is a Tate
$R$-module, and $O_D \subset F_D$ is a
c-lattice. We have the obvious embedding of
$R$-algebras $\sO (U) \hra F_D$ which
identifies $\sO (U)$ with a d-lattice in
$F_D$. So we have the standard
splittings
\eq{3.4.2}{s_{O_D}^c : O_D^\times \to
F_D^{\times\flat}, \quad s_{\sO (U)}^d :
\sO (U)^\times \to F_D^{\times\flat}.}

Let
$\{\, ,\,
\}^\flat_D$  be the symbol map  (3.1.2)  for
$F_D$. If $D$ is a disjoint union of $R$-points
$x_\alpha$ then $F_D =\Pi F_{x_\alpha}$ and
3.3\footnote{Together with Exercise (ii) in
3.1.} identifies
$\{ \, ,\, \}^\flat_D$ with the product of (the
inverse of) local tame symbols at $x_\alpha$.
If
$D$ is
\'etale over $\Spec R$ then one computes $\{
\, ,\, \}^\flat_D$ using Exercise (i) in 3.1.
So the vanishing of $\{ \, ,\, \}^\flat_D$ on
$\sO^\times (U) \subset F^\times_D$ is the
classical Weil reciprocity law. Passing to
Lie algebras  one gets, by
3.3(iii), the residue formula.

\subsection{ $\omega (F)^\times$ and the
extended Heisenberg }  Let $R$,
$F_R$ as in 3.1 and assume that locally in
flat topology of $\Spec R$ our $F_R$ is
isomorphic to a product of several copies of
$R((t))$. All our constructions  are local, so
  we tacitly assume that we work flat
locally on
$\Spec R$.

  So, localizing $R$ if needed, we get, by
3.3(iii),
  a canonical isomorphism
$\omega (F):= \omega_F \iso F^*$ of
free $F$-modules of rank 1. Denote by $\omega
(F)^\times$ the set of generators, i.e.,
invertible elements, of
$\omega (F)$.

The story of 2.16 in our
special setting looks as follows.

\medskip

  As in 2.14  set $W:= F\oplus F^* =F\oplus
\omega (F)$, so we have the orthogonal group
$\rO =\rO (W)$ and its super extension
$\rO^\flat $. One has an embedding
$F^\times \hra \rO$, $f\mapsto f^o :=$ the
diagonal matrix with entries $f, f^{-1}$.
By 2.15 the restriction of
$\rO^\flat$ to
$F^\times$ is the Heisenberg super extension
$F^\flat$.

Consider an embedding\footnote{We identify
$\nu\in
\omega (F)^{\times}$ with an  $F$-linear
isomorphism
$F\iso \omega (F)$,
$f\mapsto f\nu$.}
$\omega (F)^{\times}\hra
\rO (W)$ which identifies $\nu$ with the
anti-diagonal matrix
$\nu^o$ with entries $\nu,\nu^{-1}$. The
disjoint union of
$F^\times$ and $\omega (F)^{\times}$ is
a subgroup
$\tilde{F}^\times$ in
$\rO (W)$ which contains $F^{\times}$
as a normal subgroup of index 2; the
multiplication table is
\eq{3.5.1}{\nu^o {\nu'}^o =(\nu' /\nu )^o ,
\quad f^o
\nu^o =(f^{-1}\nu )^o =\nu^o f^o .}
Let $\tilde{F}^{\times \flat}$ be the
super extension of $\tilde{F}^\times$
induced  from $\rO^\flat (W)$; for $g\in
\tilde{F}^{\times}$ we denote the
corresponding super line by $\lambda_g$. In
particular, we have a super line
$\lambda_{\omega (F)^\times}$ on $\omega
(F)^\times$.

\medskip

Below we consider $\omega (F)^\times$ as an
$F^\times$-torsor with respect to the action
$f,\nu \mapsto f^{-1}\nu$. The multiplication
in $F^{\times\flat}$ lifts this action to an
action of $F^{\times\flat}$ on
$\lambda_{\omega (F)^\times}$.

Notice that over
$\omega (F)^{\times}$ the tensor square of
$\lambda$ is
  canonically trivialized: one
has $\lambda_{\nu^o}^{\otimes 2}=\sO_\nu$ since
$\nu^{o 2}  =1$. Therefore we have defined a
super
$\mu_2$-torsor
$\mu_{\omega (F)^{\times}}$ on
$\omega (F)^{\times}$.

The groupoid that
corresponds to
the $F^{\times\flat}$-action on
$\omega (F)^{\times}$ is a super extension
$\sG^{\flat}_{\omega (F)^{\times}}$ of
the simple transitive groupoid $\sG_{\omega
(F)^{\times}}$. Its action on $\lambda_{\omega
(F)^\times}$ identifies $\sG^{\flat}_{\omega
(F)^{\times}}$ with the
$\lambda_{\omega (F)^\times}$-split super
extension (see A2). Since
$\lambda_{\omega (F)^\times}$ comes from the
super $\mu_2$-torsor $\mu_{\omega
(F)^{\times}}$ we see that
$\sG^{\flat}_{\omega (F)^{\times}}$
comes from the
$\mu_{\omega (F)^{\times}}$-split super
$\mu_2$-extension $\sG^\mu_{\omega
(F)^{\times}}$.
  Explicitly, $\sG^\mu_{\omega
(F)^{\times}}\subset \sG^{\flat}_{\omega
(F)^{\times}}$ is the super
$\mu_2$-subextension whose action preserves
$\mu_{\omega
(F)^{\times}}\subset \lambda_{\omega
(F)^{\times}}$.

If
$1/2
\in R$ then $\mu_{\omega (F)^{\times}}$ is an
\'etale torsor and $\sG^\mu_{\omega
(F)^{\times}}$ is an
\'etale groupoid extension of
$\sG_{\omega (F)^{\times}}$. Therefore
$\sG^\mu_{\omega (F)^{\times}}$
splits canonically
  over the formal completion\footnote{An
$R'$-point of $\hat{\sG}_{\omega
(F)^{\times}}$ is a pair of $R'$-points of
$\omega (F)^{\times}$ that coincide modulo a
nilpotent ideal of $R'$.}
$\hat{\sG}_{\omega (F)^{\times}}$ of
$\sG_{\omega (F)^{\times}}$. The
corresponding splitting of $\sG_{\omega
(F)^{\times}}^{\flat}$ is called {\it the
canonical formal rigidification} of
$\sG_{\omega (F)^{\times}}^{\flat}$.

\medskip

{\it Remark.} The canonical formal
rigidification of $\sG_{\omega
(F)^{\times}}^{\flat}$ is completely
determined by the connection on
$\lambda_{\omega (F)^\times}$ defined by the
$\mu_2$-structure. The horizontal leaves of
this connection are just the orbits of the
{\it adjoint} action of $F^\times$ on
$\lambda_{\omega (F)^\times}\subset
\tilde{F}^{\times\flat}$.\footnote{Notice that
this action lifts the $F^\times$-action $f,\nu
\mapsto f^{-2}\nu$ on
$\omega (F)^\times$.}

\medskip

  The above  picture is compatible with
the standard splittings of the
Heisenberg group from 3.4. Namely, suppose we
have $O, A \subset F$ as in 3.4. Set
$O^\circ:= O^\bot, A^\circ := A^\bot $. So
$O^\circ \subset \omega (F)=F^*$ is an
$O$-submodule which is a c-lattice, $A^\circ
\subset \omega (F)$ is an
$A$-submodule which is a d-lattice.
Suppose that $O^\circ$ is generated, as an
$O$-module, by some element of $\omega
(F)^\times$,
and similarly,
$A^\circ$ is generated, as an $A$-module,
by some element of $\omega (F)^\times$.  Let
$O^{\circ\times},A^{\circ\times}\subset \omega
(F)^\times$ be the sets of such generators.
Then
$\tilde{O}^\times := O^\times \sqcup
O^{\circ\times}$, $\tilde{A}^\times :=
A^\times \sqcup A^{\circ\times}$ are subgroups
of
$\tilde{F}^\times$. Since
$\tilde{O}^\times
\subset \rO (W, O_W )$,
$\tilde{A}^\times \subset \rO (W, A_W )$
(see Remark (ii) in 2.15 and Remarks (i), (ii)
in 2.14 for notation) we have the standard
sections
\eq{3.5.2}{s_{O_W}^c : \tilde{O}^\times
\to
\tilde{F}^{\times\flat} , \quad
s_{A_W}^d :
\tilde{A}^\times
\to
\tilde{F}^{\times\flat} .} On
$O^\times ,A^\times \subset
F^\times_\omega $ these are splittings
$ s_{O}^c $, $ s_{A}^d$ from
(3.4.1) (see Remark (ii) in 2.15). The
restriction of (3.5.2) to
$O^{\circ
\times},A^{\circ\times}\subset \omega
(F)^\times$  trivialize the
restrictions to
$O^{\circ\times}$,
$A^{\circ\times}$ of the
$\mu_2$-torsor
$\mu_{\omega (F)^\times}$.
  The corresponding
splittings of the groupoid super extensions
$\sG^{\flat}_{O^{\circ\times}}$,
$\sG^{\flat}_{A^{\circ\times}}$ come
from sections $ s_{O}^c $, $ s_{A}^d$.

\medskip

{\it Remark.} With
$\omega (F)$ replaced by $F^*$ the above
constructions remain valid for every $F$ such
that $F^*$ is a free
$F$-module of rank 1. E.g.~one can take for
$F$ an algebra
$F_D$ from 3.4; the divisor $D$ need not be
\'etale over $\Spec R$.

\subsection{Remarks} We assume that $1/2 \in
R$.

(i)  Take any $\nu \in \omega (F)^\times$.
Then the canonical formal rigidification of
$\sG^{\flat}_{\omega (F)^\times}$
provides a section
$\nabla_\nu^\flat : F\to
F^{\flat}_\omega$. As folows from the first
Remark in 3.5 one has $\nabla^\flat_{f^{-2}\nu}
=\Ad_f \nabla^\flat_\nu$, i.e., by (3.1.6),
\eq{3.6.1}{\nabla_{f\nu}^\flat (a) =
\nabla_\nu^\flat (a)-
\frac{1}{2} \Res (ad\log f).}

Therefore the map $\nu \mapsto
\nabla^\flat_\nu$ identifies the
$\omega (F)=F^*$-torsor of sections  $F\to
F^\flat$ with the push-out of the
$F^\times$-torsor
$\omega (F)^\times$ by the homomorphism
$-\frac{1}{2}d\log : F^\times \to
\omega (F)$.

\medskip

(ii) According to
Remark (iii) in 2.16 the above
$\nabla^\flat_\nu$ can be also described as
follows. Consider the involution
$\fs_\nu :=$Ad$_{\nu^o}$ of $\rO (W)$ and
$\rO^\flat (W)$. It preserves
$F^\flat \subset {\frak o}^\flat (W)$
and acts on
$F \subset
  {\frak o} (W)$ as
multiplication by
$-1$. Our
$\nabla^\flat_\nu$ is an
$\fs_\nu$-invariant splitting
$ F
\to F^{\flat}$, and it is uniquely
characterized by this
property.\footnote{Indeed, if there were two
such, the difference would be a map $F\to\sO$
which is $\fs_\nu$-equivariant, where
$\fs_\nu$ acts by $-1$ on $F$ and by $+1$ on
$\sO$.}

\medskip

(iii) Set $O :=
R[[t]]\subset F=R((t))$. Consider the splitting
$s_O^c : O\to F^{\flat}$ from
(3.4.1). Then for every $a\in O$ and $f\in
F^\times$ one has
\eq{3.6.2}{
\nabla^\flat_{fdt}(a)=s_O^c (a)
-\frac{1}{2} \Res (ad\log f).} Indeed,
$\nabla^\flat_{dt} (a)=s^c_O
(a)$ by the discussion at the end of 3.5; the
general case follows from (3.6.1).

\medskip

(iv)
Let us describe the super $\mu_2$-torsor
$\mu_{\omega (F)^\times}$ for
$F=R((t))$ explicitly. We can assume that $R$
is reduced.  On  even connected components
(those of the forms
$t^{2n}dt$) our
$\mu$ is even, and it can be trivialized
choosing
$L=t^{-n}R[[t]]$ as in 2.17(i). On odd
components (those of the forms $t^{2n+1}dt$)
our
$\mu$ is odd. For a form $\nu = rt^{2n+1}dt
+..$ our $\mu_\nu$ is the
$\mu_2$-torsor of square roots of $r$ as
seen from (2.17.1) applied to
$L=t^{-n}R[[t]]$.

We see that the restriction of
$\sG^\mu_{\omega (F)^\times}$ to even
components
  is canonically
trivialized. Therefore the restriction of
$\sG^{\flat}_{\omega (F)^\times}$ to even
components is also canonically trivialized.
The restriction of
$\sG^\mu_{\omega (F)^\times}$ to odd
components has non-trivial monodromy.

\subsection{Some twists}
  Let $R$ be a commutative algebra,  $F_R$ a
topological $R$-algebra. Assume
that
$F_R$ is isomorphic to
$R((t))$, so $\omega (F)=F^*$. Let
$V_R$ be a finitely generated projective
$F_R$-module. Then $V_R$ is a Tate $R$-module
with respect to a natural topology whose base
is formed by sub
$R$-modules
$\Sigma I\cdot v_\alpha$ where $\{ v_\alpha
\}$ is a finite set of $F_R$-generators of
$V$, $I\subset F_R$ an open
$R$-submodule.\footnote{We already met this
topology in Example (ii) of 2.11.}

Consider the action of $F^\times$ on $V$ by
homotheties
$F^\times
\to \GL (V)$. Let $F^{\times\flat}_{(V)}$ be
the pull-back of the Tate super extension
$\GL (V)^\flat$. Its Lie algebra
$F_{(V)}^\flat$ is a central extension of $F$
by $\sO$.

{\it Remark.} Super
$\sO^\times$-extensions of
$F^\times$ form a Picard groupoid with respect
to Baer product (see A3). As follows from
2.10(iii),
$V\mapsto F_{(V)}^{\times \flat}$ is a
symmetric monoidal functor  (with respect to
direct sum of $V$'s).

\medskip

The following proposition-construction is
crucial for the definition of
$\varepsilon$-connection. The reader who is
willing to admit that $V$ is a free
$F$-module can use instead a  much shorter
equivalent construction from 3.9(i) below.

\medskip

Denote by $n$ the rank of $V_R$ (as of an
$F_R$-module). We can consider $n$ as a
locally constant function on $\Spec R$.
Let $F^{n\flat}$ be the Baer $n$-multiple of
the extension $F^\flat$.

\subsection{Proposition} Every
$R$-relative connection $\fd$ on the
$F_R$-line $\det_{F_R} V_R$ yields
a canonical isomorphism of
the Lie algebra extensions \eq{3.8.1}{\xi^\fd :
F^{n\flat}
\iso F^\flat_{(V)}}
such that for
$\eta \in \omega (F)$ and $a^\flat \in
F^{n\flat}$ lifting $a\in F$ one has
\eq{3.8.2}{\xi^{\fd+\eta} (a^\flat )-\xi^\fd
(a^\flat )=-\Res (a\eta )\in \sO \subset
F^\flat_{(V)}.}

\begin{proof} (a) It is convenient to use a
slightly different format. Namely, for every
$R$-relative connection $\nabla$ on
$V$ we will define a canonical isomorphism of
the Lie algebra extensions
\eq{3.8.3}{\xi^\nabla : F^{n\flat}
\iso F^\flat_{(V)}}
such that for any
$\rho \in \omega (F)\otimes_F \text{End}_F
(V)$  one has
\eq{3.8.4}{\xi^{\nabla +\rho} (a^\flat
)=\xi^\nabla (a^\flat )-\Res\,
\text{tr}(a\rho ).}

Such $\xi^\nabla$ amounts
to a datum of (3.8.1)  subject to (3.8.2).
Namely, $\xi^\fd$ that
corresponds to $\xi^\nabla$ is
determined by the condition
$\xi^{\text{tr}\nabla}=\xi^\nabla$
where tr$\nabla$ is the connection on
$\det_F V$ defined by
$\nabla$.\footnote{Notice that
$\nabla$ form an
$\omega (F)\otimes_F \text{End}_F V$-torsor
$\sC onn (V)$, and $\nabla
\mapsto
\text{tr}\nabla$ identifies the $\omega
(F)$-torsor of $\fd$'s with the $\omega
(F)$-torsor induced from $\sC onn (V)$ by the
map id$_{\omega (F)}\otimes$tr$:
\omega (F)\otimes_F
\text{End}_F V \to \omega (F)$.}

\medskip

(b) To construct $\xi^\nabla$ we use a
``geometric"
  interpretation of $F^\flat_{(V)}
$ from
\cite{BS} which we briefly recall now. The
picture below is compatible with base change,
so we skip $R$ from the notation (replacing it
by $\sO$ when necessary, see 2.10).

Set $V' := \Hom_F (V ,F)$,
$V^\circ :=
\Hom_F (V,
\omega (F))=\omega (F)\otimes_F V'$. These are
finitely generated projective
$F$-modules. As a Tate $\sO$-module $V^\circ$
is canonically isomorphic to the dual Tate
$\sO$-module to $V$: the pairing $V^\circ
\times V \to \sO$ is $v^\circ ,v\mapsto  \Res
(v^\circ (v))$. One has a usual identification
$V\otimes_F V' =\text{End}_F V$.

  The diagonal embedding $\Delta :\Spec
F\hra
\Spec (F\hat{\otimes}_\sO F)$ is a Cartier
divisor.\footnote{For
$F_R \iso R((t))$ one has $F\hat{\otimes}_R F
\iso R[[t_1 ,t_2 ]][t_1^{-1}, t_2^{-1}]$.}
Consider an exact sequence of
$F\hat{\otimes}_\sO F$-modules  \eq{3.8.5}{0 \to
V\hat{\otimes}_\sO V^\circ \to
V\hat{\otimes}_\sO V^\circ (\Delta )\to
\text{End}_F V \to 0} where the right arrow is
the residue around the diagonal along the
second variable. Let \eq{3.8.6}{ 0\to
\sO \to F^{\flat '}_{(V)}\to F \to 0} be its
pull-back by $F\hra \text{End}_F V$, $f\mapsto
f\cdot id_V$, and push-forward by
$V\hat{\otimes}_\sO V^\circ \to \sO$,
$v\otimes v^\circ
\mapsto -\Res (v^\circ (v))$. Notice that
$F^{\flat '}_{(V)}$ inherits from $V\hat{\otimes}_R
V^\circ (\Delta )$ an
  $F$-bimodule structure. It defines a Lie
bracket on
  $F^{\flat '}_{(V)}$ with the adjoint action of
$f^{\flat '}$ coming from the commutator with
$f$ with respect to the $F$-bimodule
structure. So
  $F^{\flat '}_{(V)}$ is a central Lie algebra
extension of the commutative Lie algebra $F$.

Now there is a canonical isomorphisms of
central $\sO$-extensions \eq{3.8.7}{  F^{\flat
'}_{(V)}\iso F^{\flat}_{(V)}.} Namely, for
$f^{\flat '}\in F^{\flat '}_{(V)}$ the
corresponding $f^\flat
\in F^\flat_{(V)}$ is defined as follows.
Choose $k=k(t_1 ,t_2 )dt_2 \in V\boxtimes
V^\circ (\Delta )$ that lifts $f^{\flat '}$. It
defines an ``integral oprator"\footnote{Here $k(t_1 ,t_2)
v(t_2 )dt_2 \in V\boxtimes\omega (\Delta
)$.}  \eq{3.8.8}{A_k : V \to
V,
\quad A_k (v) (t_1):=  -\Res_{t_2 =0} k(t_1
,t_2) v(t_2 )dt_2 .} It is easy to see that
$A_k \in\fgl_d (V)$ and $A_k^\infty =
f^\infty$. As was explained in the end of
2.13, such $A_k$ defines a lifting of $f$ to
$F^{\flat}_{(V)}$; this is our $f^\flat$. The
independence of the auxiliary choices is
immediate.

\medskip

(c) Now we are ready to define (3.8.3) . Let
$F^{\flat '}$ be the central extension (3.8.6)
for
$V=F$. Its
Baer $n$-multiple $F^{n\flat '}$ is the
push-forward of the exact sequence $0\to
F\hat{\otimes}_\sO
\omega (F) \to F\hat{\otimes}_\sO \omega
(F)(\Delta )\to F \to 0$ by the map
$F\hat{\otimes}_\sO \omega
(F)\to \sO$, $f\otimes \nu \mapsto -n\Res (f\nu
)$. According to (3.8.7) we can rewrite
  (3.8.3)  as an
isomorphism $\xi^{\nabla} : F^{n\flat '}\iso
F^{\flat '}_{(V)}$.

Our $\xi^\nabla$ will come from certain
  morphism of
$F\hat{\otimes}_\sO F$-modules
$\Xi^\nabla :F\hat{\otimes}_\sO \omega
(F)(\Delta )\to V\hat{\otimes}_\sO V^\circ
(\Delta )$. Since
$F\hat{\otimes}_\sO
\omega (F)(\Delta )$ is a free
$F\hat{\otimes}_\sO F$-module of rank 1, such
$\Xi^\nabla$ is multiplication by a section
$\gamma^\nabla \in V\hat{\otimes}_\sO V'$. This
$\gamma^\nabla$ must satisfy the condition
$\gamma^\nabla|_\Delta =id_V \in
V\otimes_F V'=\text{End}_F V$ in order to
assure that
$\xi^\nabla$ is well-defined. Notice that
$\xi^\nabla$ depends only on the restriction
of $\gamma^\nabla$ to the first infinitesimal
neighbourhood $\Delta^{(1)}$ of $\Delta$.

Let $\nabla'$ be the connection on $V'$ dual to
$\nabla$. Denote by $\tilde{\nabla}$ the
$\sO$-relative connection on the
$F\hat{\otimes}_\sO F$-module
$V\hat{\otimes}_\sO V'$ defined by $\nabla$ and
$\nabla'$. The restriction of the relative
K\"ahler differentials of $F\hat{\otimes}_\sO
F$ to $\Delta$ equals $\omega (F) \oplus
\omega (F)$, so the composition
of $\tilde{\nabla}$ with the restriction to
$\Delta$ is a morphism
$\bar{\nabla }:V\hat{\otimes}_\sO V' \to
(\omega (F)\oplus
\omega (F))\otimes_F
V\otimes_F V'$.

Now our $\gamma^\nabla \in
V\hat{\otimes}_\sO V'$ is any section
killed by $\bar{\nabla}$ whose restriction to
$\Delta$ equals $id_V$. Such $\gamma^\nabla$
exists since $id_V$ is a horizontal section of
End$_F V=V\otimes_F V'$, and its restriction
to $\Delta^{(1)}$ is uniquely defined. Thus we
defined  $\xi^\nabla$. Formula (3.8.4)  holds
since for $\rho =\rho (t)dt$ one
has $\gamma^{\nabla +
\rho}=\gamma^\nabla +\rho (t_1 )(t_2 -t_1 )
$ on $\Delta^{(1)}$.
\end{proof}

\subsection{Remarks} (i) Isomorphisms (3.8.1)
and (3.8.3) are compatible with direct sums of
$V$'s. Precisely, both $F^{n\flat}$ and
$F^\flat_{(V)}$ transform direct sum of
$V$'s into the Baer sum of extensions (this is
true on the group level; see 2.10(iii) for
$F^{\times\flat}_{(V)}$). With respect to
these identifications $\xi^{\oplus
\nabla_\alpha}$ equals the Baer sum of
$\xi^{\nabla_\alpha}$.

Therefore if $V$ is a {\it free}
$F$-module then isomorphisms
(3.8.1)  or (3.8.3)  can be described directly as
follows. As we mentioned above, 2.10(iii)
yields a canonical identification of super
extensions
$F^{\times\flat}_{(F^n
)}=F^{\times n\flat}$.\footnote{The latter
extension is the Baer $n$-multiple of
$F^{\times \flat}$.} So every
$F$-basis
$\phi : F^n
\iso V$ yields  an isomorphism  \eq{3.9.1}{
\xi^\phi :
  F^{\times n\flat} \iso
F_{(V)}^{\times\flat}.} The
corresponding isomorphism of Lie algebras is
$\xi^{\nabla_\phi}$ of (3.8.3)  where
$\nabla_\phi$ is the connection that
corresponds to the trivialization $\phi$.
You define then $\xi^\nabla$ for arbitrary
$\nabla$  using
$\xi^{\nabla_\phi}$ and (3.8.4).

\medskip

(ii)  Isomorphisms $\xi^\phi$
of (3.9.1) satisfy the property
\eq{3.9.2}{\xi^{\phi g} (f^\flat )=\xi^\phi
(f^\flat )\{ \det g ,f
\}^\flat . } Here $f\in F^\times$, $f^\flat
\in F^{\times n\flat}$ is a lifting of $f$,
$g\in \GL (n,F)$ is an invertible matrix,
$\det g\in F^\times$ its determinant, and
$\{\, ,\,
\}^\flat$ is (3.1.2) . Indeed, we can rewrite
our statement as the equality $\{ g,f_V
\}^\flat =\{
\det g ,f
\}^\flat$ where the l.h.s.~is the Tate
commutator pairing in $\GL (n,F) \subset
\GL_R (F^n )$ (and $f$ there means $f\cdot
id_{F^n}$) . This property is clear for the
diagonal and unipotent matrices
$g$ from 2.10(iii),(iv). This proves
our statement in case when $R$ is a local
Artinian ring (then $F=R((t))$ is also local
Artinian, so every $g$ is a product of
unipotent and diagonalizable matrices). The
general case reduces to that one as follows.

Notice that the pairing $\GL_n (F)\times
F^\times \to R^\times$, $g,f \mapsto \{ g,f
\}^\flat$, is continuous. Precisely, for every
$g,f$ there exists $N\gg 0$ such that for
every $\beta \in t^N \text{Mat}_n (R[[t]])$,
$\alpha\in t^N R[[t]]$ one has $g+\beta \in
\GL_n (F)$, $f+\alpha \in F^\times$, and $\{
g,f\}^\flat =\{ g+\beta , f+\alpha \}^\flat$.
Same is true for $\{ \det g, f\}^\flat$. So we
can assume that all but finitely many
coefficients of $f$ and entries of $g$ are
non-zero. Replacing $R$ by its subring
generated by these coefficients and the
inverses of the top non-nilpotent coefficients
of $\det g$ and $f$ we can assume that $R$ is
finitely generated over $\Z$, hence $R$ is
Noetherian. We want to prove that two elements
of $R^\times$ are equal. It suffices to check
this on every infinitesimal neighbourhood of
every point of $\Spec R$, and we are done.

\medskip

(iii) The involution $\fs_\nu$ from 3.6(ii) can
be easily described in terms of the
identification
$F^\flat =F^{\flat '}$ (case $V=F$ of (3.8.7)).
Namely, for $\nu =\nu (t)dt$ $\fs_\nu$
comes from the involution of
$F\hat{\otimes}_R \omega (F)(\Delta )$ which
sends
$f(t_1 ,t_2 )\nu (t_2 )dt_2 $ to $f(t_2 ,t_1
)\nu (t_2 )dt_2$. Therefore the
$\fs_\nu$-fixed section
$\nabla_\nu^\flat : F \to F^{\flat '}$ from
3.6(i) is
  \eq{3.9.3}{f(t) \mapsto \frac{
f(t_1)dt_2}{t_2 -t_1} +\frac{1}{2}(f'(t_1
)+f(t_1)\frac{\nu'(t_1 )}{\nu (t_1 )})dt_2 . }

\subsection{The canonical flat connections}
Let us rephrase the material of
3.5 and 3.8 in a format to be used in sect.~4.

Our input is a triple $(R,F ,V)$ where $R$ is
a $\Q$-algebra, $F=F_R$ a topological
$R$-algebra, $V=V(F)$ a finitely
generated projective
$F$-module. We assume that
locally in the flat topology of $\Spec R$ our
$F_R$ becomes a finite direct product
$\Pi F_i$ where every $F_i$ is isomorphic to
$R((t))$. All our constructions  are
$\Spec R$-local, so we can assume that the
above decomposition holds on $\Spec R$ itself.
Then $V=\Pi V_i$ where $V_i$ is a finitely
generated projective $F_i$-module.

By 3.3(iii) one
has a canonical identification $F^* =\omega
(F)=\Pi \omega (F_i )$. Recall that
$F^\times =\Pi F_i^\times$ is a formally smooth
group ind-scheme. We have an
$F^\times$-torsor $F^{*\times}=\omega
(F)^\times =\Pi \omega (F_i )^\times$.

As in 3.7 we notice that $V$ is a Tate
$R$-module in a natural way, so the
$F^\times$-action on $V$ by homotheties
$F^\times \hra \GL (V)$
defines\footnote{Pulling back  the Tate
extension $\GL^\flat (V)$.} the super extension
$F^{\times\flat}_{(V)}$ of
$F^\times$.  The embeddings $F_i^\times \hra
F^\times$ lift canonically to mutually
commuting embeddings
$F^{\times\flat}_{i (V_i )}\hra
F^{\times\flat}_{(V)}$  (see 2.10(iii)). For
$f\in F^\times$ we denote its super line in
$F^{\times\flat}_{(V)}$ by $\lambda_f^{(V)}$,
so  $\lambda_{(f_i )}^{(V)}=\otimes
\lambda_{f_i}^{(V_i )}$.

Let $\sE$ be any
super line bundle
$\sE$ on $\omega
(F)^\times$ equipped with an
action of $F^{\times\flat}_{(V)}$
which lifts the action $f,\nu
\mapsto f^{-1}\nu$ of $F^\times$ on
$\omega (F)^\times$.\footnote{By definition,
such
$\sE$ is a rule that assigns:
(i) to every point $\nu \in
\omega (F)^{\times}(R')$ a  super
$R'$-line $\sE_\nu$;
(ii) to every $f\in F^\times
(R)$ an isomorphism
$\lambda_{f}^{(V)}\cdot
\sE_\nu \iso
\sE_{f^{-1}\nu}$; (iii)
to every morphism
$r : R'
\to R''$ an isomorphism $\sE_{r\nu} =\sE_\nu
\otimes_{R'}R''$.
This datum should satisfy the obvious
compatibilities.} We are not interested in its
nature at the moment.
  Notice that an
$F^{\times \flat}_{(V)}$-action on $\sE$
amounts to a collection of mutually commuting
actions of
$F_{i (V_i )}^{\times \flat}$  lifting
the
$F_i^\times$-action along $\omega (F_i
)^\times$.

\medskip

The key structure on $\sE$ that arises
automatically is:

{\it  a rule that assigns to every
connection $\fd$ (relative to $R$) on
$\det V$  a flat connection $\nabla^\fd$ on
$\sE$.}

Let us define $\nabla^\fd$  assuming that
$F$ is isomorphic to
$R((t))$. In general case this will
  define $\nabla^\fd$ separately in every
$\omega (F_i )^\times$-direction. Since $\omega
(F)^\times =\Pi \omega (F_i )^\times$ this
determines our connection; its flatness
is immediate from the construction.

Denote by $\hat{F}^\times$ the formal
completion of $F^\times$. Let
$\hat{F}^{\times\flat}$,
$\hat{F}^{\times \flat}_{(V)}$ be the
restrictions  of the super
extensions
$F^{\times\flat}$, $F^{\times
\flat}_{(V)}$ to
$\hat{F}^\times$; these are central
$\sO^\times$-extensions of $\hat{F}^\times$.
Since we are in characteristic 0, the
isomorphism of Lie algebras (3.8.1)  amounts to
an isomorphism $\xi^\fd :
\hat{F}^{\times n\flat}\iso
\hat{F}^{\times \flat}_{(V)}$. The
$F^{\times \flat}_{(V)}$-action on $\sE$
yields, via $\xi^\fd$, an action of
$\hat{F}^{\times n\flat}$ on $\sE$ which
is the same as an action of the corresponding
groupoid
$\hat{\sG}_{\omega (F)^\times}^{n\flat}$
(see 3.5). Now the canonical formal
rigidification
$
\hat{\sG}_{\omega
(F)^\times} \to\hat{\sG}_{\omega
(F)^\times}^{\flat}$ from 3.5 provides an
$\hat{\sG}_{\omega (F)^\times}$-action on
$\sE$. According to Grothendieck \cite{Gr},
such action amounts to a flat connection. This
is our
$\nabla^\fd$.

The construction of $\nabla^\fd$ is compatible
with base change.

\medskip

{\it Remarks.} (i) The connection $\nabla^\fd$
acts on $\sE$ by the
structure $F^{\times\flat}_{(V)}$-action via a
family of splittings $\nabla_{\nu}^\fd : F \to
F^\flat_{(V)}$, $\nu \in \omega (F)^\times$,
where
\eq{3.10.1}{ \nabla_\nu^\fd := \xi^\fd
\nabla_\nu^\flat .} Here $\nabla_\nu^\flat
:F\to F^\flat \to F^{n\flat}$ was defined in
3.6(i).

(ii) For $\chi \in \omega (F)$ and $a\in F$
one has (see (3.8.2)) \eq{3.10.2}{\nabla^{\fd
+\chi }_\nu (a)=\nabla^\fd_\nu (a)-\Res (a\chi
).}

(iii) Suppose that $F=R((t))$ and $V = V_R
((t))$ where $V_R$ is a projective $R$-module
of  rank $n$. Set
$O:=R[[t]]\subset F$, $V(O ) :=
V_R [[t]]\subset V$. According to (2.13.3)
we have a section $s_{V(O)}^c : O^\times
\to F^{\times\flat}_{(V)}$, so
$\sE$ is an
$O^\times$-equivariant bundle.  Set
$\chi^\fd :=
\fd (\gamma )/\gamma
\in
\omega (F)$ where $\gamma$ is any
trivialization\footnote{$\gamma$ exists locally
on
$\Spec R$.} of
$\det V_R$. Then for
$a\in O$, $f\in F^\times$ one has
\eq{3.10.3}{\nabla^\fd_{fdt} (a)
=s_{V(O)}^c (a) -\Res\, (a\chi^\fd
+\frac{n}{2}a\, d\log f ).} Indeed, if $\fd$
came from a trivialization of $V_R$ then it is
formula (3.6.2); the general case follows from
(3.10.2).

In particular, if $\fd$ is non-singular
with respect to
$\det V(O)\subset \det V$
(i.e., $\chi^\fd \in \omega (O)
$) then the
restriction of $\nabla^\fd$ to the
$O^\times$-torsor $\omega (O)^\times
\subset \omega (F)^\times$ coincides with the
``constant" connection defined by the
$O^\times$-action on $\sE$.

(iv) Of course, to define $\nabla^\fd$ we need
only the action of
$\hat{F}^{\times\flat}$ on
$\sE$. The whole $F^{\times\flat}$-action
will be used in the next subsection.

\subsection{The global picture}
Now we have $X, D, U$  as in 3.4. We assume
that $D$ is
\'etale over
$\Spec R$.   Let
$V$ be a vector bundle on $U$, i.e., we have a
finitely generated projective $\sO (U)$-module
$V(U)$. Then $F:=F_D$ and $V(F)
:=F\otimes_{\sO (U)}V(U)$ satisfy the
assumptions of 3.10.

Set $F^{\times\flat}_{(V)}:=
F^{\times\flat}_{(V(F) )}$. Since
$V(U)\subset V(F)$ is a d-lattice, (2.13.3)
lifts the embedding
$\sO (U)^\times
\hra F^\times$ to a splitting\footnote{For
$V=\sO_U$ this is (3.4.2).}
\eq{3.11.1}{s_{V(U)}^d : \sO (U)^\times \to
F^{\times\flat}_{(V)}.}

Consider any $\sE$  as in 3.10. Then $\sO
(U)^\times$ acts on $\sE$ via the splitting
$\bar{s}_{V(U)}^d : \sO (U)^\times \to
F^{\times\flat}_{(V)}$ defined by
(3.11.1).

Let $\omega =\omega_{U/R}$ be the
canonical bundle. {\it We assume that
$\omega (U)$ is a free $\sO (U)$-module} (of
rank 1), i.e., the set of generators $\omega
(U)^\times \subset \omega (U)$ is non-empty.
Then $\omega (U)^\times \subset \omega
(F)^\times$ is an $\sO (U)^\times$-torsor.

Consider the restriction $\sE_{\omega
(U)^\times}$ of
$\sE$ to $\omega (U)^\times$. The  $\sO
(U)^\times$-action yields a ``constant" flat
connection on $\sE_{\omega (U)^\times}$ which
we denote by $\nabla^0$.

Suppose we have a connection $\fd$ on the
line bundle $\det V$ on $U$. Let
$\nabla^{\fd_F}$ be the connection on $\sE$
that corresponds, by the above, to the
  restriction $\fd_F$ of $\fd$ to $F$.

\subsection{Lemma}  The restriction of
$\nabla^{\fd_F}$ to $\omega (U)^\times$
is equal to $\nabla^0$.

\begin{proof} (a) By (3.10.2) and the residue
formula the restriction of
$\nabla^{\fd_F}$ to $\omega (U)^\times$ does
not depend on $\fd$.

(b) Assume that $ V$ is a trivial vector
bundle on $U$. Fix an isomorphism $\sO_U^n
\iso V$ and let $\fd$ be the connection that
corresponds to the trivialization $\det\phi$
of $\det V$. Then $\phi$ yields an
identification $\xi^\phi : F^{\times n\flat}
\iso F_{(V)}^{\times\flat}$ of (3.9.1), hence
we get an
$F^{\times n\flat}$-action on $\sE$. By
3.9(i) our $\nabla^{\fd_F}$ comes from this
action via the canonical formal rigidification
(see 3.5). We are done by the discussion at
the end of 3.5, since
$\xi^\phi$  identifies sections
$s_{\sO (U)^\circ}^d$ and $s_{V(U)^\circ}^d$.

(c) Let us reduce
the general situation to the case of trivial
$V$. Our
$X,D,V$, some
$\nu\in\omega (U)^\times$, and  $\sE_\nu$
are defined over a finitely generated
subring of $R$. Since  $\sE$ is an equivariant
bundle over a torsor, we see that the whole our
datum is defined over this subring. Therefore
we can assume that
$R$ is a finitely generated
$\Q$-algebra. Then in order to check that two
connections are the same it suffices to do
this on ifinitesimal neighbourhoods of points
in
$\Spec R$. So, by base change, we can assume
that $R$ is a local Artinian $\Q$-algebra.
Then one can choose
$D' = D\sqcup D''$ such that $V$ is trivial on
$U':=X\smallsetminus D' \subset U$. It remains
to show that our statement for $(U,V_U )$
follows from the one for $(U', V_{U'})$.

Set $F=F_D$, $F'=F_{D'}$, $F'' =F_{D''}$;
let $O'' \subset F''$ be the
product of local rings,
so we have
$V(O'')\subset V(F'')$. The projection
$V(F') \twoheadrightarrow V(F )\times
V(F'')/V(O'')$ is preserved by
$F^\times \times O^{''\times}$-action, so the
restriction of $F^{''\times\flat}_{(V)}$ to
$F^\times \times O^{''\times}$
  equals the pull-back of
$F_{(V)}^{\times\flat}$.\footnote{i.e., we
trivialize $F^{''\times\flat}_{(V)}$ over
$O^{''\times}$ by means of $s_{V(O)}^c$, see
Remark (iii) in 3.10.} The exact sequence
$0\to V(U)\to V(U') \to V(F'')/V(O'') \to 0$
shows that this isomorphism identifies the
restriction of $s_{V(U')}^d$ to $\sO (U)^\times
\subset \sO (U')^\times$ with $s_{\sO (U)}^d$.
Let $\sE'$ be an
$F^{'\times \flat}_{(V)}$-equivariant
super line on $\omega (F')^\times =\omega
(F)^\times
\times \omega (F'')^\times$ whose restriction
to
$\omega (F)^\times
\times
\omega (O'' )^\times$  coincides with the
pull-back of $\sE$. Now on $\omega (F)^\times
\times
\omega (O'' )^\times$ the $\sO
(U)^\times $-action
on $\sE'$ (which comes from the $\sO
(U')^\times$-action) coincides with the
$\sO (U)^\times$-action on $\sE$, and for
every connection $\fd$ on $\det V$ the
connection
$\nabla^{\fd_{F''}}$ on $\sE'$ coincides with
  the connection $\nabla^{\fd_F}$ on $\sE$. We
are done.
\end{proof}

\section{The $\varepsilon$-factors}

 From now on we work over a base field $k$ of
characteristic 0. ``Absolute connection"
means a connection relative to $k$.

\subsection{The $\varepsilon$-lines} Let $K$ be
a field,
$F:=K((t))$, $O:= K[[t]]\subset F$; we consider
$F$ as a topological
$K$-algebra. Let $\omega =\omega (F ):=
\omega (F/K)$ be the 1-dimensional $F$-vector
space of K\"ahler differentials, $\omega
=Fdt$,  Der$(F/K)=\Hom_F (\omega  ,F)=
F\partial_t$ the one of vector
fields.

Suppose $V$ is an $F$-vector space of dimension
$n$ equipped with a connection $\nabla : V\to
\omega \otimes_F V$. For a non-zero $\nu \in
\omega$ let $\tau_\nu := \nu^{-1}$ be
the corresponding vector field. Consider
$V$ as a Tate $K$-vector space (see 2.11), so
we have a continuous operator
$\tau_\nu  = \nabla (\tau_\nu ) : V\to V$. It
is well-known\footnote{and follows from
the fact that
$V$ is generated by a single vector, see
below.} that the corresponding asymptotic
operator
$\tau^\infty_\nu$ acting on $V^\infty$ (see
2.13) is invertible. So we have the
$\Z$-graded super
$K$-line
\eq{4.1.1}{\sE (V,\nabla )_\nu =\sE (V)_\nu :=
\det (V^\infty ,\tau_\nu^\infty )} called
{\it the $\varepsilon$-line of $(V,\nabla )$
at
$\nu$}.

{\it Remark.} The
determinant super $K$-line of the de Rham
cohomology $\det H^\cdot_{dR}(V)$ is {\it
canonically trivialized} (see 5.9(a)(iv)).

\medskip

Sometimes it is convenient to compute $\sE
(V)_\nu$ as follows. Set $\tau =\tau_\nu$.
Denote by
$\sD_F$ the associative
$K$-algebra of differential operators; let
$\sD_F^{\le i}$ be its standard filtration
by degree. So $F$ and $V$ are left
$\sD_F$-modules. It is well-known (see
e.g.~\cite{D1} p.42 or \cite{M} III 1.1) that
$V$ is generated by a single vector $e$. Then
$e,
\tau e,..,\tau^{n-1}e$ is an $F$-base in
$V$, so   $\tau^n e=
a_{n-1}  \tau^{n-1} e+..+ a_0 e$ for certain
$a_i = a_i (t)\in F$. Set $A_\tau
:= (-1)^{n-1}\tau^n +
\mathop\sum\limits_{i\le n-1}
(-1)^{i}\tau^i a_i
\in \sD_F$. This is a non-zero differential
operator, so the corresponding asymptotic
operator\footnote{We consider $F$ as a
Tate $K$-vector space.}
  is invertible, and we have the
$\Z$-graded super $K$-line $\det
(F^\infty ,A^\infty_\tau )$.

\subsection{Lemma.} There is a natural
isomorphism of $\Z$-graded super lines
\eq{4.2.1}{\sE (V,\nabla )_{\nu}=\det
(F^\infty ,A^\infty_\tau ).}

\begin{proof} By definition, $\sE (V)_{\nu}=
\det (V/L_V ,V/L_V , \tau^\infty )$ where
$L_V \subset V$ is any c-lattice. Take for
$L_V$ the sum
$ L_F e+\tau (L_F
e)+..+\tau^{n-1}(L_F e)$ where
$L_F \subset F$ is any c-lattice in $F$. To
compute our determinant line consider the Tate
submodule $V' := Fe +\tau
(Fe)+..+\tau^{n-2}(Fe)$.
Since
$\tau$ induces  isomorphisms $V' \iso
\tau (V' )$ and $L_{V}\cap V' \iso
L_V \cap \tau (V')$ we see
that\footnote{Use (2.9.1) for
$M=N=V/L_V$, $f^\infty =\tau^\infty$, and
1-term filtrations
$ N_1 :=V'/ L_{V}\cap V' \subset N$,
$M_1 :=\partial_t (V')/L_V \cap \tau
(V')\subset M$; then apply (2.6.2) to
trivialize the line
$\det (M_1 ,N_1 ,
\tau^\infty|_{N_1 })$ by $\det (\tau|_{N_1})$.}
$\sE (V)_{\nu}=\det (V/(L_V \oplus \tau (V')) ,
V/(L_V +V'),\tau^\infty )$. The identifications
$F\iso V/V'$, $f\mapsto \tau^{n-1}(fe)$,
and $F\iso V/\tau (V')$, $f\mapsto fe$,
are compatible with c-lattices, and
$\tau (\tau^{n-1}(fe))$ is equal to
$A_\tau (f)e$ modulo
$\tau (V')$. Therefore $\sE
(V)_{\nu}=\det (F/L_F ,F/L_F, A_\tau^\infty )$,
which yields (4.2.1). We leave it to the
reader to check that our identification does
not depend on the choice of $L_F$.
\end{proof}

\subsection{Corollary} The degree of the
$\Z$-graded super line $\sE (V,\nabla )_\nu$
is equal to $i(\nabla ) + (v(\nu )+1) n$
where $i (\nabla )$ is the irregularity of
$\nabla$,\footnote{See e.g.~\cite{D1}
  or \cite{M} IV, sect.~4.} and
$v(\nu )$ the valuation of
$\nu$.\footnote{i.e., $\nu =ct^{v(\nu)}dt
+$ higher order terms,
$c\neq 0$.}

\begin{proof}
Write $\nu = fdt/t$, so
$\tau_\nu = f^{-1}t\partial_t$. The degree of
$\sE (V)_\nu$ is the index of the operator
$\tau^\infty$ acting on $V^\infty$. The
index of $(f^{-1})^{\infty}$ acting on
$V^\infty$ is equal to $v(f)n$. Thus it
suffices to check 4.3 for $\nu =dt/t$.

By 4.2 the degree
of $\sE (V)_{dt/t}$ is equal to the index of
the operator $A_{t\partial_t}^\infty$. On the
other hand, acording to \cite{D1} p.110
or \cite{M} IV(4.6), one has
$i(\nabla )= \sup (0,-v(a_i ))$ where $a_i \in
F$ are coefficients of $A_{t\partial_t}$.
Notice that
$A_{t\partial_t}$ shifts the filtration $t^j
K[[t]]$ on $F$ by $i(\nabla )$, and the
corresponding operator between the associated
graded quotients $t^j K \to t^{j-i(\nabla )}K$
is non-zero for almost all
$j$'s. So the index of
$A^\infty_{t\partial_t}$ equals $i(\nabla )$,
q.e.d.
\end{proof}

\subsection{Families of $\varepsilon$-lines}
Let us consider the picture of 4.1 depending
on parameters. So we have $(R,F,V,\nabla )$
where $(R,F,V)$ are as in 3.10, $\nabla$ is an
$R$-relative connection on $V$.

For $\nu \in \omega (F)^\times$ (it exists
flat locally on $\Spec R$) let
$\tau_\nu :=
\nu^{-1}\in \text{Der} (F/R)$ be the
corresponding vector field, so we have
$\nabla (\tau_\nu ): V\to V$.  We say that
our $\Spec R$-family is {\it
$\varepsilon$-nice} if the asymptotic operator
$\nabla (\tau_\nu )^\infty : V^\infty \to
V^\infty$ is invertible (here we consider $V$
as a Tate
$R$-module, see 2.11, 2.13). This
property does not depend on the choice of
$\nu$ and it is local with respect to flat
topology of $\Spec R$.

In $\varepsilon$-nice
situation for $\nu\in \omega (F)^\times$  we
have
$\Z$-graded super
$R$-line  \eq{4.4.1}{\sE (V,\nabla
)_\nu =\sE (V)_\nu :=
\det (V^\infty ,\nabla (\tau_\nu )^\infty )}
called {\it the
$\varepsilon$-line}. The base change of an
$\varepsilon$-nice family is
$\varepsilon$-nice, and
  $\varepsilon$-lines are
compatible with the base change.

{\it Remark.} The above construction is
compatible with disjoint union of families of
formal discs. Namely, assume that we have a
finite collection of families $(F_{\alpha
R},V,\nabla_\alpha )$. Set
$F_R :=\Pi F_{\alpha R}$,
$V:= \Pi V_\alpha$, $\nabla :=\Pi
\nabla_\alpha$. Then
$\omega (F)^\times  =\Pi
\omega (F_\alpha )^\times$,
$(F,V,\nabla )$ is
$\varepsilon$-nice if and only if such is
every $(F_\alpha ,V_\alpha ,\nabla_\alpha )$,
and, by (2.6.3), there is a canonical
isomorphism \eq{4.4.2}{\sE (V,\nabla
)=\boxtimes \sE (V_\alpha ,\nabla_\alpha ).}

\medskip

How  can one determine if a family is
$\varepsilon$-nice? We know (see 4.1) that
this is always true if $R$ is
a field, hence if $R$ is an Artinian algebra
(see Remark (b) in 2.5). Thus every family is
$\varepsilon$-nice at the generic point.

According to  2.5(ii), a $\Spec R$-family
is
$\varepsilon$-nice if and only if for a
c-lattice\footnote{We consider $V$ as a Tate
$R$-module.}
$L\subset V$ the quotient $R$-module $V/(L+
\nabla (\tau_\nu )(V)) $ is finitely generated,
and the irregularity function $i: \Spec
R\to\Z$, $x
\mapsto i (\nabla_x )$, is locally constant
(see 4.3, 2.4). The latter
condition is superfluous if $R$ is Noetherian
(see Remark (c) in 2.5). The absence of jumps
of irregularity alone does not imply
the family
is $\varepsilon$-nice.\footnote{E.g.,
consider
$R=\Spec k[x]$, $F=R((t))$, $V=F$, $\nabla
(\partial_t )=\partial_t + x/t$.}

\medskip

{\it Remarks.} (i) If $R$ is a local ring then
we are not aware of any example of a family
with constant irregularity which is not
$\varepsilon$-nice.

(ii) If $\nabla$ comes from an absolute
connection on $V$ then it looks probable that
the absence of jumps of irregularity implies
that the family is $\varepsilon$-nice.

\subsection{The $\varepsilon$-connection}
Assume that we are in situation of 4.4 and our
family is
$\varepsilon$-nice.
The compatibility with base change shows that
$\sE (V)_\nu$ form a $\Z$-graded super line
bundle over the $F^\times$-torsor
$\omega (F)^\times$. Let us show that $\sE
(V)$ is canonically rigidified with respect to
infinitesimal variations of $\nu$, i.e., $\sE
(V)$ carries  a canonical flat
connection
$\nabla^\varepsilon$ relative to
$\Spec R$ called {\it the
$\varepsilon$-connection}.

Our problem is local with respect to flat
topology of $\Spec R$, so we can assume that
$F_R$ is isomorphic to a product of several
copies of
$R((t))$. Consider the super extension
$F_{(V)}^{\times\flat}$ of $F^\times$
(see 3.7). Now the
$F^\times$-action on
$\omega (F)^\times$, $f,\nu \mapsto f^{-1}\nu$,
lifts canonically to an
$F_{(V)}^{\times\flat}$-action on
$\sE (V)$. Namely, for
$f\in F^\times$, $\nu
\in
\omega (F)^\times$ one has
$\tau_{f^{-1}\nu}=f\tau_\nu$, hence
$\nabla (\tau_{f^{-1}\nu})= f\nabla
(\tau_\nu )$. Our action is the product map
\eq{4.5.1}{\det (V^\infty
,f^\infty )\cdot\det (V^\infty ,
\nabla (\tau_{\nu})^\infty )\iso \det
(V^\infty ,\nabla (\tau_{f^{-1}\nu})^\infty ).}

So our $\sE$ fits into the setting of 3.10. Let
$\fd :=$tr$\nabla$ be the connection on
$\det V$ defined by
$\nabla$. Our
$\nabla^\varepsilon$ is the connection
$\nabla^\fd$ from 3.10.

{\it Remark.}
$\nabla^\varepsilon$ is functorial with
respect to isomorphisms of $(F,V,\nabla )$, and
compatible with base change and disjoint
sum identification (4.4.2).

\medskip

  Consider the group ind-scheme
$\Aut (F )$ on
$\Spec R$,
$\Aut (F)(R'):=\Aut (F_{R'})$ (see 3.2). Let
$\Aut (F)\,\hat{}$ be its formal
completion; this is a formal group whose Lie
algebra is the Lie algebra of vector fields
$\Theta (F) :=\text{Der} (F_R /R )$. The group
$\Aut (F)$ acts on the ind-scheme
$\omega (F)^\times$.

  This action can be lifted
to an action of the formal completion $\Aut
(F)\,\hat{}$ on $\sE (V)$ in the following two
ways:

(a) Using $\nabla$ one lifts
the action of $\Aut (F)\,\hat{}$ on $F$ to
$V$. This action is compatible with $\nabla$,
so $\Aut (F)\,\hat{}$ acts on $\sE (V)$ by
transport of structure.

(b) The $\varepsilon$-connection
$\nabla^\varepsilon$ lifts the $\Aut
(F)\,\hat{}\,$-action on $\omega (F)^\times$ to
$\sE (V)$.

\subsection{Proposition} The above two $\Aut
(F)\,\hat{}\,$-actions on $\sE (V)$ coincide.

\begin{proof}
It suffices to check that the two actions of
the Lie algebra $\Theta (F)$ on $\sE (V)$
coincide. For $\theta\in\Theta (F)$ let
$\theta^{(a)}, \theta^{(b)}$ be the
corresponding actions of $\theta$ on $\sE
(V)$. Then $\kappa (\theta ):=  \theta^{(a)}-
\theta^{(b)}$ is an $\sO$-linear endomorphism
of $\sE (V)$, i.e., a function on $\omega
(F)^\times$. We want to show that $\kappa
:\Theta (F) \to \sO (\omega (F)^\times )$
vanishes.

Action
(a) preserves $\nabla^\varepsilon$ (see Remark
  in 4.5), so for every
$\theta,\theta' \in \Theta (F)$ one has
$[\theta^{(a)},\theta^{'(b)}]=
[\theta,\theta']^{(b)}$ which is
$[\theta^{(b)},\theta^{'(b)}]$. Hence
$\theta '(\kappa (\theta ))=0$, i.e., the
image of $\kappa$ belongs to the subspace of
$\Theta (F)$-invariant functions.

Therefore the two actions coincide on the
commutator $[\Theta (F), \Theta (F)]$. We are
done since
$\Theta (F)$ is a perfect topological Lie
algebra.
\end{proof}

{\it Remark.} Using definition (a) we see
that the
$\Aut (F)\,\hat{}\,$-action on $\sE (V)$ is
compatible with the $\hat{F}^{\times
\flat}_{(V)}$-action from 4.5, i.e., $\sE
(V)$ carries a natural action of the
semi-direct product of $\Aut
(F)\,\hat{}\,$ and $\hat{F}^{\times
\flat}_{(V)}$.\footnote{$\Aut
(F)\,\hat{}\,$ acts on $\hat{F}^{\times
\flat}_{(V)}$ since it acts on $V$ via
$\nabla$.}

\subsection{More on the
$\Theta (F)$-action on $\sE (V)$}  First let us
discuss the
$\Theta (F)$-action on
$\omega (F)^\times$. Consider a formally smooth
morphism
  \eq{4.7.1}{\fr:
\omega (F)^\times \to \A^1_R , \quad \fr (\nu
):=\Res (\nu ),}  invariant with
respect to the
$\Aut (F)$-action.

Take any $\nu \in \omega (F)^\times$. The
tangent space to $\nu$ equals $F$.\footnote{
$a\in F$ corresponds to a derivative $\phi
\mapsto
\partial_\epsilon \phi ((1+\epsilon a)\nu )$.}
  So the
$\Theta (F)$-action yields a morphism
$\zeta_\nu : \Theta (F )\to F$.

\medskip

{\bf Lemma.} (i) The image of $\zeta_\nu$
equals
$\tau_\nu (F)\subset F$. It coincides with the
tangent space to the fiber of $\fr$ at $\nu$,
hence
$\Aut
(F)\hat{}\,$ acts formally transitively along
the fibers of $\fr$.

(ii) The kernel of $\zeta_\nu$, i.e., the
stabilizer of $\nu$ in $\Theta (F)$, is
$\sO \tau_\nu \subset\Theta (F)$.

\begin{proof}
The differential
to $\fr$ at $\nu$ equals the
functional
$r_\nu : F\to R$, $f\mapsto \Res (f\nu )$, so
the tangent space to the fiber of $\fr$ at
$\nu$ equals $\Ker r_\nu$. For
$\theta\in\Theta (F)$ one has
$\zeta_\nu (\theta )=\sL ie_\theta (\nu)
/\nu =-d(\nu \theta )/\nu$. Therefore
the multiplication by $\nu$ isomorphism
identifies $\zeta_\nu$ with the de
Rham differential
$F
\to\omega (F)$. This implies our assertions.
\end{proof}

Let us describe the
  $\Theta (F)$-action on $\sE (V)$
directly in terms of the
$F^{\flat}_{(V)}$-action on $\sE$ (see
(4.5.1)). Since
$\sE (V)$ is a $F^{\times
\flat}_{(V)}$-torsor, the $\Theta (F)$-action
can be considered as a rule that assigns to
$\nu \in \omega (F)^\times$ a lifting
$\zeta_\nu^\flat : \Theta (F )\to
F^{\flat}_{(V)}$ of
$\zeta_\nu :\Theta (F )\to F$.

\medskip

{\bf Lemma.} The kernel of $\zeta_\nu^\flat$
coincides with the kernel of $\zeta_\nu$, and
its image coincides with image of $\tau_\nu $
acting on $ F^{\flat}_{(V)}$. Thus
$\zeta_\nu^\flat$ is the composition of
$\zeta_\nu$ and the inverse to the isomorphism
$\tau_\nu (F^{\flat}_{(V)}) \iso \tau_\nu
(F)$.

\begin{proof} Consider the section
$\nabla^\varepsilon_\nu :F\to
F^{\flat}_{(V)}$. According to 4.6 one
has
$\zeta_\nu^\flat =\nabla^\varepsilon_\nu
\zeta_\nu$, so $\Ker \zeta^\flat_\nu
=\Ker \zeta_\nu$. Since
$\nabla_\nu^\varepsilon$ is
$\tau_\nu$-invariant (see Remark in 4.6) and
$\tau_\nu $ kills $\sO \subset
F^{\flat}_{(V)}$, one has $\tau_\nu
(F^{\flat}_{(V)})=\nabla^\varepsilon_\nu
(\tau_\nu (F))$. Since
$\zeta_\nu (\Theta (F) )=\tau_\nu (F)$ by the
previous Lemma, we are done.
\end{proof}

\subsection{The absolute
$\varepsilon$-connection} Suppose that $V$
carries an absolute flat connection
$\nabla^{abs} $ that extends $\nabla$. Let us
show that the super line bundle $\sE (V,\nabla
)$ on
$\omega (F)^\times$ acquires then an absolute
flat connection  which extends the
relative
$\varepsilon$-connection.

So let $T$ be a test $k$-algebra, $I\subset T$
a nilpotent ideal, and $t,t' \in
\omega (F)^\times (T)$ two  $T$-points that
coincide mod $I$. We want to define a
canonical isomorphism \eq{4.8.1}{\sE
(V,\nabla )_t
\iso \sE (V,\nabla )_{t'}} which is the
  identity  mod $I$, satisfies a
  transitivity property, and is compatible
with the base change (see \cite{Gr}).

Let $t_R ,t'_R :R\to T$ be the
corresponding $T$-points of $\Spec R$. Let
$F_T := F\hat{\otimes}_{t_R} T$, $V_T$ ,
$\nabla_T$ be the
$t_R$-pull-back of $F$, $V$, $\nabla $;
this is an $\varepsilon$-nice $T$-family.
Define
$F'_T$, $V'_T$, $\nabla'_T$ using $t'_R$
instead of $t_R$. The two families coincide
mod $I$. Since $F$ is formally smooth the
identification $F_{T/I}=F'_{T/I}$ can be
extended to an isomorphism of topological
$T$-algebras $\phi : F_T
\iso F'_T$. Now the absolute connection
$\nabla^{abs}$ lifts $\phi$ to an isomorphism
$\phi_V : V_T \iso V'_T$ which extends the
identity isomorphism $V_{T/I}=V'_{T/I}$. It is
automatically compatible with $\nabla_T$,
$\nabla'_T$.

Let $\phi_\omega :\omega^\times_{F_T} \iso
\omega^\times_{F'_T}$ be the isomorphism of
$T$-ind-schemes defined by $\phi$. We have the
$\phi_\omega$-isomorphism of
$\varepsilon$-lines defined by
$\phi_V$
\eq{4.8.2}{\phi_V^\varepsilon :\sE (V_T
,\nabla_T )\iso
\sE (V'_T ,\nabla'_T )}
which extends the identity isomorphism mod
$I$ and is compatible with the
$\varepsilon$-connections. We can consider
$t$, $t'$ as $T$-points of
$\omega_{F_T}^\times$,
$\omega_{F'_T}^\times$. They coincide
mod $I$, so $\phi_V^\varepsilon$ and
$\nabla^\varepsilon$ yield an identification
\eq{4.8.3}{ \sE (V_T ,\nabla_T )_t \iso \sE
(V'_T ,\nabla'_T )_{\phi_\omega (t)}
\iso
\sE (V'_T ,\nabla'_T )_{t'}.} The super
lines from (4.8.3) equal the super lines from
(4.8.1) by base change, so we have defined the
promised isomorphism (4.8.1). It does not
depend on the auxiliary choice of $\phi$:
indeed, various $\phi$ are
$\Aut (F_T
)\hat{}\,$-conjugate, so the assertion follows
from 4.6. The transitivity and base change
properties are clear.

\subsection{The standard isomorphisms}
Here is a list (cf.~1.1):

\medskip

(i) {\it Direct sums.}  Let $F_R$ be as in 4.5,
and
$(V_\alpha ,\nabla_\alpha )$ is a  finite
collection of
$F$-modules with relative connection as in 4.5
such that every every $(V_\alpha
,\nabla_\alpha )$ is $\varepsilon$-nice in the
sense of loc.~cit. Then $(V,\nabla ):=(\oplus
V_\alpha ,\oplus \nabla_\alpha )$ is
$\varepsilon$-nice, and, by (2.6.3), we have  a
canonical isomorphism of $\Z$-graded super
lines
\eq{4.9.1}{\sE (V,\nabla )=\otimes \sE
(V_\alpha ,\nabla_\alpha )} on
$\omega (F)^\times$ compatible with
$\varepsilon$-connections. This isomorphism is
compatible with the constraints, so $\sE$ is
a symmetric monoidal functor from the
category of $(V,\nabla )$ with respect to
direct sums to that of super lines with
connection on
$\omega (F)^\times$. If
our $\nabla$'s come from absolute connections
then (4.9.1) is compatible with the absolute
$\varepsilon$-connections.

\medskip

(ii) {\it Filtrations.} Let $(F_R ,V,\nabla )$
be an
$\varepsilon$-nice family and assume we have a
filtration on $V$ compatible with $\nabla$ such
that $(F_R ,\text{gr}V, \text{gr}\nabla )$ is
$\varepsilon$-nice.\footnote{Probably, the
latter condition is automatic.} Then (2.9.1)
yields a canonical isomorphism \eq{4.9.2}{\sE
(V,\nabla )=\sE (\gr V,\gr \nabla )}
compatible with  $\varepsilon$-connections.
If
our $\nabla$'s come from absolute connections
then (4.9.2) is compatible with the absolute
$\varepsilon$-connections.

\medskip

(iii) {\it Induction.} Let $(F'_R ,V',\nabla'
)$ be an $\varepsilon$-nice family, $F_R$
another topological $R$-algebra as in 4.5, and
$F_R \to F'_R$ a morphism of topological
$R$-algebras. It is automatically \'etale.
Denote by $V$ our $V'$ considered as an
$F$-module; it carries the induced
connection $\nabla$. One has an
embedding $\omega (F)^\times \hra
\omega (F')^\times$. Now there is a canonical
isomorphism\footnote{Take any $\nu \in
\omega (F)^\times \subset \omega (F')^\times$;
let $\tau_\nu \in \text{Der}(F/R)$, $\tau'_\nu
\in \text{Der} (F'/R)$ be the corresponding
vector fields. Then
$\nabla (\tau_\nu )=\nabla' (\tau'_\nu )$,
hence the $\varepsilon$-lines are the same.}
\eq{4.9.3}{\sE (V',\nabla'
)|_{\omega (F)^\times} \iso \sE (V,\nabla)}
compatible with
$\varepsilon$-connections. If $\nabla'$ comes
from an absolute connection, then it defines
an absolute connection on $V$, and (4.9.3)
is compatible with the absolute
$\varepsilon$-connections.

\medskip

(iv) {\it Non-singular situation.} Assume
that $F_R$ is isomorphic to $R((t))$ and
$(V,\nabla )$ is non-singular, i.e.,
$F\otimes_R V^\nabla \iso V$. Let $\omega (F)^0
\subset \omega (F)^\times$ be the $0^{\rm th}$
connected component, i.e., the component of
$dt$. Then there is a canonical isomorphism of
$\Z$-graded super lines
\eq{4.9.4}{\sE (V,\nabla )|_{\omega (F)^0}
=(\det V^\nabla )_{\omega (F)^0 }} compatible
with the relative connections (the
$\varepsilon$- and the trivial one
respectively). If $\nabla$ comes from
an absolute connection then (4.9.4) is
compatible with the absolute connections.

To establish (4.9.4) we choose an isomorphism
$R((t))\iso F_R$; let $O\subset F$ be the
image of $R[[t]]$. Set $V(O ):= O\otimes_R
V^\nabla \subset V$. Then $\nabla
(\partial_t)$ preserves $V(O)$, the
induced operator $\nabla (\partial_t )$ on
$V/V(O)$ is injective, and
$V^\nabla \iso
t^{-1}V(O )/V(O) \iso \Coker (\partial_t
:V/V(O)
\to V/V(O) )$ where the left arrow is
multiplication by $t^{-1}$. The composition
yields
\eq{4.9.5}{\sE (V,\nabla )_{dt} =\det (V/V(O)
,\nabla (\partial_t )^\infty )\iso \det
V^\nabla .}

Our connection comes\footnote{Locally on
$\Spec R$.} from a trivialization of $V$.  By
  3.9(i) and discussion at the end of
3.5\footnote{Notice that $\omega
(F)^0_{red}=\omega (O)_{red}$.} the connection
$\nabla^\varepsilon$ is ``constant" along the
fibers of
$\omega (F)^0 \to \Spec R$. We
define (4.9.4) as the horizontal morphism
which is (4.9.5) at
$dt$.\footnote{The fibers of $\omega (F)^0 \to
\Spec R$ are connected, so this defines
(4.9.4) uniquely.}

Let us check that our isomorphism does not
depend on the auxiliary choice of $t$. Since
$\nabla$ is non-singular, the action of $\Aut
(F)$ on $F$ lifts to an action on $V=
F\otimes_R V^\nabla$ which preserves $\nabla$.
Thus $\Aut (F)$ acts by transport of
structure on $\sE (V,\nabla )$. This action
preserves $\nabla^\varepsilon$, and the
induced action on $\sE (V,\nabla
)^{\nabla^\varepsilon}_{\omega (F)^0}$ is
trivial by 4.6.\footnote{Recall that $\Aut
(F)$ is connected.} Since $\Aut (F)$ acts
transitively on the set of $t$'s, and the
action of $g\in \Aut (F)$ on $\sE
(V)$ sends the isomorphism (4.9.4) defined
by means of $t$ to that defined by means
of $g(t)$, we are done.

\medskip

(v) {\it The product formula}.
Assume we are in a global situation as
discussed in
  3.11, so we have a family $X$ of smooth
projective curves over $\Spec R$, a relative
divisor $D\subset X$ such that the
projection $D\to \Spec R$ is \'etale and
surjective, and a vector bundle $V$ on $U:=
X\smallsetminus D$. Let
$\nabla$ be a (relative) connection on $V$.
Consider the morphism
$\nabla : V(U)
\to (\omega\otimes V)(U)$ of projective
$R$-modules.\footnote{The $R$-projectivity
follows since $\sO (U)$ is a
projective $R$-module and $V(U)$ is a
projective $\sO (U)$-module.} We say that our
situation is {\it
$\varepsilon$-nice} if
$\nabla$ is a Fredholm morphism. Then
one has  the corresponding
$\Z$-graded super $R$-line
\eq{4.9.6}{\det (R\Gamma_{dR}(U, V)[1]) :=
\det ((\omega
\otimes V)(U), V(U),\nabla^\infty ).}

Since $V(U)\subset V(F)$ is a d-lattice, i.e.,
$V(U)^\infty \iso V(F)^\infty$, we see that
$(U,V,\nabla )$ is $\varepsilon$-nice if and
only if $(F,V(F) ,\nabla )$ is
$\varepsilon$-nice in the sense of 4.4. Below
we assume this.

\subsection{Proposition}
There is a canonical isomorphism of
$\Z$-graded super lines on $\omega (U)^\times$
(the   product formula)
\eq{4.10.1}{\sE (V(F) ,\nabla )|_{\omega
(U)^\times}\iso\det (R\Gamma_{dR}(U,V
)[1])_{\omega (U)^\times}} where the r.h.s.~is
the pull-back of the determinant
$R$-line to $\omega (U)^\times$. This
isomorphism is compatible with the
relative connections (the l.h.s.~carries the
$\varepsilon
$-connection, the r.h.s.~the trivial one).

  If $V$ carries an absolute integrable
connection which extends $\nabla$ then (4.10.1)
is compatible with the corresponding absolute
connections (the $\varepsilon$- and the
Gau\ss-Manin one). Thus (4.10.1)
describes the determinant of the Gauss-Manin
connection in terms of the absolute
$\varepsilon$-connections.

\begin{proof}
A trivialization $\nu \in\omega
(U)^\times$  identifies $\nabla
:V(U)\to (\omega \otimes V)(U)$ with $\nabla
(\tau_\nu ): V(U)\to V(U)$ where
$\tau_\nu :=\nu^{-1} \in \text{Der}(U/R)$, so
we get
$\det (
R\Gamma_{dR}(U,V)[1])_\nu = \det
(V(U),V(U),\nabla (\tau_\nu )^\infty )$. The
latter super line is identified with $\sE
(V(F),\nabla )_\nu$ via the asymptotic
isomorphism
$V(U)^\infty \iso V(F)^\infty$.

The compatibility with relative connections
follows from 3.12. It implies immediately the
compatibility with absolute connections.
\end{proof}

\subsection{Remark}  The above
standard isomorphisms are mutually compatible
in the obvious sense; they are compatible with
the disjoint sum identifications (4.4.2). In
particular, if $D$ is a disjoint union of
sections $x_i$ then $F=\Pi F_i$ (the
product of local fields at $x_i$) and $\sE (V
(F),\nabla )_\nu = \otimes\sE (V(F_i
),\nabla )_\nu $, hence the name
``product formula" for (4.10.1).

\subsection{Duality} Suppose that we have a
family
$(F,V,\nabla )$ such that
$H^\cdot_{dR}(V,\nabla )=0$, i.e., $\nabla : V
\to
\omega \otimes V$ is an isomorphism.
Let $V'$ be the
dual vector bundle  equipped with the dual
connection
$\nabla'$. It satisfies the
same property. For $\nu \in\omega (F)^\times$
the Tate $R$-module $V'$ identifies with the
dual to $V$ via the pairing $v,v' \mapsto \Res
(v'(v)\nu )$, and $\nabla' (\tau_\nu )$ is
the operator adjoint to $-\nabla (\tau_\nu
)=\nabla (\tau_{-\nu})$. So (2.13.5) yields a
canonical identification of super lines
\eq{4.11.1}{\sE (V' ,\nabla'
)_\nu =\sE (V,\nabla )^{-1}_{-\nu}.}  We leave
it to the reader to check that it is compatible
with the $\varepsilon$-connections. 

Suppose
now we are in the global situation of
4.9(v) and the  local de Rham
cohomology vanish, i.e., $\nabla : V(F)\iso
(\omega\otimes V)(F)$. Thus $H^0_{dR}
(U,V)=H^0_{dR} (U,V')=0$ and the global
Poincar\'e duality identifies $H^1_{dR}(U,V')$
with the dual to $H^1_{dR}(U,V)$. Passing to
determinants we get a canonical identification
$\det R\Gamma_{dR}(U, V' )=
\det R\Gamma_{dR}(U, V)^{-1}$. The
product formula is compatible with these
identifications.

{\it Remark.} We have formulared 
compatibility (4.12.1) in the situation when
the local de Rham cohomology
$H^\cdot_{dR}(V,F)$ vanish. If the
base ring $R$ is a field (or, more generally,
an Artinian algebra) then (4.12.1) holds for
arbitrary $(V,\nabla )$ by (2.13.6) since the
local super line $\det (V,\nabla (\tau_\nu )
)=\det H^\cdot_{dR}(V)$ is canonically
trivialized (see 5.9(iv)).\footnote{Probably,
this is true for arbitrary $R$ under the
assumption that both
$\nabla$ and $\nabla'$ are $\varepsilon$-nice.
It also looks probable that $\nabla$ is
$\varepsilon$-nice iff such is $\nabla'$.}

\subsection{Questions} (a) What
happens when  singular points merge together,
i.e., we do not assume that $D$ is \'etale
over $\Spec R$?

(b)  Is it
true that
$(\sE ,\nabla^\varepsilon )$ is uniquely
determined by compatibilities (4.4.2),
(4.9.1)--(4.9.4), (4.10.1)?\footnote{One can
show that the connection on our $\sE$
(defined in (4.4.1)) is uniquely determined by
the condition that it is compatible with the
standard isomorphisms.}

\section{Some formula (with apologies to
Charlotte){\protect \footnote{
The authors refer to ``Charlotte's Web", an
acclaimed novel by E.~B.~White ruminating over
the problems of a person in quest for
tenure.} }}

In this section we assume, for technical
reasons, that our parameter
space\footnote{This is
$\Spec R$ of the previous section.} is
$\Spec K$ where
$K$ is a field. The case of an arbitrary
smooth parameter space reduces to this
situation (since a line bundle with connection
on a smooth variety is determined by its
restriction to the generic point).

\subsection{Dependence on $\nu$} Let $F$ be  a
local
$K$-field,
$V_F$ an $F$-vector space of dimension $n$
equipped with an absolute connection $\nabla$.
We have the
$\varepsilon$-line $\sE =\sE (F,V_F )$; this is
a super line with an absolute
integrable connection
$\nabla^\varepsilon$ on the ind-scheme of
invertible relative 1-forms
$\omega (F)^\times$. Our aim is to compute its
fiber $\sE_\nu =\sE (F,V_F )_\nu$ over a
$K$-point
$\nu$ of $\omega (F)^\times$ which is a super
$K$-line equipped with an integrable
connection.

When doing computations it is convenient to
choose $\nu$ in a special way, e.g.~to make
it fixed by horizontal vector fields.  One can
pass then to arbitrary
$\nu$ using Proposition below.

\medskip

Choose a parameter $t\in F$. Then $F=K' ((t))$
where a finite extension $K'/K$ is the integral
closure of $K$ in $F$. Set $O:=K'[[t]]\subset
F$. Below we make vector fields on $K$ act on
$F$ so that for $\theta \in \Theta_K =$Der$K/k$
one has
$\theta (t)=0$. Hence $\Theta_K$ acts on
$\omega (F)^\times$. Differential 1-forms and
vector fields on $\Spec F$ and $\omega
(F)^\times$ decompose into ``horizontal" and
``vertical" components accordingly. For
1-forms on $\Spec F$ these components are
denoted by lower indices
$x$ and $t$. E.g., for
$f\in F$ we write $df=d_x f+d_t f$ where $d_x
f\in F\otimes
\Omega^1_{K/k}$, $d_t f\in \omega
(F)=\Omega^1_{F/K}$.

\medskip

Choose an $O$-lattice $V_O\subset V_F$. It
yields  a section
$s^c_{V_O}:O^\times
\to F^{\times\flat}_{(V_F )}$ (see (2.13.3),
3.7). Since
$\sE$ is an
$F^{\times \flat}_{(V_F )}$-equivariant line
bundle (see 4.5), $s^c_{V_O}$ defines an
$O^\times$-equivariant structure on
$\sE$.\footnote{Recall that the
$O^\times$-action on $\sE$ lifts the
$O^\times$-action
$u,\nu
\mapsto u^{-1}\nu$ on $\omega (F)^\times$.}

Assume we have $\nu ,\nu' \in \omega (F)(K)$
such that $v(\nu )=v(\nu')=:\ell$. Then
$u:=\nu'/\nu
\in O^\times$, i.e., $u=u_0 +u_1 t+..$, $u_i
\in K'$, $u_0 \neq 0$. So the equivariant
structure provides an isomorphism
\eq{5.1.1}{\alpha^u
=\alpha^u_{V_O}:=s^c_{V_O}(u):
\sE_{\nu'}\iso \sE_{\nu} .}
We want to see to what extent $\alpha^u$ is
horizontal with respect to connections
$\nabla^\varepsilon$, i.e., to compute
$\nabla^\varepsilon \log
\alpha^u := (\nabla^\varepsilon \alpha^u
)/\alpha^u
\in
\Omega^1_{K/k}$.

\medskip

Set $\chi =\chi_x +\chi_t := (\text{tr}\nabla
)\log
\gamma\in
\Omega^1_{F/k}
$ where $\gamma$ is any generator of $\det
V_O$. This is a closed 1-form; its class
modulo $\Omega^1_{O/k}$ does not depend on the
choice of $\gamma$.

\subsection{Proposition} One has
\eq{5.2.1}{\nabla^\varepsilon \log\alpha^u
  =tr_{K'/K} (\frac{n\ell}{2}d\log u_0 +\Res
(d\log u \wedge \chi ) ).}

\begin{proof} By base change we can assume
that $K'=K$.

Take any  $a\in K[[t]]=O\subset F$, $a=a_0 +a_1
t +..$. Consider $a$ as a vertical vector
field on $\omega (F)^\times$ (the corresponding
derivative is
$\phi
\mapsto a\phi$, $a\phi (\nu ):=
  \partial_\epsilon \phi ((1+\epsilon a)\nu)$).
It acts on $\sE$ according to
$\nabla^\varepsilon$: namely,
$\nabla^\varepsilon (a)$ is the action of
$\nabla^\fd (a) \in F^\flat_{(V)}$ on $\sE$
(see 4.5). Our
$a$ also acts on
$\sE$ according to the $O^\times$-equivariant
structure, i.e., as $s^c_{V(O)} (a)\in
F^\flat_{(V)}$.

 From now on we work on the connected component
$\omega (F)^\ell_{red}$ of the reduced scheme
$\omega (F)^\times_{red}$ (the
action of $a$ preserves $\omega
(F)^\ell_{red}$). By  (3.10.3)
and 4.5 one has\footnote{Here the l.h.s.~is
a (constant) function on $\omega
(F)^\ell_{red}$, in the r.h.s.~$a$ is
considered as an element of $F$.}
  \eq{5.2.2}{\nabla^\varepsilon
(a)-s^c_{V(O)}(a)= -\frac{n\ell}{2}
a_0 -
\Res (a \chi_t ).}

Let $r$ be a non-zero $O^\times$-invariant
section of $\sE$.  Then $r_{\nu}$,
$r_{\nu'}$ are trivializations of,
respectively,
$\sE_{\nu}$ and $\sE_{\nu'}$ which are
identified by the isomorphism $\alpha^u$. Set
$\psi := \nabla^\varepsilon \log r$. One has
$\nabla^\varepsilon \log \alpha^u
=\psi_{\nu}-\psi_{\nu'}$.

Our $\psi$ is a closed 1-form on the
$O^\times$-torsor $\omega (F)_{red}^\ell$. By
(5.2.2) the vertical part of $\psi$ is an
$O^\times$-invariant form whose value on
the vertical vector field $a\in O$ equals
$-\frac{n\ell}{2}  a_0 -
\Res (a\chi_t )$. If $\psi'$ is any other
1-form with these properties then $\psi'-\psi$
is the pull-back of a 1-form on $\Spec K$,
hence $\nabla^\varepsilon \log \alpha^u
=\psi'_{\nu}-\psi'_{\nu'}$.

Let us construct such $\psi'$. Fix
$\nu'$ and consider $\nu$ as a variable point
of $\omega (F)_{red}^\ell$. The
scheme
$\omega (F)_{red}^\ell \hat{\otimes}_K O$
carries a canonical invertible function $u$:
for $\nu \in\omega (F)_{red}^\ell$ the
restriction of $u$ to the fiber over $\nu$ is
the function $\nu' /\nu$. The pull-back of
$\chi$ from $F$ to $\omega (F)_{red}^\ell
\hat{\otimes}_K F$ is a closed 1-form
which we denote also by
$\chi$. Set $\psi' :=
\frac{n\ell}{2}d\log u_0 +\Res
(d\log u \wedge \chi )$ where $\Res$
is the residue along the fibers of  $\omega
(F)_{red}^\ell
\hat{\otimes}_K F \to \omega (F)_{red}^\ell$.

  Incerting $\nu$, $\nu'$ into
$\psi'$ we get (5.2.1).
\end{proof}

\subsection{Case of regular singularity}
  We  use  notation from 5.1.

Suppose that (the vertical
part of)
$\nabla$ has regular singularities. Then one
can find an $O$-lattice $V_O \subset V_F$ such
that  $V_O$ is preserved by the action of
$\nabla (t\partial_t )$ and horizontal vector
fields $\Theta_K$.

Consider a $K'$-vector space $V_0 :=
V_O /tV_O$. The action of horizontal vector
fields define a connection $\nabla_0$ on
$V_0$. The action of $\nabla (t\partial_t )$
defines a $\nabla_0$-horizontal
endomorphism of $V_0$ which we denote by
$\kappa$.\footnote{So $\kappa$ is the residue
of the polar part of $\nabla$ with respect to
$V_O$.}

\medskip

{\it Remark.}  Recall a standard way to
construct such $V_O$. Let $\bar{K'}$ be an
algebraic closure of $K'$. Choose any
(set-theoretic) section $s:\bar{K'}/\Z \to
\bar{K'}$ that
commutes with the Galois action. Let $V_{K'}$
be the sum of generalized $s(\lambda
)$-eigenspaces of
$\nabla (t\partial_t )$, $\lambda \in
\bar{K'}/\Z$. Then $V_{K'}$ is a $K'$-structure
on
$V_F$. Our lattice is
$V_O :=V_{K'}[[t]]$.

\medskip

Consider a Fredholm endomorphism $\nabla
(t\partial_t )$ of $V_F /V_O$. Its
determinant line is canonically trivialized.
Namely, $V_F /V_O$ decomposes into a direct
sum of generalized eigenspaces of $\nabla
(t\partial_t )$; for each eigenvalue
$\lambda$ the corresponding subspace
$P_\lambda$ is finite-dimensional. So the
embedding $P_0 \hra V_F /V_O$ induces an
isomorphism on cohomology hence  an
isomorphism of the determinant lines. And the
determinant line for $P_0$ is trivialized
since $\dim P_0 <\infty$.

\medskip

Take any $\nu \in \omega (F)^\times (K)$,
$v(\nu )=:\ell$. So $\nu = ut^\ell dt$, $u\in
K'[[t]]^\times$, $u=u_0 +u_1 t+..$, $u_i \in
K'$, $u_0 \neq 0$. We want to compute
$\sE_\nu =\sE (F,V)_\nu$ in terms of $V_O$.

There is a natural isomorphism of super
lines
\eq{5.3.1}{\beta_\nu :\sE_\nu \iso (\det{}_K
V_0 )^{\otimes \ell +1}.} Namely, recall that
$\sE_\nu = \det (V^\infty_F ,\nabla
(\tau_\nu )^\infty  )=\det
(V_F /V_O ,\nabla
(\tau_\nu )^\infty  )$ where $\tau_\nu
:=\nu^{-1}= u^{-1}t^{-\ell}\partial_t $. Let us
present $\nabla (\tau_\nu )^\infty$ as a
composition $V_F /V_O
\buildrel{\nabla(t\partial_t
)}\over\longrightarrow V_F /V_O
\buildrel{t^{-\ell -1}}\over\longrightarrow V_F
/t^{-\ell -1}V_O
\buildrel{u^{-1}}\over\longrightarrow V_F
/t^{-\ell-1}V_O
\buildrel{id^\infty}\over\longrightarrow V_F
/V_O$. We already trivialized the determinant
line of the first operator; that of the
middle operators is trivialized since they are
isomorphisms. Thus $\sE_\nu = \det (V_F
/V_O ,V_F /t^{-\ell-1}V_O , id^\infty )=(\det
V_0 )^{\otimes \ell+1}.$\footnote{The latter
identification comes from isomorphisms
$t^{-i} :t^i V_O /t^{i+1}V_O \iso V_0$.}

The r.h.s.~of (5.3.1) carries a connection
defined by $\nabla_0$. It yields, via
$\beta_\nu$, a connection on $\sE_\nu$ which
we also denote  by
$\nabla_0$.

\subsection{Proposition}
$\nabla^\varepsilon -\nabla_0 =
-tr_{K'/K}( \Tr\kappa -n(\ell/2 +1))d\log u_0
).$

\begin{proof} For $\nu =t^\ell dt$, i.e.,
$u=1$, the formula is clear. Indeed, any
horizontal vector field $\theta$ commutes with
$\tau_\nu$, so
$\nabla (\theta )$ acts on the determinant
line of $\tau_\nu$ by transport of structure,
and
$\nabla^\varepsilon (\theta )$ coincides with
this action (see 4.8).

For arbitrary $\nu$ notice that the isomorphism
$\beta_\nu^{-1}\beta_{t^\ell dt}:\sE_{t^\ell
dt}\iso
\sE_\nu$ is inverse to the isomorphism
$\alpha^{u}$ of (5.1.1) defined by
the lattice $t^{-\ell-1}V_O$. The
corresponding form $\chi$ from 5.1 is $(\Tr
\kappa -n(\ell +1))dt/t$. Now use (5.2.1).
\end{proof}

\medskip

{\it Remark.} Choose a $k$-structure
$V_k \subset V_F$ such that $V_O $ is the
$O$-module generated by $V_k$.
Let
$\nabla^0$ be the absolute connection on $V_F
, V_{K'}$ defined by
$V_k$, $\nabla^0 (V_k )=0$. Then $\nabla
-\nabla^0 =\eta + t^{-1}g dt$ where $\eta \in
\Omega^1_{K'/k}\otimes \text{End}
V_{K'}[[t]]$, $g\in \text{End} V_{K'}[[t]]$.
Write $\eta =\eta_0 +\eta_1 t +..$, $g= g_0 +
g_1 t +..$, $\eta_i \in
\Omega^1_{K'/k}\otimes \text{End} V_{K'}$, $g_i
\in \text{End} V_{K'}$. Then $V_0 =V_{K'}$,
$\nabla_0 =\nabla^0 +\eta_0$, $\kappa =g_0$,
so we can rewrite formula 5.4 as \eq{5.4.1}{
\nabla^\varepsilon - \nabla^0 =
tr_{K'/K}((\ell+1)\Tr
\eta_0 - (\Tr g_0 - n(\ell/2 +1))d\log u_0 ).}

\subsection{Case of an irregular admissible
connection.} As above, we follow notation from
5.1.

Suppose that there is an $O$-lattice
$V_O \subset V_F$ such that $\nabla$ is {\it
$m$-admissible} with respect to $V_O$ for
some $m\ge 2$ (see \cite{BE3}). This means that
the vertical part of $\nabla$ has pole
of ``exact" order $m$, i.e., $t^m \nabla
(\partial_t ) :V_O \iso V_O$, and the
horizontal part of $\nabla$ has pole of order
$\le m-1$.

Choose a $k$-structure $V_k \subset V_F$ such
that $V_O$ is the $O$-submodule generated by
$V_k$. Let $\nabla^0$ be the corresponding
absolute connection on $V_F$, $\nabla^0 (V_k
)=0$. Set $\nabla -\nabla^0 =\sA =\sA_t +\sA_x
\in \Omega^1_{F/k}\otimes\text{End} V$. Write
  $\sA_t =t^{-m}gdt $,
$\sA_x =t^{-m+1}\eta$. Admissibility
  means that
$g\in \text{Aut} V_{K'}[[t]]$, $\eta \in
\Omega^1_{K'/k}\otimes  \text{End}
V_{K'}[[t]]$. So
$g =g_0 + g_1 t+..$, $\eta =\eta_0 +\eta_1 t
+..$ where $g_i \in  \text{End}
V_{K'}$,
$\eta_i
\in
\Omega^1_{K'/k}\otimes \text{End}V_{K'}$, $g_0$
is invertible.

\medskip

Consider a 1-form $\nu =ut^\ell dt \in \omega
(F)^\times$,
$u\in O^\times$, and the corresponding
  $\sE_\nu =\sE (F,V_F )_\nu$.
We have a natural isomorphism of
super lines \eq{5.5.1}{
\beta_\nu :\sE_\nu \iso (\det{_K}
(V_{K'}))^{\otimes \ell+m}} similar to (5.3.1).
Namely, notice that the operator
$\nabla (\tau_\nu )$ yields an isomorphism
$V_F /V_O \iso  V_F /t^{-\ell-m}V_O$ which
trivializes the super line $\det
(V_F /t^{-\ell-m}V_O , V_F /V_O ,\nabla (\tau
)^\infty )$. Our identification is the
  composition $\sE_\nu
=
\det (V_F /V_O ,V_F /t^{-\ell-m}V_O , id^\infty
)\cdot
\det (V_F /t^{-\ell-m}V_O , V_F /V_O ,\nabla
(\tau )^\infty )=\det (V_F /V_O ,V_F
/t^{-\ell-m}V_O , id^\infty )=(\det{}_K
(V_{K'}) )^{\otimes \ell+m}.$

The $k$-structure $V_k$ on $V_{K'}$ defines an
integrable connection on $V_{K'}$, hence on the
r.h.s.~of (5.5.1). Denote by $\nabla^0_{\sE}$
the corresponding connection on $\sE_\nu$.

\subsection{Theorem} One has
$\nabla^\varepsilon -\nabla^0_\sE
=tr_{K'/K}(\phi_1 +\phi_2 +\phi_3 )$ where
  $\phi_i \in \Omega^1_{K'/K}$ are:

$\phi_1 = \Res\Tr (g^{-1}d g
\wedge
\sA_x ) -\frac{m}{2}d\log\det (g_0 )$

$\phi_2 = (\ell+m)\Res
\Tr (\sA_x \wedge dt/t ) $

$\phi_3 =
  n(\ell/2 +m)d\log
u_0 -\Res (d
\log u \wedge \Tr \sA ).$

\begin{proof} (a) By base change we can assume
that $K'=K$, so $F=K((t))$.

\medskip

(b) It suffices to treat the
case of $\nu =t^\ell dt$. Indeed, the formula
for arbitrary $\nu$ follows then from (5.2.1)
since
$\beta_\nu^{-1}\beta_{t^\ell dt}:\sE_{t^\ell
dt}\iso
\sE_\nu$ is inverse to the isomorphism
$\alpha^{u}$ of (5.1.1) for the lattice
$t^{-\ell-m}V_O$. The form $\chi$ from 5.1
is $\Tr \sA -(\ell +m)ndt/t$. We get the term
$\phi_3$.

\medskip

(c) From now on $\nu =t^\ell dt$. Set $h:=
1+t^m g^{-1}\nabla^0 (\partial_t )\in \GL (V_F
)$. Our $\sE_\nu$ is the determinant line of
the invertible operator
\eq{5.6.1}{\nabla (\tau_\nu )= t^{-\ell-m}gh}
acting on a Tate $K$-vector space $V_F =V_K
((t)).$ Notice that both factors
$g$ and $h$
preserve the lattice $V_O$, so we have
$\tilde{g}:= s^c_{V_O}(g)$,
$\tilde{h}:=s^c_{V_O}(h) \in F_{(V_F
)}^{\times\flat}$ (see (2.13.3)). The operator
$t$ is defined over $k$, i.e., $t\in \GL (V_k
((t)) )\subset
\GL (V_K ((t)))$. Let $\tilde{t}\in \GL^\flat
(V_k ((t)) )\subset \GL^\flat (V_K ((t)) )$ be
any lifting of $t$ which is also defined over
$k$. So we have \eq{5.6.2}{
\tilde{t}^{-\ell-m}\tilde{g}\tilde{h}\in\sE_\nu
.} It follows from the construction of
$\nabla^0_\sE$ that  \eq{5.6.3}{\nabla^0_\sE
(\tilde{t}^{-\ell-m}\tilde{g}\tilde{h})=0.}
It remains to compute
$\phi :=\nabla^\varepsilon \log
(\tilde{t}^{-\ell-m}\tilde{g}\tilde{h})\in\Omega^1_{K/k}
$.

\medskip

(d) Since $\tau_\nu$ is preserved by the
action of horizontal vector fields our
$\nabla^\varepsilon$ coincides with the
$\nabla$-action of horizontal vector fields
(see 4.8).

Explicitly, for
$\theta
\in\Theta_K \subset \Theta_F$ the infinitesimal
automorphism $1+\epsilon \nabla (\theta )$
acts on the Tate $K$-vector space $F=K((t))$
hence on $\GL (V_F )$ and its super extension
$\GL^\flat (F)$. This action fixes  $\nabla
(\tau_\nu )\in \GL (V_F )$; the
  action on its super line
$\sE_\nu \subset \GL^\flat (V_F )$ is
$r\mapsto r+\epsilon\nabla^\varepsilon (\theta
)r$.

Notice that the $\nabla$-action of horizontal
vector fields preserves, in fact, both $t$ and
$gh
\in \GL (V_F )$. So $(1+\epsilon \nabla
(\theta ))
\tilde{t}^{-\ell-m}=(1+\epsilon \phi_2 (\theta
))\tilde{t}^{-\ell-m}$ and
$(1+\epsilon \nabla (\theta ))
(\tilde{g}\tilde{h})=(1+\epsilon \phi_1 (\theta
))\tilde{g}\tilde{h}$ for some $\phi_i (\theta
)\in K\subset \fgl (V_F )$. Thus
\eq{5.6.4}{\phi = \phi_1 +\phi_2 .}

We will compute $\phi_i$
separately:

\medskip

(e) Since $\nabla^0 (\theta )\tilde{t}=0$  one
has
$1+\epsilon
\phi_2 (\theta )=\Ad_{1+\epsilon\sA_x (\theta
)} (\tilde{t}^{-\ell-m})/\tilde{t}^{-\ell-m}$
or
  $\phi_2 (\theta )=-(\ell+m)(\sA_x (\theta
)^\flat -\Ad_{\tilde{t}}\sA_x (\theta )^\flat
)$ where $\sA_x (\theta
)^\flat$ is any lifting of $\sA_x (\theta
)\in \text{End}_F (V_F )\subset \fgl (V_F )$ to
$\fgl^\flat (V_F )$. By (3.1.6), 3.3(iii) one
has $\phi_2 (\theta )=(\ell+m)\Res \Tr (\sA_x
(\theta )dt/t)$ or
\eq{5.6.5}{\phi_2 =(\ell+m)\Res \Tr (\sA_x
\wedge dt/t).}

(f) Since $ \nabla^0 (\theta ) (V_O )\subset
V_O
$ the action of $1+\epsilon \nabla^0 (\theta )$
preserves the subgroup $\GL (V_F ,V_O
)\subset\GL (V_F )$ and commutes with
$s_{V_O}^c$ (see (2.13.3)). The adjoint action
of $1+\epsilon
\sA_x (\theta )$ does not have this property.
However it sends both $g$ and $h$ to $\GL (V_F
,V_O )$:  this is clear for $h$, and for
$g$ you  notice that,
  due to integrability of
$\nabla$, one has \eq{5.6.6}{[\sA_x (\theta
),g]= t\partial_t (\eta (\theta ))-(m-1)\eta
(\theta )-\theta (g)\in \text{End}V_K [[t]].}

So $\Ad_{1+\epsilon \sA_x (\theta
)}\tilde{g}=(1+\epsilon \phi' (\theta
))s^c_{V_O}(\Ad_{1+\epsilon \sA_x (\theta
)}g)$ for certain $\phi' (\theta )\in K$. We
have defined $\phi' \in \Omega^1_{K/k}$;
replacing $g$ by $h$ we get $\phi'' \in
\Omega^1_{K/k}$.

We see that $(1+\epsilon\nabla (\theta
))\tilde{g}= (1+\epsilon \phi' (\theta
))s^c_{V_O}(\Ad_{1+\epsilon \nabla (\theta
)}g)$; same for $g$ replaced by $h$ and
$\phi'$ by $\phi''$. The product of
these identities gives
  \eq{5.6.7}{\phi_1
=\phi' +\phi''.}

(g) Let us rewrite
$\phi'$, $\phi''$  in the following
general manner:

Assume we have $a\in \fgl (V_F )$ and
$q\in \GL (V_F )$ such that $q(V_O )=V_O$ and
$[a,q](V_O )\subset V_O$. Set $b:= a- \Ad_q
a$, $b^\flat := a^\flat -\Ad_q a^\flat$ where
$a^\flat \in \fgl^\flat (V_F )$ is any
lifting of $a$. Notice that $b^\flat$ is
well-defined\footnote{i.e., it does not depend
on the choice of $a^\flat$.} and
$b(V_O )\subset V_O$. Set
$c(a,q):=b^\flat
-s^c_{V_O}b \in K\subset \fgl^\flat
(V_F )$. We have
\eq{5.6.8}{\phi' (\theta )= c(\sA_x (\theta
),g),
\quad \phi'' (\theta )=c(\sA_x (\theta ),h).}

The following simple lemma helps to compute
$c$:

{\bf Lemma.}    $c(a,q)=\Tr (V_O
,(a^\sharp -\Ad_q a^\sharp )|_{V_O})$ where
$a^\sharp
\in\fgl (V_F )$ is any endomorphism with open
kernel such that
$(a-a^\sharp )(V_O )\subset V_O$.

  {\it Proof of Lemma.}
Notice that
$c(a,q)$ vanishes if $a(V_O )\subset
V_O$.\footnote{To see this consider $a^\flat
=s^c_{V_O}a$.} So we can assume that
$a=a^\sharp$. Then $b^\flat = s_d (b)$ (see
  the end of 2.13), and we are done since
$s_d b - s^c_{V_O}b =\text{Tr}(V_O ,
b|_{V_O})$ for every $b\in \fgl (V_F )\cap
\fgl (V_F ,V_O )$.
\hfill$\square$

\medskip

To compute $\phi'$, $\phi''$ we apply
the lemma to $q=g,h$ and $a^\sharp =\sA_x \Pi$
where
$\Pi :V_F \to V_F$ is a projector  with
kernel
$t^{m-1}V_K [[t]]$ and image
$t^{m-2}V_K [t^{-1}]$.

\medskip

(h) Let us begin with $\phi'$. One has
$\phi' =
\Tr (V_O ,
\sA_x
\Pi -g\sA_x \Pi g^{-1})=\Tr (V_O ,(\sA_x
-g\sA_x g^{-1})\Pi )+ \Tr (V_O , g\sA_x
(g^{-1}\Pi -
\Pi g^{-1}))$.

The first summand vanishes. Indeed, $B:=\sA_x
-g\sA_x g^{-1} \in \text{End} V_K [[t]];$
our summand is $(m-1)\Tr (B_0 )$.
But $\Tr (B_0 )=(\Tr B)_0$, and $\Tr B =0$.

The
second summand can be rewritten as
$\Tr (V_O ,\sA_x g^{-1}(\Pi g -g\Pi ))=-\Tr
(V_O ,\sA_x g^{-1}(1-\Pi)g\Pi )$. An immediate
calculation identifies it with
$-\Res\Tr_F ( gd_t (\sA_x g^{-1}))=-\Res \Tr_F
(\sA_x g^{-1}d_t g)=\Res\Tr_F (g^{-1}d_t g
\wedge
\sA_x )$.

\medskip

  One has $\phi'' = \Tr (V_O , \sA_x
\Pi -h\sA_x \Pi h^{-1})= -\Tr (V_O
,t^m g^{-1}\partial_t \sA_x
\Pi )= \frac{m(m-1)}{2}\Tr g_0^{-1}\eta_0$. As
follows from (5.6.6)\footnote{Indeed, by
(5.6.6) one has $(\sA_x -
g\sA_x g^{-1})_0 = -((m-1)\eta_0 + d_x g_0
)g_0^{-1} $.  Hence
$\Tr ((m-1)\eta_0 + d_x g_0
)g_0^{-1}) =0$.} it equals
$-\frac{m}{2}\Tr g_0^{-1} dg_0
=-\frac{m}{2}d\log (\det g_0 )$. So, by
(5.6.7),
\eq{5.6.9}{\phi_1 = \Res\Tr_F (g^{-1}d_t g
\wedge
\sA_x ) -\frac{m}{2}d\log (\det g_0 ).}
  We
are done.
\end{proof}

\subsection{Example} Suppose that our
$(V_F ,\nabla )$ is such that the
corresponding relative connection
$\nabla_{/K}$ is a successive extension of a
single irregular relative connection
$(L_F ,\nabla_{/K})$ of rank 1.

Let us show that $(V_F ,\nabla )$ is
admissible. We need to construct an
$O$-lattice $V_O
\subset V_F$ as in 5.5.  Notice that the
relative connection on $L$ can be extended to
an absolute integrable connection so that
$L^{\otimes n}$ is isomorphic to $\det V$ as
lines with  absolute connection. Set $P:=
V\otimes L^{-1}$. As a relative connection our
$P$ is a successive extension of trivial
connections. There is a unique $O$-lattice
$P_O \subset P$ preserved by
$\nabla (t\partial_t )$ and such that
$\nabla (t\partial_t )$ acts on the fiber
$P_0 :=P_O /tP_O$ as a nilpotent matrix. Our
$P_O$ is automatically preserved by horizontal
vector fields. Set $V_O :=
P_O\otimes L_O
\subset V$ where $L_O
\subset L$ is any $O$-lattice.

To check the conditions from 5.5  choose a
$k$-structure in $V_O$. Let $\nabla^0$ be the
corresponding connection, so $\nabla
=\nabla^0 +\sA_t +\sA_x$. Then $\sA_t =
  f dt/t^m +
hdt/t $ where  $f\in
O^\times$,
$h=h_0 +h_1 t +..$, $h_i \in K' \otimes
\text{End}V_k$, $h_0$ is nilpotent, and
$\sA_x = a_{q}t^{q} + a_{q+1}t^{q+1}
+..$,
$a_i \in
\Omega^1_{K'/k}\otimes\text{End}V_k$,
$a_q \neq 0$. We want to show that
$q> -m$. Indeed, if
$q\le -m$ then the integrability condition
implies that
$qa_q +[h_0, a_q]=0$ which
contradicts nilpotency of $h_0$.

\subsection{Lemma} There is an isomorphism of
super lines with connection \eq{5.8.1}{\sE
(F,V)
\iso \sE (F,L)^{\otimes n}} where $L$ is the
$F$-line with an absolute connection defined
above.

\begin{proof} A
straightforward comparison of formula 5.6 for
both parts of (5.8.1) (notice that the term
$\Tr (h_0 a_{1-m})$ for l.h.s.~vanishes since
$a_{1-m}$ is a scalar matrix\footnote{Indeed,
  one has $[h_0 ,
a_{1-m}]= (m-1) a_{1-m} +d_x f_0$ due to
integrability. Since $h_0$ is nilpotent this
implies that
$a_{1-m}=(1-m)^{-1}d_x f_0 $.} and
$h_0$ is nilpotent).
\end{proof}

\subsection{A general digression} In the next
subsection we will show that computation of
$\sE_\nu$ for arbitrary $(V_F ,\nabla )$
can be reduced, in principle, to the situation
of 5.7. Before doing this let us list some
(well-known) general properties of connections
on
$F=K((t))$.  For a $K$-relative connection
$(V,\nabla_{/K} )$ we denote by
$H^\cdot_{dR}(V)$ the de Rham cohomology.

\medskip

(a) For every
$(V,\nabla_{/K} )$ as above
$H^0_{dR} (V )$ is the space of horizontal
morphisms
$F\to V$, and
$H^1_{dR}(V)$ is dual to the space of
horizontal morphisms $V\to F$.

Some immediate corollaries:

(i) If $(V,\nabla_{/K})$ is a
non-trivial irreducible relative connection
then $H^\cdot_{dR} (V)=0$. For the
trivial connection $F=(F,\nabla^0_{/K})$ one
has $H^0_{dR}(F)=H^1_{dR}(F)=F$.

(ii) For arbitrary $(V,\nabla_{/K})$ the Euler
characteristic $\chi_{dR}(V)$ vanishes.

(iii) Since $\text{Ext}^\cdot
(V,V')=H^\cdot_{dR}(\sH om (V,V'))$ we see
that every indecomposable\footnote{i.e., $V$
cannot be represented as a direct sum of
connections of smaller rank.}
$V=(V,\nabla_{/K})$ is isomorphic to
$V^{irr}\otimes F^{(n)}$ where $V^{irr}$ is an
irreducible relative connection and $F^{(n)}$
is the nilpotent Jordan
block of length $n$.\footnote{i.e., $F^{(n)}$
is the connection $\nabla_{/K}$ on
$F^n$ with $\nabla_{/K} -\nabla^0_{/K}
=dt/t \cdot$ the nilpotent Jordan block of
length $n$.}

In particular, every $V=(V,\nabla_{/K})$ admits
a  canonical {\it decomposition by
$\nabla_{/K}$-isotypical
components}

\eq{5.9.1}{V=\oplus V^L .}
This decomposition is
labeled by isomorphism classes of relative
irreducible connections
$L$ and characterized by property that  every
irreducible subquotient of
$V^L$ is isomorphic to $L$.

(iv) The  $\Z$-graded super line $\det
H^\cdot_{dR}(V)$ is canonically
trivialized.\footnote{ The argument below works
also in case when our base ring is an
arbitrary Artinian algebra.}

This trivialization is uniquely
characterized by  two properties:
it is compatible with exact sequences of
$V$, and for irreducible
$V$ it comes from
(i) above. Here is a
construction.

By (iii)
  $V$ splits canonically into a
direct sum $V=V^{un}\oplus V^{nun}$ where
$V^{un}$ is the $\nabla_{/K}$-isotypical
component of the trivial connection in
$V$. By (i) one has
$H^\cdot_{dR}(V^{un})=
H^\cdot_{dR}(V)$.

There is a
unique
$K[[t]]$-sublattice $V^{un}_O\subset V^{un}$
preserved by
$\nabla (t\partial_t )$ and such that the
operator
$\kappa$ induced by
$\nabla (t\partial_t )$ on the fiber $V^{un}_0
:= V^{un}_O /tV^{un}_O$ is nilpotent. The de
Rham complex
$V^{un}\to
\omega
\otimes V^{un}$ is quasiisomorphic to its
logarithmic subcomplex\footnote{Here
$\omega_{O\log}=t^{-1}K[[t]]dt \subset
\omega$.}
$V^{un}_O\to\omega_{O\log}\otimes V^{un}_O$.
The latter projects
quasiisomorphically\footnote{By the
projections $V^{un}_O \to V^{un}_0$, $\Res
:\omega_O \to K$.} to the complex
$V^{un}_0 \buildrel\kappa\over\to V^{un}_0$
whose determinant is $\det V^{un}_0 /\det
V^{un}_0 =K$. This is our trivialization.

{\it Remark.}  If
$\nabla_{/K}$ comes from an absolute
connection then our trivialization is
horizontal with respect to the Gau\ss-Manin
connection.

\medskip

(b) Let $F'/F$ be a finite extension, $V_F
=(V_F ,
\nabla_{/K} )$ a relative connection on $F$.

(i)  $V_F$ is semi-simple if and only if
its pull-back $V_{F'}$ to $F'$ is
semi-simple.

Indeed, to see that the
pull-back of an irreducible $V_F$ is
semi-simple notice that the maximal semi-simple
subobject of
$V_{F'}$ descends to a subobject of $V_F$,
hence it equals $V_{F'}$. The rest follows
from (a)(iii).

(ii) If $V_F$ is irreducible and $V_{F'}$ is a
direct sum of several copies of the same rank
1 connection $L_{F'}$ then $V_F$ has rank 1.

Indeed, the class of $L_{F'}$ is $Gal
(F'/F)$-invariant, hence $L_{F'}$ comes from a
line with connection $L_F$ on
$F$.\footnote{Isomorphism classes of lines
with connection on $F'$ form the group $\omega
(F' )/d\log ({F'}^\times )=(\omega (F' )
/\omega (O'))/\Z $. Every Galois
invariant there lifts to the one in $\omega
(F')$, i.e., comes from $\omega (F)$.}
Replacing
$V_F$ by $V_F \otimes L_F^{-1}$ we are
reduced to the situation when $L_{F'}$ is a
trivial connection. Then $V_F$ has regular
singularities and all eigenvalues of
$\nabla_{/K} (t\partial_t )$ are in $\Q$.
Every irreducible connection with these
properties has rank 1.

\medskip

(c) If $(V_F ,\nabla_{/K})$ is irreducible then
one can find a finite extension $F'/F$ such
that the pull-back $V_{F'}$ of $V_F$ to $F'$ is
isomorphic to a direct sum of connections of
rank 1.

Indeed, we have the desired
decomposition over $F''=\bar{K}((t^{1/e}))$,
where
$\bar{K}$ is an algebraic closure of $K$, due
to the
Levelt-Turritin theorem (see e.g.~\cite{M}
III(1.2)) and (b)(i). It is defined over some
$F'=K'((t^{1/e}))$, $K'$ a finite
extension of
$K$, for the following reason. First, the rank
1 connections $L_{F''}$ that occur in the
Levelt-Turritin decomposition are defined over
some $F'$ as above. Since
  $F''\otimes \Hom (L_{F'},
V_{F'})\iso \Hom (L_{F''},V_{F''})$, one has
$\oplus\, L_{F'}\otimes\Hom (L_{F'},V_{F'})\iso
V_{F'}$.

\medskip

(d)  Now let $V_F$ be an $F$-vector space
equipped with an absolute integrable connection
$\nabla $; denote by $\nabla_{/K}$ the
corresponding relative connection. Let
$V_F =\oplus V_F^\alpha$ be a decomposition
which is
$\nabla_{/K}$-horizontal. Assume that for
every $V_F^\alpha \neq V_F^{\alpha'}$ there is
no non-trivial $\nabla_{/K}$-horizontal
morphisms
$V_F^\alpha \to V_F^{\alpha'}$. Then our
decomposition is
$\nabla$-horizontal.\footnote{Indeed, for every
horizontal vector field $\theta$ the $(\alpha
,\alpha' )$-component of
$\nabla (\theta )$ is a
$\nabla_{/K}$-horizontal morphism $V_F^\alpha
\to V_F^{\alpha'}.$}

Therefore for every $(V_F ,\nabla )$
the decomposition (5.9.1) by
$\nabla_{/K}$-isotypical components is
$\nabla$-horizontal. Thus
every indecomposable
$(V_F ,\nabla )$ is $\nabla_{/K}$-isotypical.

\subsection{A reduction to the rank 1 case.} We
use the notation from 5.1, so $F=K'((t))$, etc.

Let us show that computation of the
$\varepsilon$-connection for arbitrary $(V_F
,\nabla )$ can be reduced, in principle, to
the situation considered in 5.7, hence, by
(5.8.1), to the rank 1 situation.
Due to the multiplicativity property of
$\sE$ we can, and will, assume that
$(V_F ,\nabla )$ is indecomposable.
Furthemore, we assume that
$\nabla$ does not have regular
singularities (as a relative connection, see
5.3).

\medskip

{\bf Proposition.} For such $(V_F ,\nabla )$
there exists a  finite
extension $F'/F$ such that $(V_F ,\nabla )$ is
isomorphic to the push-forward of a
connection $(V'_{F'}, \nabla' )$ on
$\Spec F'$ of type considered in 5.7.

One has $\sE (F,V_F )_\nu =\sE
(F',V'_{F'})_\nu$ by (4.9.4); this is the
promised reduction.

\begin{proof}
Below, as in 5.9, we denote by $\nabla_{/K}$
the relative connection that corresponds to
$\nabla$ (the vertical part of $\nabla$).

By 5.9 (d)  $(V,\nabla_{_/K})$ is
$\nabla_{/K}$-isotypical. By 5.9 (c) applied
to the irreducible constituent of
$(V,\nabla_{/K})$ there is a finite Galois
covering
$F''/F$ such that every irreducible
subqoutient of $(V_{F''},\nabla_{/K})$ has
rank 1. Consider the
$\nabla_{/K}$-isotypical decomposition
$V_{F''}=\oplus V_{F''}^L$; its components
are labeled by irreducible
$L=(L_{F''},\nabla_{/K})$'s of rank 1.

The Galois group $Gal (F''/F)$ action on
$V_{F''}$ permutes components
$V_{F''}^L$. The action on the set of
components is transitive since $V_F$
is indecomposable.

Pick one of $L$'s that occur in
our decomposition; denote $V_{F''}^L$ by
$V'_{F''}$ and its connection by $\nabla'$. Let
$F'\subset F''$ be  the invariant field
of the stabilizer of $V_{F''}^L$
in $Gal (F''/F)$. We get an $F'$-vector
space
$V'_{F'}:=(V'_{F''})^{Gal (F''/F')}$ with
connection
$\nabla'$. Its push-forward to $F$ is
identified with
$(V_F ,\nabla )$  in the
obvious way. Thus $(V'_{F'},\nabla'
)$ is indecomposable, hence
$\nabla_{/K}$-isotypical.  By 5.9 (b)(ii)
applied to the irreducible constituent of
$(V'_{F'},\nabla'_{/K})$ our
$(V'_{F'},\nabla' )$ is of type considered in
  5.7. We are done.
\end{proof}

\bibliographystyle{plain}
\renewcommand\refname{References}

\end{document}